\documentclass[preprint,12pt]{elsarticle}

\usepackage[a4paper, total={6.6in, 9in}]{geometry}
\usepackage{lipsum}
\usepackage{amsfonts}
\usepackage{graphicx}
\usepackage{epstopdf}
\usepackage{tikz,graphicx}
\usetikzlibrary{3d,perspective} 
\usetikzlibrary{backgrounds,calc,shadows,positioning,fit}
\usepackage{amssymb}
\usepackage{amsmath}
\usepackage{enumerate} 
\usepackage{enumitem}
\usepackage[colorlinks=true,linkcolor=black,citecolor=blue,urlcolor=blue,]{hyperref}
\usepackage{algorithm}
\allowdisplaybreaks[4] 
\graphicspath{{Fig}}

\usepackage{amsthm}
\newtheorem{remark}{Remark}

\newtheorem{lemma}{Lemma}

\newtheorem{problem}{Problem}

\newtheorem{theorem}{Theorem}
\newtheorem{definition}{Definition}

\journal{Journal of Computational and Applied Mathematics}

\begin{document}

\begin{frontmatter}



\title{Mixed Finite Element Method for Multi-layer Elastic Contact Systems}


\author[Southeast]{Zhizhuo Zhang}
\ead{zhizhuo\_zhang@163.com}
\author[UPVD]{Mikaël Barboteu}
\ead{barboteu@univ-perp.fr}
\author[Southeast]{Xiaobing Nie}
\ead{xbnie@seu.edu.cn}
\author[UPVD]{Serge Dumont}
\ead{serge.dumont@unimes.fr}
\author[Sohag,Ahlia]{Mahmoud Abdel-Aty}
\ead{mabdelaty@zewailcity.edu.eg}
\author[Southeast]{Jinde Cao\corref{cao}}
\cortext[cao]{corresponding author: Jinde Cao (Email addresses: jdcao@seu.edu.cn)}

\affiliation[Southeast]{organization={School of Mathematics, Southeast University},
            city={Nanjing},
            postcode={211189}, 
            country={China}}
\affiliation[UPVD]{organization={Laboratoire de Mathématiques et Physique, Université de Perpignan Via Domitia},
            addressline={52 Avenue Paul Alduy}, 
            city={Perpignan},
            postcode={66860}, 
            country={France}}
\affiliation[Sohag]{organization={Department of Mathematics, Faculty of Science, Sohag University},
            city={Sohag},
            postcode={82524}, 
            country={Egypt}}
\affiliation[Ahlia]{organization={Deanship of Graduate Studies and Scientific Research, Ahlia University},
            city={Manama},
            postcode={10878}, 
            country={Bahrain}}

\begin{abstract}
With the development of multi-layer elastic systems in the field of engineering mechanics, the corresponding variational inequality theory and algorithm design have received more attention and research. In this study, a class of equivalent saddle point problems with interlayer Tresca friction conditions and the mixed finite element method are proposed and analyzed. Then, the convergence of the numerical solution of the mixed finite element method is theoretically proven, and the corresponding algebraic dual algorithm is provided. 
Finally, through numerical experiments, the mixed finite element method is not only compared with the layer decomposition method, but also its convergence relationship with respect to the spatial discretization parameter $H$ is verified.
\end{abstract}




\end{frontmatter}



\section{Introduction}

With the application of multi-layer physical models in engineering mechanics, multi-layer elastic contact systems and their numerical algorithms based on variational inequality theory have gradually received attention and research \cite{ma2021analytical, ma2022dynamic}.
Non-smooth contact boundary conditions have always posed a challenge and have been a focal point in algorithm design for contact problems \cite{haslinger2014domain,han2020numerical}. In multi-layer contact problems, as the number of contact surfaces increases, the coupling effect between displacement fields presents additional challenges for the design of finite element algorithms.
In the unilateral contact problem, there are currently many excellent numerical algorithms used to solve the corresponding variational inequalities, including the Lagrange method \cite{simo1992augmented, hild2010stabilized}, the non-smooth multi-scale method \cite{krause2009nonsmooth}, the primitive dual active set algorithm \cite{hueber2008primal, abide2022unified}, Nitsche method \cite{chouly2022nitsche} and others.
However, methods that can be generalized to solve the problem of contact between two objects are relatively scarce, and algorithms that can handle multi-layer contact problems are even less studied \cite{haslinger1982approximation, haslinger2014domain, bayada2008convergence}.
It is worth emphasizing that multi-layer contact problems have a wide range of applications in engineering practice, such as analyzing the mechanical response of pavement under vehicle load and studying the mechanical properties of multi-layer composite materials \cite{hu2015three, ma2023probability}. In terms of theoretical research, addressing the nonlinear problem arising from the coupling effect of the boundary displacement field under the contact function in a multi-layer contact system is also deserving of in-depth discussion.

In previous research, the existence and uniqueness of theoretical solutions for multi-layer elastic and viscoelastic contact systems have been established within the framework of variational inequalities. Additionally, the convergence of finite element numerical solutions has been demonstrated \cite{zhang2022variational, zhang2024variational}.
On the basis of this theory, a Layer Decomposition (LD) algorithm is proposed, which is based on the domain decomposition method, to calculate the numerical solution of the multi-layer elastic system \cite{zhang2024layer}. 
Different from the design concept of decomposing and calculating the two-object contact problem in the domain decomposition method \cite{haslinger2014domain}, the mixed finite element method deforms the original problem equivalently and directly solves the global problem \cite{haslinger1982approximation}.
Therefore, in this study, the Mixed Finite Element (MFE) method is extended to solve the multi-layer elastic contact system with interlayer Tresca friction conditions. 
The proposed MFE method can not only be compared and verified with the LD method but also provide more theoretical support and algorithmic tools for mechanical modeling of multi-layer physical models in engineering mechanics.

Specifically, the paper is structured as follows:

In Section 2, the partial differential equations and variational inequalities of the multi-layer elastic contact system will be presented. Subsequently, equivalent saddle point problems and mixed problems will be proposed and validated.
In Section 3, the convergence analysis and error order of the finite element numerical solution for this mixed problem will be presented and proven.
In Section 4, the implementation details of the MEF method will be given, including the algebraic form and algebraic dual form of the mixed problem, which will provide support for the design of specific algorithms.
In Section 5, a numerical simulation experiment was conducted on a three-dimensional three-layer elastic contact system. The experiment compared the LD method with the MFE method and also verified the convergence properties of the MFE method. Additionally, a numerical experiment was conducted on a two-dimensional four-layer contact system to test the convergence rate of the numerical solution.
Finally, in Section 6, all studies in the paper are summarized.

\section{Setting of the problem}

\subsection{The physical model and variational inequality}

\begin{figure}[!t]
\centering
\begin{minipage}{0.49\linewidth}
\begin{tikzpicture}[3d view={-30}{30}, scale=0.7]
    \draw[->] (0,0,0) -- (0.5,0,0) node[pos=1,right]{x};
    \draw[->] (0,0,0) -- (0,0.5,0) node[pos=1,left]{y};
    \draw[->] (0,0,0) -- (0,0,0.5) node[pos=1,above]{z};
    \foreach \n in {0,1} 
    {
        \draw (0,0,0+\n) -- (8,0,0+\n) -- (8,4,0+\n) -- (0,4,0+\n) -- (0,0,0+\n);
        \draw (0,0,1+\n) -- (8,0,1+\n) -- (8,4,1+\n) -- (0,4,1+\n) -- (0,0,1+\n);
        \draw (0,0,0+\n) -- (0,0,1+\n);
        \draw (8,0,0+\n) -- (8,0,1+\n);
        \draw (8,4,0+\n) -- (8,4,1+\n);
        \draw (0,4,0+\n) -- (0,4,1+\n);
    }
    \foreach \n in {0,1} 
    {
        \draw (0,0,3+\n) -- (8,0,3+\n) -- (8,4,3+\n) -- (0,4,3+\n) -- (0,0,3+\n);
        \draw (0,0,4+\n) -- (8,0,4+\n) -- (8,4,4+\n) -- (0,4,4+\n) -- (0,0,4+\n);
        \draw (0,0,3+\n) -- (0,0,4+\n);
        \draw (8,0,3+\n) -- (8,0,4+\n);
        \draw (8,4,3+\n) -- (8,4,4+\n);
        \draw (0,4,3+\n) -- (0,4,4+\n);
    }
    \foreach \n in {0.1,0.2,...,0.9}
    { \fill (0,0,2+\n) circle (0.5pt);
      \fill (8,0,2+\n) circle (0.5pt);
      \fill (8,4,2+\n) circle (0.5pt);
      \fill (0,4,2+\n) circle (0.5pt);
    }
    \begin{scope}[canvas is xy plane at z=5]
        \node at (4,2) [color = blue, transform shape,scale=1.5] {$\Gamma_2^1$};
    \end{scope}
    \filldraw[fill opacity=0.1,fill=blue] (0,0,5) -- (8,0,5) -- (8,4,5) -- (0,4,5) -- cycle ;
    \draw [color=blue] (0,0,5) -- (8,0,5) -- (8,4,5) -- (0,4,5) -- (0,0,5);
    \node at (8,0,4.5) [color = blue,scale=1.5,right] {$\Omega^{1}$};
    
    \begin{scope}[canvas is xy plane at z=4]
        \node at (4,2) [color = purple, transform shape,scale=1.5] {$\Gamma_2^2/\Gamma_3^1$};
    \end{scope}
    \filldraw[fill opacity=0.1,fill=purple] (0,0,4) -- (8,0,4) -- (8,4,4) -- (0,4,4) -- cycle ;
    \draw [color=purple] (0,0,4) -- (8,0,4) -- (8,4,4) -- (0,4,4) -- (0,0,4);
    \node at (8,0,3.5) [color = purple, scale=1.5,right] {$\Omega^{2}$};
    
    \begin{scope}[canvas is xy plane at z=3]
        \node at (4,2) [color=orange, transform shape,scale=1.5] {$\Gamma_2^3/\Gamma_3^2$};
    \end{scope}
    \filldraw[fill opacity=0.1,fill=orange] (0,0,3) -- (8,0,3) -- (8,4,3) -- (0,4,3) -- cycle ;
    \draw [color=orange] (0,0,3) -- (8,0,3) -- (8,4,3) -- (0,4,3) -- (0,0,3);
    
    \begin{scope}[canvas is xy plane at z=2]
        \node at (4,2) [color=cyan, transform shape,scale=1.5] {$\Gamma_2^{n-1}/\Gamma_3^{n-2}$};
    \end{scope}
    \filldraw[fill opacity=0.1,fill=cyan] (0,0,2) -- (8,0,2) -- (8,4,2) -- (0,4,2) -- cycle ;
    \draw [color=cyan] (0,0,2) -- (8,0,2) -- (8,4,2) -- (0,4,2) -- (0,0,2);
    \node at (8,0,1.5) [color = cyan, scale=1.5,right] {$\Omega^{n-1}$};
    
    \begin{scope}[canvas is xy plane at z=1]
        \node at (4,2) [color = gray, transform shape,scale=1.5] {$\Gamma_2^{n}/\Gamma_3^{n-1}$};
    \end{scope}
    \filldraw[fill opacity=0.1,fill=gray] (0,0,1) -- (8,0,1) -- (8,4,1) -- (0,4,1) -- cycle ;
    \draw [color=gray] (0,0,1) -- (8,0,1) -- (8,4,1) -- (0,4,1) -- (0,0,1);
    \node at (8,0,0.5) [color = gray, scale=1.5,right] {$\Omega^{n}$};
    
    \begin{scope}[canvas is xy plane at z=0]
        \node at (4,2) [color = teal, transform shape,scale=1.5] {$\Gamma_1^{n}(\Gamma_3^{n})$};
    \end{scope}
    \filldraw[fill opacity=0.1,fill=teal] (0,0,0) -- (8,0,0) -- (8,4,0) -- (0,4,0) -- cycle ;
    \draw [color=teal] (0,0,0) -- (8,0,0) -- (8,4,0) -- (0,4,0) -- (0,0,0);
    
    \begin{scope}[canvas is yz plane at x=0]
        \node at (1,0.5) [rotate=180,yscale=-1,transform shape,scale=1.2] {$\Gamma_1^{n}$};
        \node at (1,1.5) [rotate=180,yscale=-1,transform shape,scale=1.2] {$\Gamma_1^{n-1}$};
        \node at (1,3.5) [rotate=180,yscale=-1,transform shape,scale=1.2] {$\Gamma_1^{2}$};
        \node at (1,4.5) [rotate=180,yscale=-1,transform shape,scale=1.2] {$\Gamma_1^{1}$};
    \end{scope}
    \begin{scope}[canvas is xz plane at y=0]
        \node at (7,0.5) [xscale=1,transform shape,scale=1.2] {$\Gamma_1^{n}$};
        \node at (7,1.5) [xscale=1,transform shape,scale=1.2] {$\Gamma_1^{n-1}$};
        \node at (7,3.5) [xscale=1,transform shape,scale=1.2] {$\Gamma_1^{2}$};
        \node at (7,4.5) [xscale=1,transform shape,scale=1.2] {$\Gamma_1^{1}$};
    \end{scope}
    \begin{scope}[canvas is yz plane at x=8]
        \node at (3,0.5) [rotate=180,yscale=-1,transform shape,scale=1.2] {$\Gamma_1^{n}$};
        \node at (3,1.5) [rotate=180,yscale=-1,transform shape,scale=1.2] {$\Gamma_1^{n-1}$};
        \node at (3,3.5) [rotate=180,yscale=-1,transform shape,scale=1.2] {$\Gamma_1^{2}$};
        \node at (3,4.5) [rotate=180,yscale=-1,transform shape,scale=1.2] {$\Gamma_1^{1}$};
    \end{scope}
    \begin{scope}[canvas is xz plane at y=4]
        \node at (1,0.5) [xscale=1,transform shape,scale=1.2] {$\Gamma_1^{n}$};
        \node at (1,1.5) [xscale=1,transform shape,scale=1.2] {$\Gamma_1^{n-1}$};
        \node at (1,3.5) [xscale=1,transform shape,scale=1.2] {$\Gamma_1^{2}$};
        \node at (1,4.5) [xscale=1,transform shape,scale=1.2] {$\Gamma_1^{1}$};
    \end{scope}
    \foreach \x in {3.5,3.83,4.16,4.5} \foreach \y in {1.5,1.83,2.16,2.5} {
    \draw[->,color=red, line width=0.2pt] (0.15+\x,\y,5.4)--(\x,\y,5);}
    \node at (4,2,6.3) [color=red,scale=1.6] {$f_1$};
\end{tikzpicture}
\caption{Three-dimensional n-layer contact system}
\label{4:fig1}
\end{minipage}
\begin{minipage}{0.49\linewidth}
\centering
\includegraphics[width=0.95\textwidth]{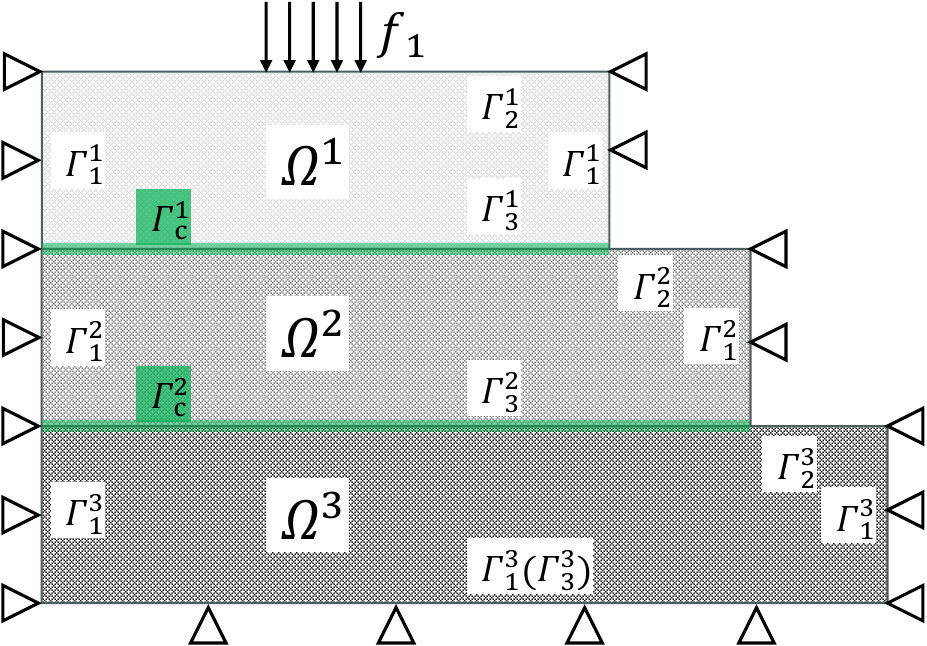}
\caption{Two-dimensional 3-layer contact system}
\label{4:fig:model.2D}
\end{minipage}
\end{figure}

Before presenting the target problem, it is necessary to introduce the specific physical model of the problem and some function spaces.
Referring to the physical models in pavement mechanics studies, the three-dimensional multi-layer contact system is constructed as shown in Fig.\ref{4:fig1}.
\textcolor{black}{Furthermore, to more clearly illustrate the layered structure, a two-dimensional, three-layer contact system is depicted in Fig.\ref{4:fig:model.2D}.}

In this system, the bounded open domain occupied by each layer of an elastic body is denoted as $\Omega^{i}\subset \mathbb{R}^{d}$, where $d=2,3$ and $i=1,2,\ldots,n$, with $n$ representing the number of elastic bodies in the system. Therefore, $\Omega=\cup_{i=1}^{n} \Omega^{i}$ is used to represent the domain occupied by the system. The Lipschitz boundary of $\Omega^{i}$ is denoted by $\partial\Omega^{i} = \Gamma^{i} = \bar{\Gamma}^{i}_{1}\cup\bar{\Gamma}^{i}_{2}\cup\bar{\Gamma}^{i}_{3}$, where $\bar{\Gamma}^{i}_{1}$, $\bar{\Gamma}^{i}_{2}$ and $\bar{\Gamma}^{i}_{3}$ are disjoint, non-empty, and open in $\Gamma^{i}$ as shown in Fig.\ref{4:fig1}.
For $x\in\Gamma$, the unit outward normal vector at $x$ is denoted by $\boldsymbol{\nu}$. To distinguish, $\boldsymbol{\alpha}^{i}$ and $\boldsymbol{\beta}^{i}$ represent the unit outward normal vectors at $x\in\Gamma^{i}_{2}$ and $x\in\Gamma^{i}_{3}$, respectively.
\textcolor{black}{In addition, the space of second-order symmetric tensors defined on $\mathbb{R}^d$ is denoted by $\mathbb{S}^d$, and the inner product (double dot product) on this space is denoted by "$:$".}

Obviously, the physical model pertains to the contact problem model in solid mechanics, where the contact zone with friction between $\Omega^{i}$ and $\Omega^{i+1}$ is represented by $\Gamma_{c}^{i} = \Gamma^{i}_{3} \cap \Gamma^{i+1}_{2}$, $i=1,\ldots,n-1$. Subsequently, the boundary conditions of the multi-layer elastic system are established as follows: on boundary $\Gamma_{1}^{i}$, the displacement field is zero; on boundary $\Gamma_{2}^{1}$, the elastic body $\Omega^{1}$ is subjected to a surface force $\boldsymbol{f}_{1}$. 
\textcolor{black}{On boundaries $\Gamma^{i}_{3}$ and $\Gamma^{i+1}_{2}$ ($i=1,\ldots,n-1$), regions outside of the contact zone $\Gamma_{c}^{i}$ are unconstrained by boundary conditions. Consequently, in subsequent discussions, we continue to denote $\Gamma^{i}_{3}\cap\Gamma_{c}^{i}$ and $\Gamma^{i+1}_{2}\cap\Gamma_{c}^{i}$ simply as $\Gamma^{i}_{3}$ and $\Gamma^{i+1}_{2}$, respectively.}
Moreover, the system is affected by a force field, and the body force on the elastic body $\Omega^{i}$ is denoted as $\boldsymbol{f}^{i}_{0}$.

It is our goal to solve the displacement field for the multi-layer elastic system with interlayer friction coupling under mechanical action. Therefore, the vector function space and symmetric second-order tensor function space defined on $\Omega^{i}$ are denoted as:
$$
\begin{aligned}
& V^i=\left\{\boldsymbol{v}=\left(v_k\right) \in\left(\mathcal{H}^1\left(\Omega^i\right)\right)^d \big| \boldsymbol{v}=\mathbf{0} \text { on } \Gamma_1^i\right\}, \\
& Q^i=\left\{\boldsymbol{\tau}=\left(\tau_{k l}\right) \in\left(L^2\left(\Omega^i\right)\right)^{d \times d} \big| \tau_{l k}=\tau_{l k}, 1 \leqslant k, l \leqslant d\right\}, \\
& Q_1^i=\left\{\boldsymbol{\tau} \in Q^i \big| \nabla \cdot \boldsymbol{\tau} \in\left(L^2\left(\Omega^i\right)\right)^d\right\},
\end{aligned}
$$
where $i=1,\ldots,n$, $\mathcal{H}^k(\cdot) = W^{k,2} (\cdot)$ is Sobolev space and $\nabla$ is gradient operator. 
In particular, let $\mathcal{H}^0(\cdot) = L^{2} (\cdot)$. Then, the vector field and second-order tensor field of the total system can be denoted as $\boldsymbol{v} = \left(\boldsymbol{v}^1, \boldsymbol{v}^2, \ldots, \boldsymbol{v}^n\right)$ and $\boldsymbol{\tau}=\left(\boldsymbol{\tau}^1, \boldsymbol{\tau}^2, \ldots, \boldsymbol{\tau}^n\right)$, respectively. Hence, the spaces of $\boldsymbol{v}$ and $\boldsymbol{\tau}$ are defined as:
$$
V=V^1 \times V^2 \times \cdots \times V^n \text { and } Q_1=Q_1^1 \times Q_1^2 \times \cdots \times Q_1^n.
$$
Since the above spaces are Hilbert spaces, the inner product on them can be defined as follows:
\begin{align*}
& \left(\boldsymbol{u}^i, \boldsymbol{v}^i\right)_{s}=\sum_{k=1}^d\left(u_k^i, v_k^i\right)_{\mathcal{H}^{s}\left(\Omega^i\right)}, && (\boldsymbol{u}, \boldsymbol{v})_{s}=\sum_{i=1}^n \left(\boldsymbol{u}^i, \boldsymbol{v}^i\right)_{s}, \\
& \left(\boldsymbol{\sigma}^i, \boldsymbol{\tau}^i\right)_{s}=\sum_{k, l=1}^d\left(\sigma_{k l}^i, \tau_{k l}^i\right)_{\mathcal{H}^{s}\left(\Omega^i\right)}, && (\boldsymbol{\sigma}, \boldsymbol{\tau})_{s}=\sum_{i=1}^n \left(\boldsymbol{\sigma}^i, \boldsymbol{\tau}^i\right)_{s},
\end{align*}
where $s=0$ or $1$.
Based on the above definitions of inner product, the standard norm $\|\cdot\|_{s}$ can be defined for \textcolor{black}{these spaces}.
Then the strain tensor $\boldsymbol{\varepsilon}\left(\boldsymbol{v}^i\right) \in Q^i$ and stress tensor $\boldsymbol{\sigma}^{i} \in Q^i_{1}$ are defined as: 
$$
\boldsymbol{\varepsilon}\left(\boldsymbol{v}^i\right)=\frac{1}{2}\left(\nabla \boldsymbol{v}^i+\left(\nabla \boldsymbol{v}^i\right)^{\top}\right) \text{ and } \boldsymbol{\sigma}^{i} = A^{i} : \boldsymbol{\varepsilon}\left(\boldsymbol{v}^i\right),
$$
where $A^{i}$ is a symmetric fourth-order tensor, \textcolor{black}{which satisfies} the following properties:
\begin{equation}\label{4:character:A}
\left\{\begin{aligned}
&\text { (a) } {A}^{i}: L^{\infty}\left(\Omega^{i}\right)^{d\times d\times d\times d},~i=1,\ldots,n; \\ 
&\text { (b) } {A}^{i} \text{ is the symmetric fourth-order tensor, that is } A^{i} = \left\{a^{i}_{klrs}\right\} \text{ and }\\ 
&~~~~~~~~  a^{i}_{klrs} = a^{i}_{lkrs} = a^{i}_{rskl}, ~k,l,r,s = 1,\ldots,d;\\ 
&\text{ (c) \textcolor{black}{There exist} positive constants } a^{i}_{m}>0 \text{ and } a^{i}_{M}>0 \text { such that } \\ 
&~~~~~~~~  a^{i}_{m}\boldsymbol{\tau}^{i}:\boldsymbol{\tau}^{i} \leqslant A^{i}:\boldsymbol{\tau}^{i}:\boldsymbol{\tau}^{i} \leqslant a^{i}_{M}\boldsymbol{\tau}^{i}:\boldsymbol{\tau}^{i}.
\end{aligned}\right.
\end{equation}

Since $meas_{2}(\Gamma^{i}_{1})>0$, the following Korn's inequality holds:
\begin{equation}\label{4:eq:korn}
\left\|\boldsymbol{v}^i\right\|_{1} \leqslant c_k^i\left\|\boldsymbol{\varepsilon}\left(\boldsymbol{v}^i\right)\right\|_{0} \quad \forall \boldsymbol{v}^i \in V^i.
\end{equation}
where the constant $c^{i}_{k}>0$ depends on $\Omega^{i}$ and $\Gamma^{i}_{1}$. Therefore, another inner product on $V^{i}$ can be defined as: 
$$
\left(\boldsymbol{u}^i, \boldsymbol{v}^i\right)_{\textcolor{black}{V^{i}}}=\left(\boldsymbol{\varepsilon}\left(\boldsymbol{u}^i\right), \boldsymbol{\varepsilon}\left(\boldsymbol{v}^i\right)\right)_{0} \quad \forall \boldsymbol{u}^i, \boldsymbol{v}^i \in V^i.
$$
Furthermore, the norm $\|\cdot\|_{\textcolor{black}{V^{i}}}$ induced by this inner product is equivalent to $\|\cdot\|_{1}$ on $V^{i}$. Similarly, another inner product on $V$ can be defined as
$$
(\boldsymbol{u}, \boldsymbol{v})_V=\sum_{i=1}^n \left(\boldsymbol{u}^i, \boldsymbol{v}^i\right)_{\textcolor{black}{V^{i}}} \quad \forall \boldsymbol{u}, \boldsymbol{v} \in V,
$$
and the norm $\|\cdot\|_{V}$ is equivalent to $\|\cdot\|_{1}$ on $V$.

It is worth noting that the Coulomb friction contact condition is considered to be close to the real contact state. However, its theoretical analysis and algorithm design are very complex. Therefore, approximating it through the simpler Tresca friction contact condition is a universal and effective solution. Existing research has shown that by designing and solving a fixed-point algorithm that calculates the unilateral contact problem with Tresca's friction law at each step, the unilateral contact problem with Coulomb's friction law can be approximately solved \cite{laborde2008fixed, eck1998existence}.
Based on this, the research will focus on the multi-layer elastic contact system with \textcolor{black}{the Tresca friction} contact condition. The partial differential equations problem under the action of external forces can be expressed as follows:
\begin{problem}[$\mathcal{P}_{0}$]\label{prb:4:p_0}
Find a displacement field $\boldsymbol{u}^{i}:\Omega^{i} \rightarrow \mathbb{R}^{d}$ and the stress field $\boldsymbol{\sigma}^{i}: \Omega^{i} \rightarrow \mathbb{S}^{d}$ ($i=1,2,\ldots,n$) such that:
\begin{align*}
&\operatorname{Div} \boldsymbol{\sigma}^{i}+\boldsymbol{f}_{0}^{i}=\mathbf{0} && \text { in } \Omega^{i},\\
&\boldsymbol{u}^{i}=\mathbf{0} && \text { on } \Gamma_{1}^{i}, \label{3:pde.3}\\
&\boldsymbol{\sigma}^{1} \cdot \boldsymbol{\alpha}^{1}=\boldsymbol{f}_{1} && \text { on } \Gamma_{2}^{1}, \\
&\sigma^{i}_{\alpha} = -\sigma^{i-1}_{\beta},~
\boldsymbol{\sigma}_{\tau}^{i} = \boldsymbol{\sigma}_{\eta}^{i-1},
&&\text { on } \Gamma_{2}^{i}, ~i\ne 1,\\
&\left.
\begin{aligned}
& \left[{u}_{N}^{i}\right] \leqslant 0,~ \sigma^{i}_{N}\leqslant 0,~ \sigma^{i}_{N}\cdot[u_{N}^{i}] = 0,~ |\boldsymbol{\sigma}_{T}^{i}|\leqslant g^{i}(x)\\
& |\boldsymbol{\sigma}_{T}^{i}| < g^{i}(x) \Rightarrow
[\boldsymbol{u}_{T}^{i}] = 0\\
& |\boldsymbol{\sigma}_{T}^{i}| = g^{i}(x) \Rightarrow
\boldsymbol{\sigma}_{T}^{i} = -\lambda[\boldsymbol{u}_{T}^{i}], ~\lambda\geqslant 0
\end{aligned}
\right\} &&  \text { on } \Gamma_{c}^{i},~i\ne n,
\end{align*}
where $g^i({x})$ is the interlayer friction bounded function and $[\cdot]$ is the jump operator on $\Gamma^{i}_{c}$.
\end{problem}

Before characterizing the displacement function on the contact boundary zone $\Gamma^{i}_{c}$, the following trace operators are introduced:
\begin{definition}[]\label{4:def:1}
Let $\gamma^{i}_{r}: V^{i}\to L^{2}(\Gamma^{i}_{r})^{d}$ ($r=2$ or $3$) be the trace operator, and 
$$
\mathcal{H}^{1/2}\left(\Gamma_r^i\right)^{d}=\gamma^i_r\left(V^i\right),
$$
where $\mathcal{H}^{1/2}\left(\Gamma_r^i\right)^{d}$ is a Hilbert space with norm:
$$
\|\boldsymbol{\varphi}^{i}\|_{1/2, \Gamma_r^i}=\inf_{\substack{\boldsymbol{v}^i \in V^{i}\\ \gamma^i_r \boldsymbol{v}^i=\boldsymbol{\varphi}^{i}}}\left\|\boldsymbol{v}^i\right\|_{V^i}, ~~\boldsymbol{\varphi}^{i}\in \mathcal{H}^{1/2}\left(\Gamma^{i}_r\right).
$$
\textcolor{black}{Since $Q^{i}_{1} = \mathcal{H}^{1}(\Omega^{i})^{d\times d}\cap Q^{i}$, the trace operator on it can also be denoted by $\gamma^{i}_{r}: Q^{i}_{1} \to L^{2}(\Gamma^{i}_{r})^{d\times d}$, and it possesses similar properties as mentioned above.}
\end{definition}
Then, according to the trace theorem on Sobolev space \cite{adams2003sobolev}, the following inequality holds:
\begin{equation}\label{4:eq:trace}
\left\|\gamma_r^i \boldsymbol{v}^i\right\|_{1/2, \Gamma_r^j} \leqslant c^{i}_{t}\left\|\boldsymbol{v}^i\right\|_{V^i}, \forall \boldsymbol{v}^i \in V^i, r=2 \text { or } 3,
\end{equation}
where the constant $c^{i}_{t}$ dependent on $\Omega^{i}$ and $\Gamma^{i}$.

In order to address the non-penetration condition and friction condition at the contact boundary zone, the normal and tangential components of the boundary function need to be defined as follows:
\begin{align*}
&v_{\beta}^{i}=\gamma^{i}_{3}\boldsymbol{v}^{i} \cdot \boldsymbol{\beta}^{i}, && \boldsymbol{v}_{\eta}^{i}= \gamma^{i}_{3}\boldsymbol{v}^{i} - v_{\beta}^{i}\cdot\boldsymbol{\beta}^{i}, && \textcolor{black}{\boldsymbol{v}^{i}\in V^{i}};\\
&v_{\alpha}^{i}=\gamma^{i}_{2}\boldsymbol{v}^{i} \cdot \boldsymbol{\alpha}^{i}, && \boldsymbol{v}_{\tau}^{i}= \gamma^{i}_{2}\boldsymbol{v}^{i} -v_{\alpha}^{i}\cdot\boldsymbol{\alpha}^{i}, && \textcolor{black}{\boldsymbol{v}^{i}\in V^{i}};\\
&\sigma_{\beta}^{i}=\boldsymbol{\beta}^{i}\cdot \gamma^{i}_{3}\boldsymbol{\sigma}^{i} \cdot \boldsymbol{\beta}^{i}, && \boldsymbol{\sigma}_{\eta}^{i}= \gamma^{i}_{3}\boldsymbol{\sigma}^{i} \cdot \boldsymbol{\beta}^{i}-\sigma_{\beta}^{i}\cdot\boldsymbol{\beta}^{i}, && \textcolor{black}{\boldsymbol{\sigma}^{i}\in Q^{i}_1};\\
&\sigma_{\alpha}^{i}=\boldsymbol{\alpha}^{i}\cdot \gamma^{i}_{2}\boldsymbol{\sigma}^{i} \cdot \boldsymbol{\alpha}^{i}, && \boldsymbol{\sigma}_{\tau}^{i}= \gamma^{i}_{2}\boldsymbol{\sigma}^{i} \cdot \boldsymbol{\alpha}^{i}-\sigma_{\alpha}^{i}\cdot\boldsymbol{\alpha}^{i}, && \textcolor{black}{\boldsymbol{\sigma}^{i}\in Q^{i}_1}.
\end{align*}
Then, the jump operator $[\cdot]$ on $\Gamma^{i}_{c}$ can be defined as
$$
\left[v_N^i\right]=v_{\beta}^{i} + v_{\alpha}^{i+1},~~
\left[\boldsymbol{v}_T^i\right]=\boldsymbol{v}_\eta^i - \boldsymbol{v}_\tau^{i+1}.
$$
In a physical model, the two elastic bodies cannot penetrate each other in the contact zone, so the displacement field of the system satisfies $\left[v_N^i\right]\leqslant 0$, $i = 1,\ldots,n-1$. To characterize these contact conditions, let $\delta^{i}_{r}: V^{i} \to L^{2}\left( \Gamma^{i}_{r}\right) \times L^{2}\left( \Gamma^{i}_{r}\right)^{d}$, and $\delta^{i}_{r} \boldsymbol{v}^{i} = \left(v^{i}_{n}, \boldsymbol{v}^{i}_{t}\right)$, where
\begin{align*}
& v^{i}_{n} = v^{i}_{\alpha} \text{ and } \boldsymbol{v}^{i}_{t} = \boldsymbol{v}^{i}_{\tau}, \text{ on } \Gamma^{i}_{2},\\
& v^{i}_{n} = v^{i}_{\beta} \text{ and } \boldsymbol{v}^{i}_{t} = \boldsymbol{v}^{i}_{\eta}, \text{ on } \Gamma^{i}_{3}.
\end{align*}
And let
$$
\delta^{i}_{r} V^{i} = W^{i}_{r} = W^{i}_{rN} \times  W^{i}_{rT} \text{ on } \Gamma^{i}_{r}.
$$
It is worth explaining that $\boldsymbol{v}^{i}_{t}$ is the projection of $\boldsymbol{v}^{i}$ onto the tangent plane of $\Gamma^{i}_{r}$, so $W^{i}_{rT}$ and $\mathcal{H}^{1/2}\left(\Gamma_r^i\right)^{d-1}$ are isomorphic. Therefore, it is not difficult to find that $W^{i}_{r}$ and $\mathcal{H}^{1/2}\left(\Gamma_r^i\right)^{d}$ are isomorphic. Indeed, let $\boldsymbol{\varphi}^{i} \in \mathcal{H}^{1/2}\left(\Gamma_r^i\right)^{d}$, then
$$
\boldsymbol{\omega}^{i} = \left(\varphi^{i}_{n}, \boldsymbol{\varphi}^{i}_{t}\right) = \varrho^{i}_{r} \boldsymbol{\varphi}^{i} \in W^{i}_{r}, 
$$
where operator $\varrho^{i}_{r}: \mathcal{H}^{1/2}\left(\Gamma_r^i\right)^{d} \to W^{i}_{r}$ is inversable. Hence, the norm in $W^{i}_{r}$ can be defined as
$$
\|\boldsymbol{\omega}^{i}\|_{W^{i}_{r}} = \|(\varrho^{i}_{r})^{-1}\boldsymbol{\omega}^{i}\|_{1/2,\Gamma^{i}_{r}},~~\boldsymbol{\omega}^{i}\in W^{i}_{r}.
$$
Moreover, let $W_{r} = W^{1}_{r} \times \cdots W^{n}_{r}$ with the norm 
$$
\|\boldsymbol{\omega}_{r}\|^{2}_{W_{r}} = \sum_{i=1}^{n} \|\boldsymbol{\omega}^{i}_{r}\|^{2}_{W^{i}_{r}}.
$$

Let $\mathcal{H}^{-1/2}\left(\Gamma_r^i\right)$, $W^{i\prime}_{rN}$ and $W^{i\prime}_{rT}$ be the dual spaces to $\mathcal{H}^{1/2}\left(\Gamma_r^i\right)$, $W^{i}_{rN}$ and $W^{i}_{rT}$ ($r=2$ or $3$). Then, the dual space of $W^{i}_{r}$ can be denoted by $W^{i\prime}_{r} = W^{i\prime}_{rN} + W^{i\prime}_{rT}$. For $\boldsymbol{\omega}^{i*} = \left(\omega^{i*}_{n}, \boldsymbol{\omega}^{i*}_{t}\right) \in W^{i\prime}_{r}$, $\boldsymbol{\omega}^{i} = \left(\omega^{i}_{n}, \boldsymbol{\omega}^{i}_{t}\right) \in W^{i}_{r}$, we write
$$
\left\langle\boldsymbol{\omega}^{i*}, \boldsymbol{\omega}^{i} \right\rangle_{r} = \left\langle{\omega}^{i*}_{n}, {\omega}^{i}_{n} \right\rangle_{r} + \left\langle\boldsymbol{\omega}^{i*}_{t}, \boldsymbol{\omega}^{i}_{t} \right\rangle_{r}.
$$
But, it should be noted that $Ker(W^{i}_{rT}) = \left\{ \boldsymbol{\omega}^{i*}_{t} \in W^{i\prime}_{rT} \big| \left\langle\boldsymbol{\omega}^{i*}_{t}, \boldsymbol{\omega}^{i}_{t}\right\rangle _{r}=0, ~\forall \boldsymbol{\omega}^{i}_{t}\in W^{i}_{rT},~\boldsymbol{\omega}^{i*}_{t}\neq \boldsymbol{0} \right\} \neq \emptyset$. In the subsequent discussion, the dual space of $W^{i}_{rT}$ defaults to $W^{i\prime}_{rT}-Ker(W^{i}_{rT})$, and is still denoted as $W^{i\prime}_{rT}$ for the sake of notation simplicity.
Then, $\|\boldsymbol{\omega}^{i*}\|_{W^{i\prime}_{r}}$ is defined as the usual dual norm of $\boldsymbol{\omega}^{i*}$, that is
$$
\|\boldsymbol{\omega}^{i*}\|_{W^{i\prime}_{r}} = \sup_{\boldsymbol{\omega}^{i}\in W^{i}_{r}} \frac{\left\langle\boldsymbol{\omega}^{i*}, \boldsymbol{\omega}^{i} \right\rangle_{r}}{\|\boldsymbol{\omega}^{i}\|_{W^{i}_{r}}}.
$$
Similarly, let $W^{\prime}_{r} = W^{1\prime}_{r} \times \cdots W^{n\prime}_{r}$ with the norm $\|\boldsymbol{\omega}^{*}_{r}\|^{2}_{W^{\prime}_{r}} = \sum_{i=1}^{n} \|\boldsymbol{\omega}^{i*}_{r}\|^{2}_{W^{i\prime}_{r}}$.
\textcolor{black}{Then, for given $\boldsymbol{v}^{i} \in V^{i}$ and $\boldsymbol{\tau}^{i} \in Q^{i}_{1}$, a Green's formula of the following form can be established \cite{haslinger1982approximation,aubin2007approximation}:
\begin{equation}\label{4:eq:Green}
\begin{aligned}
&\left(\boldsymbol{\tau}^{i}, \boldsymbol{\varepsilon}(\boldsymbol{v^{i}})\right)_{0} + \left(\nabla\cdot\boldsymbol{\tau}^{i}, \boldsymbol{v^{i}}\right)_{0} = \left\langle \gamma^{i}_{2}\boldsymbol{\tau}^{i}\cdot\boldsymbol{\alpha}^{i}, \gamma^{i}_{2} \boldsymbol{v}^{i} \right\rangle_{2} + \left\langle \gamma^{i}_{3}\boldsymbol{\tau}^{i}\cdot\boldsymbol{\beta}^{i}, \gamma^{i}_{3} \boldsymbol{v}^{i} \right\rangle_{3} \\
=& \left\langle {\tau}^{i}_{\alpha}, {v}^{i}_{\alpha} \right\rangle_{2} + \left\langle \boldsymbol{\tau}^{i}_{\tau}, \boldsymbol{v}^{i}_{\tau} \right\rangle_{2} + \left\langle {\tau}^{i}_{\beta}, {v}^{i}_{\beta} \right\rangle_{3} + \left\langle \boldsymbol{\tau}^{i}_{\eta}, \boldsymbol{v}^{i}_{\eta} \right\rangle_{3},
\end{aligned}
\end{equation}
where $\gamma^{i}_{2} \boldsymbol{\tau}^{i}\cdot\boldsymbol{\alpha}^{i} \in \mathcal{H}^{-1/2}(\Gamma^{i}_{2})^d$, $\gamma^{i}_{3} \boldsymbol{\tau}^{i}\cdot\boldsymbol{\beta}^{i} \in \mathcal{H}^{-1/2}(\Gamma^{i}_{3})^d$, ${\tau}^{i}_{\alpha} \in W^{i\prime}_{2N}$, $\boldsymbol{\tau}^{i}_{\tau} \in W^{i\prime}_{2T}$, ${\tau}^{i}_{\beta} \in W^{i\prime}_{3N}$, $\boldsymbol{\tau}^{i}_{\eta} \in W^{i\prime}_{3T}$ and $\left\langle\cdot,\cdot\right\rangle_{r}$ is the scalar product in $L^{2}(\Gamma^{i}_{r})$ ($r=2$ or $3$).}
It should be emphasized that the normal stresses $\sigma_{\beta}^{i}$, $\sigma_{\alpha}^{i+1}$ and normal displacement $[v^{i}_{N}]$ on contact boundary zone $\Gamma_{c}^{i}$ satisfy the following conditions:
$$
\sigma_{\beta}^{i} \cdot [v^{i}_{N}] = 0,~~\sigma_{\beta}^{i} = -\sigma_{\alpha}^{i+1}.
$$
Therefore, only friction boundary conditions on the contact zone $\Gamma_{c}^{i}$ need to be considered. For convenience, on $\Gamma_{c}^{i}$ let $\sigma_{N}^{i} = \sigma_{\beta}^{i}$, $\boldsymbol{\sigma}_{T}^{i} = \boldsymbol{\sigma}_{\eta}^{i}$, $v_{N}^{i} = v_{\beta}^{i}$, $\boldsymbol{v}_{T}^{i} = \boldsymbol{v}_{\eta}^{i}$ and $\delta_c^i: V^i \rightarrow L^2\left(\Gamma_c^i\right) \times L^2\left(\Gamma_c^i\right)^d$. 
Analogously, $\delta_c^i \boldsymbol{v}^i=\left(v_N^i, \boldsymbol{v}_T^i\right)$,
$$
\delta_c^i V^i=W_c^i=W_{c N}^i \times W_{c T}^i \text { on } \Gamma_c^i, 
$$
and the dual spaces of $W_{c N}^i$, $W_{c T}^i$ and $W_{c}^i$ can be denoted by $W^{i\prime}_{c N}$, $W^{i\prime}_{c T}$ and $W^{i\prime}_{c}$. 
Differently, let $W_{c} = W^{1}_{c} \times\cdots\times W^{n-1}_{c}$ and $W^{\prime}_{c} = W^{1\prime}_{c} \times\cdots\times W^{n-1\prime}_{c}$. The norm in the corresponding space is defined in the same way.

Moreover, let 
$$
W_{+}^{i\prime}=\left\{{\omega}_{N}^{i*} \in W^{i\prime}_{cN} \mid \left\langle {\omega}_{N}^{i*}, v^{i}_N\right\rangle_{c} \geqslant 0 ~~\forall \boldsymbol{v}^{i} \in V^{i}, v^{i}_{N} \geqslant 0 \text { on } \Gamma_c^{i}\right\}
$$
and
$$
W_{-}^{i\prime}=-W_{+}^{i\prime}
$$
be the closed convex cones of non-negative and non-positive functionals in $W^{i\prime}_{c N}$, respectively, where $\langle\cdot, \cdot\rangle_c$ is the scalar product in $L^2\left(\Gamma_c^i\right)^{d}$.

\begin{remark}[]\label{4:rem:dual.space}
\begin{itemize}
\item Introduce the following variational problem:
\begin{equation}\label{4:eq:aux.var}
\left(\boldsymbol{\varepsilon}(\boldsymbol{u}^{i}), \boldsymbol{\varepsilon}(\boldsymbol{v}^{i})\right)_0 + \left(\boldsymbol{u}^{i}, \boldsymbol{v}^{i}\right)_0 = \left\langle\boldsymbol{\varphi}^{i*}_{2}, \gamma^{i}_{2} \boldsymbol{v}^{i}\right\rangle _{2} + \left\langle\boldsymbol{\varphi}^{i*}_{3}, \gamma^{i}_{3} \boldsymbol{v}^{i}\right\rangle _{3} ~~\forall \boldsymbol{v}^{i} \in V^{i}.
\end{equation}
For any fixed $\boldsymbol{\varphi}^{i*}_{2}\in \mathcal{H}^{-1/2}(\Gamma^{i}_{2})^{d}$ and $\boldsymbol{\varphi}^{i*}_{3}\in \mathcal{H}^{-1/2}(\Gamma^{i}_{3})^{d}$, the above variational equation has precisely one solution $\boldsymbol{u}^{i}(\boldsymbol{\varphi}^{i*}_{2},\boldsymbol{\varphi}^{i*}_{3})$ and $\boldsymbol{\tau}^{i} = \boldsymbol{\varepsilon}\left(\boldsymbol{u}^{i}(\boldsymbol{\varphi}^{i*}_{2},\boldsymbol{\varphi}^{i*}_{3})\right) \in Q^{i}_{1}$. Combining this and Green's formula (\ref{4:eq:Green}), it can be deduced that $\boldsymbol{\tau}^{i}\cdot\boldsymbol{\alpha}^{i} = \boldsymbol{\varphi}^{i*}_{2}$ on $\Gamma^{i}_{2}$ and $\boldsymbol{\tau}^{i}\cdot\boldsymbol{\beta}^{i} = \boldsymbol{\varphi}^{i*}_{3}$ on $\Gamma^{i}_{3}$, that is, operator $G^{i}_{2}: \boldsymbol{\tau}^i \mapsto \boldsymbol{\tau}^i\cdot\boldsymbol{\alpha}^{i}$ maps $Q^{i}_1$ onto $\mathcal{H}^{-1/2}(\Gamma^{i}_{2})^d$ and operator $G^{i}_{3}: \boldsymbol{\tau}^i \mapsto \boldsymbol{\tau}^i\cdot\boldsymbol{\beta}^{i}$ maps $Q^{i}_1$ onto $\mathcal{H}^{-1/2}(\Gamma^{i}_{3})^d$. Similarly, it can be proven that
\begin{align*}
&G^{i}_{2N}: \boldsymbol{\tau}^i \mapsto {\tau}^i_{\alpha} \text{ maps } Q^{i}_{1} \text{ onto } W^{i\prime}_{2N};\\
&G^{i}_{2T}: \boldsymbol{\tau}^i \mapsto \boldsymbol{\tau}^i_{\tau} \text{ maps } Q^{i}_{1} \text{ onto } W^{i\prime}_{2T};\\
&G^{i}_{3N}: \boldsymbol{\tau}^i \mapsto {\tau}^i_{\beta} \text{ maps } Q^{i}_{1} \text{ onto } W^{i\prime}_{3N};\\
&G^{i}_{3T}: \boldsymbol{\tau}^i \mapsto \boldsymbol{\tau}^i_{\eta} \text{ maps } Q^{i}_{1} \text{ onto } W^{i\prime}_{3T}.
\end{align*}
\item The dual spaces $\mathcal{H}^{-1/2}(\Gamma^{i}_{r})^{d}$ and $W^{i\prime}_{r}$ are mutually isomorphic. Based on the above discussion, let $\boldsymbol{\mu}^{i*}_{2} = \left( {\tau}^i_{\alpha},\boldsymbol{\tau}^i_{\tau} \right) \in W^{i\prime}_{2}$ and $\boldsymbol{\mu}^{i*}_{3} = \left( {\tau}^i_{\beta},\boldsymbol{\tau}^i_{\eta} \right)\in W^{i\prime}_{3}$, where $\boldsymbol{\tau}^i$ satisfies $\boldsymbol{\tau}^{i}\cdot\boldsymbol{\alpha}^{i} = \boldsymbol{\varphi}^{i*}_{2}$ on $\Gamma^{i}_{2}$ and $\boldsymbol{\tau}^{i}\cdot\boldsymbol{\beta}^{i} = \boldsymbol{\varphi}^{i*}_{3}$ on $\Gamma^{i}_{3}$. Thus, $\boldsymbol{\mu}^{i*}_{2}$ and $\boldsymbol{\mu}^{i*}_{3}$ are determined solely by $\boldsymbol{\varphi}^{i*}_{2}$ and $\boldsymbol{\varphi}^{i*}_{3}$, respectively. Therefore, let $\boldsymbol{\mu}^{i*}_{r} = \varrho^{i*}_{r} \boldsymbol{\varphi}^{i*}_{r}$, $\varrho^{i*}_{r}: \mathcal{H}^{-1/2}(\Gamma^{i}_{r})^{d} \to W^{i\prime}_{r}$, $r=2$ or $3$, where $\boldsymbol{\mu}^{i*}_{r}$ is independent of the choice of $\boldsymbol{\tau}^{i} \in Q^{i}_{1}$. Finally, based on the fact that $W_r^i$ and $\mathcal{H}^{1 / 2}\left(\Gamma_r^i\right)^d$ are isomorphic, the conclusion is verified.
\item Furthermore, it can be proved that 
$$
\left\|\boldsymbol{\varphi}^{i*}_{r}\right\|_{-1/2,\Gamma^{i}_{r}} = \left\| \boldsymbol{\tau}^{i} \right\|_{1} = \left\| \boldsymbol{u}^{i}(\boldsymbol{\varphi}^{i*}_{r}) \right\|_{1}, ~\forall \boldsymbol{\varphi}^{i*}_{r}\in \mathcal{H}^{-1/2}(\Gamma^{i}_{r})^{d}
$$
and
$$
\left\|\varrho^{i*}_r\boldsymbol{\varphi}^{i*}_{r}\right\|_{W^{i\prime}_{r}} = \left\|\boldsymbol{\varphi}^{i*}_{r}\right\|_{-1/2,\Gamma^{i}_{r}} = \sup_{\boldsymbol{v}^{i}\in V^{i}} \frac{ \left\langle \boldsymbol{\varphi}^{i*}_{r}, \gamma^{i}_{r} \boldsymbol{v}^{i} \right\rangle _{r}}{\left\| \boldsymbol{v}^{i} \right\|_{V^i}}, ~ \forall \boldsymbol{\varphi}^{i*}_{r}\in \mathcal{H}^{-1/2}(\Gamma^{i}_{r})^{d},
$$
where let $\boldsymbol{\varphi}^{i*}_{5-r} = \boldsymbol{0}$ in (\ref{4:eq:aux.var}).
The specific proof process can be found in \cite{haslinger1981mixed}. This conclusion can also be extended to $\Gamma^{i}_{c}$.
\item It is worth noting that the above conclusions can be extended to the contact zone $\Gamma^{i}_{c}$, and the corresponding operators are defined following the rules mentioned above.
\end{itemize}
\end{remark}

According to previous research \cite{zhang2022variational}, the mechanical response model of this multi-layer elastic system can be expressed as the following variational inequality problem:
\begin{problem}[$\mathcal{P}_{v}$]\label{4:prb:var.ieq}
Find a displacement $\boldsymbol{u} \in \mathcal{K}$ which satisfies:
\begin{equation}\label{4:ieq:var}
a(\boldsymbol{u}, \boldsymbol{v}-\boldsymbol{u})+j(\boldsymbol{v})-j(\boldsymbol{u}) \geqslant L(\boldsymbol{v}-\boldsymbol{u}), \forall \boldsymbol{v} \in \mathcal{K},
\end{equation}
where
\begin{align*}
& a(\boldsymbol{v}, \boldsymbol{w})=\sum_{i=1}^n a^i\left(\boldsymbol{v}^i, \boldsymbol{w}^i\right) =\sum_{i=1}^n \int_{\Omega^i} A^i: \boldsymbol{\varepsilon} \left(\boldsymbol{v}^i\right): \boldsymbol{\varepsilon}\left(\boldsymbol{w}^i\right) \mathrm{d} x,\\
& L(\boldsymbol{v})=\sum_{i=1}^n \int_{\Omega^i} \boldsymbol{f}_0^i \cdot \boldsymbol{v}^i d x+\int_{\Gamma_2^1} \boldsymbol{f}_1 \cdot \boldsymbol{v}^1 d s,\\
& j(\boldsymbol{w})=\sum_{i=1}^{n-1} j^i\left(\boldsymbol{w}^i, \boldsymbol{w}^{i+1}\right)=\sum_{i=1}^{n-1} \int_{\Gamma_c^i} g^i({x})\left|\left[\boldsymbol{w}_T^i\right]\right| d s, \quad \forall \boldsymbol{w} \in V,\\
& \mathcal{K}=\left\{\boldsymbol{v} \in V \mid\left[v_N^i\right] \leqslant 0 \text {, a.e. on } \Gamma_c^i, i=1,2, \ldots, n-1\right\} \subset V.
\end{align*}
\end{problem}
Then, this variational inequality can be equivalently translated into \textcolor{black}{an optimization} problem:
\begin{problem}[$\mathcal{P}_{o}$]\label{4:prb:optimization}
Find a displacement $\boldsymbol{u} \in \mathcal{K}$ which satisfies:
\begin{equation}\label{4:eq:optimization}
\boldsymbol{u} = \arg \min_{\boldsymbol{v} \in \mathcal{K}} J(v),
\end{equation}
where
\begin{equation}\label{4:eq:J}
J(\boldsymbol{v})=\frac{1}{2} a(\boldsymbol{v}, \boldsymbol{v}) + j(\boldsymbol{v}) - L(\boldsymbol{v}).
\end{equation}
\end{problem}
The equivalence of problem \nameref{4:prb:var.ieq} and problem \nameref{4:prb:optimization} can be proved by Ref. \cite{han2002quasistatic}.

\subsection{Equivalent formulation of \texorpdfstring{$\mathcal{P}_{v}$}{}}

However, the optimization problem is difficult to solve directly because the displacement field of the elastic bodies in the multi-layer elastic system will be coupled through the interlayer contact condition. Therefore, the mixed formulation of the problem \nameref{4:prb:optimization} will be introduced. To achieve this, the sets $\Lambda^{i} = \Lambda^{i}_{N} \times \Lambda^{i}_{T} \subset W^{i\prime}_{c}$ need to be defined as follows: 
\begin{align}
& \Lambda^{i}_N=W_{+}^{i\prime} \label{4:eq:LamN}\\
& \Lambda^{i}_T=\left\{\boldsymbol{\omega}_{T}^{i} \in W^{i\prime}_{cT} \Big| \left|\boldsymbol{\omega}_{T}^{i}\right| \leqslant 1 \text { a.e. on } \Gamma^{i}_{g} \text{ and } \boldsymbol{\omega}_{T}^{i} = 0 \text{ on } \Gamma^{i}_{c}-\Gamma^{i}_{g}\right\}, \label{4:eq:LamT}
\end{align}
where $\Gamma^{i}_{g} = supp\left(g^{i}(x)\right) \cap \Gamma^{i}_{c}$ and $supp\left(g^{i}(x)\right)$ is the support set of $g^{i}(x)$.
It can be verified that $\Lambda^{i}$ is closed and convex in $W^{i\prime}_{c}$. Then, let $\Lambda = \Lambda^{1} \times \cdots\times \Lambda^{n-1}$.

Then, the following saddle-point problem can be introduced:
\begin{problem}[$\mathcal{P}_{l}$]\label{4:prb:Lagrangian}
Find a solution $\left\{\boldsymbol{w}, \boldsymbol{\lambda}\right\} \in V \times \Lambda$ which satisfies the following saddle-point problem:
\begin{equation}\label{4:eq:saddle}
\mathcal{L}(\boldsymbol{w}, \boldsymbol{\mu}) \leqslant \mathcal{L}(\boldsymbol{w}, \boldsymbol{\lambda}) \leqslant \mathcal{L}(\boldsymbol{v}, \boldsymbol{\lambda}) \quad \forall \boldsymbol{v} \in V, \forall \boldsymbol{\mu} \in \Lambda,
\end{equation}
where the Lagrangian functional $\mathcal{L}: V \times \Lambda \to \mathbb{R}$ is defined as:
\begin{equation}\label{4:eq:Lagrangian}
\mathcal{L}\left( \boldsymbol{v}, \boldsymbol{\mu} \right)=\frac{1}{2} a(\boldsymbol{v}, \boldsymbol{v}) + \mathcal{G}(\boldsymbol{\mu}, \boldsymbol{v}) - L(\boldsymbol{v})
\end{equation}
and
\begin{equation}\label{4:eq:def:g}
\mathcal{G}(\boldsymbol{\mu}, \boldsymbol{v}) = \sum_{i=1}^{n-1} \left(\left\langle \mu^{i}_{N}, [v^{i}_{N}] \right\rangle_{c} + \left\langle g^{i}({x})\boldsymbol{\mu}^{i}_{T}, [\boldsymbol{v}^{i}_{T}] \right\rangle_{c} \right).
\end{equation}
\end{problem}
Based on optimization theory, it can be verified that problem \nameref{4:prb:Lagrangian} is equivalent to the following mixed problem:
\begin{problem}[$\mathcal{P}_{m}$]\label{4:prb:mix}
Find a solution $\left\{\boldsymbol{w}, \boldsymbol{\lambda}\right\} \in V \times \Lambda$ which satisfies:
\begin{equation}\label{4:eq:mix}
\left\{\begin{aligned}
&a(\boldsymbol{w}, \boldsymbol{v})+\mathcal{G}(\boldsymbol{\lambda}, \boldsymbol{v})=L(\boldsymbol{v}),~~ \forall \boldsymbol{v} \in V \\
&\mathcal{G}(\boldsymbol{\mu}-\boldsymbol{\lambda}, \boldsymbol{w}) \leqslant 0, ~~ \forall \boldsymbol{\mu} \in \Lambda
\end{aligned}\right. .
\end{equation}
\end{problem}

Then, the equivalence between problem \nameref{4:prb:var.ieq} and \nameref{4:prb:mix} needs to be verified. However, before presenting the relevant theorem, the following lemma should be explained.
\begin{lemma}[]\label{4:lem:sigma}
In problem \nameref{4:prb:var.ieq}, the stress ${\sigma}^{i}_{N}$ and $\boldsymbol{\sigma}^{i}_{T}$ on contact zone $\Gamma^{i}_{c}$ satisfy:
\begin{equation}
\left\{\begin{aligned}
& \sigma_N^i \in W^{i\prime}_{-}, ~~\left|\boldsymbol{\sigma}_T^i\right| \leqslant g^i(x) \\
& g^i(x)\left|\left[\boldsymbol{u}_T^i\right]\right|+\boldsymbol{\sigma}_T^i \cdot\left[\boldsymbol{u}_T^i\right]=0
\end{aligned} \right. \text{ a.e. on } \Gamma_{c}^{i}.
\end{equation}
\end{lemma}

\begin{proof}
Based on the variational inequality (\ref{4:ieq:var}) and \textcolor{black}{Green's formula}, it can be obtained that
\begin{equation}\label{4:ieq:Gamma.c}
\sum_{i=1}^{n-1} \left( \left\langle \sigma_N^i , \left(\left[v_N^i\right]-\left[u_N^i\right]\right) \right\rangle_{c} + \left\langle \boldsymbol{\sigma}_T^i, \left(\left[\boldsymbol{v}_T^i\right]-\left[\boldsymbol{u}_T^i\right]\right) \right\rangle_{c} \right) + \sum_{i=1}^{n-1} \int_{\Gamma_c^i} g^i(x)\left(\left|\left[\boldsymbol{v}_T^i\right]\right| - \left|\left[\boldsymbol{u}_T^i\right]\right| \right) d s \geqslant 0.
\end{equation}
Let $\boldsymbol{v} = \boldsymbol{u} \pm \boldsymbol{w}$, where $[w_{N}^{i}]= 0$ on $\Gamma^{i}_{c}$. Since $\boldsymbol{u}\in \mathcal{K}$, it can be verified that $\boldsymbol{v}\in \mathcal{K}$. Then, the inequality (\ref{4:ieq:Gamma.c}) can be written as:
\begin{equation}\label{4:proof:lem.1.1}
\sum_{i=1}^{n-1} \left\langle \boldsymbol{\sigma}_T^i, \pm\left[\boldsymbol{w}_T^i\right] \right\rangle _{c} + \sum_{i=1}^{n-1} \int_{\Gamma_c^i} g^i(x)\left(\left|\left[\boldsymbol{u}_T^i\right] \pm \left[\boldsymbol{w}_T^i\right]\right| - \left|\left[\boldsymbol{u}_T^i\right]\right| \right) d s \geqslant 0.
\end{equation}
As consequence,
\begin{align*}
&\sum_{i=1}^{n-1} \pm \left\langle \boldsymbol{\sigma}_T^i, \left[\boldsymbol{w}_T^i\right] \right\rangle _{c} + \sum_{i=1}^{n-1}\int_{\Gamma_c^i} g^i(x)\left|\left[\boldsymbol{w}_T^i\right]\right| d s \geqslant 0 \\
\Rightarrow & \sum_{i=1}^{n-1}\left| \left\langle \boldsymbol{\sigma}_T^i, \left[\boldsymbol{w}_T^i\right] \right\rangle _{c} \right| \leqslant \sum_{i=1}^{n-1}\int_{\Gamma_c^i} g^i(x)\left|\left[\boldsymbol{w}_T^i\right]\right| d s.
\end{align*}
Since $\forall \boldsymbol{w} \in V$ such that $[w_{N}^{i}]= 0$, the above relation holds, we have
$$
\boldsymbol{\sigma}^{i}_{T} \in L^{\infty}(\Gamma^{i}_{c}),~~ \left| \boldsymbol{\sigma}^{i}_{T} \right| \leqslant g^{i}(x) \text{ a.e on } \Gamma_{c}^{i}.
$$
and
\begin{equation}\label{4:proof:lem.1.2}
\boldsymbol{\sigma}_T^i \cdot \left[\boldsymbol{u}_T^i\right] + g^i(x) \left| \left[\boldsymbol{u}_T^i\right] \right| \geqslant 0 \text { a.e on } \Gamma_c^i.
\end{equation}

Let $\boldsymbol{w}_{k}\in V$ be a sequence of function, satisfying $\left[w_{k N}^i\right]=0$ on $\Gamma_{c}^{i}$ and $\left[\boldsymbol{w}_{k T}^i\right] \to -\left[\boldsymbol{u}_T^i\right]$ or $0$ in $L^{1}\left(\Gamma_{c}^{i}\right)$. Then, from (\ref{4:proof:lem.1.1}) it can be obtained that
\begin{align*}
& \sum_{i=1}^{n-1}\left\langle\boldsymbol{\sigma}_T^i, \left[\boldsymbol{w}_T^i\right]\right\rangle_c+\sum_{i=1}^{n-1} \int_{\Gamma_c^i} g^i(x)\left(\left|\left[\boldsymbol{u}_T^i\right] +\left[\boldsymbol{w}_T^i\right]\right|-\left|\left[\boldsymbol{u}_T^i\right]\right|\right) d s \geqslant0 \\
\Rightarrow & - \left\langle\boldsymbol{\sigma}_T^i, \left[\boldsymbol{u}_T^i\right]\right\rangle_c - \int_{\Gamma_c^i} g^i(x)\left|\left[\boldsymbol{u}_T^i\right]\right| d s \geqslant 0, ~~i=1,\ldots,n-1.
\end{align*}
Based on the above inequality and formula (\ref{4:proof:lem.1.2}), it can be found that
$$
\boldsymbol{\sigma}_T^i \cdot \left[\boldsymbol{u}_T^i\right] + g^i(x) \left| \left[\boldsymbol{u}_T^i\right] \right| = 0 \text { a.e on } \Gamma_c^i.
$$
Then, from (\ref{4:ieq:Gamma.c}), we have
\begin{equation*}
\sum_{i=1}^{n-1} \left( \left\langle \sigma_N^i , \left(\left[v_N^i\right]-\left[u_N^i\right]\right) \right\rangle _{c} + \left\langle \boldsymbol{\sigma}_T^i, \left[\boldsymbol{v}_T^i\right] \right\rangle_{c} \right) + \sum_{i=1}^{n-1} \int_{\Gamma_c^i} g^i(x)\left(\left|\left[\boldsymbol{v}_T^i\right]\right| \right) d s \geqslant 0.
\end{equation*}
By inserting ${v}^{i}_{N} = 0$, ${v}^{i}_{N} = 2{u}^{i}_{N}$ and $\boldsymbol{v}^{i}_{T} = 0$ into the above inequality, we can find that
$$
\left\langle \sigma_N^i,\left[u_N^i\right] \right\rangle _{c}=0 \Rightarrow \left\langle \sigma_N^i,\left[v_N^i\right] \right\rangle _{c} \geqslant 0.
$$
Since $\left[v_N^i\right]\leqslant 0$, it can be proved that $\sigma_N^i \in W_{-}^{i \prime}$.
\end{proof}

In fact, Lemma \ref{4:lem:sigma} represents the frictional contact conditions at the contact zone $\Gamma_{c}^{i}$. Based on this lemma, the following theorem can be formulated.

\begin{theorem}[]\label{4:thm:saddle.point}
There exists a unique saddle-point $\left\{ \boldsymbol{w}, \boldsymbol{\lambda} \right\}$ of $\mathcal{L}(\boldsymbol{v},\boldsymbol{\mu})$ such that
$$
\boldsymbol{w} = \boldsymbol{u}, ~\lambda^{i}_{N} = -\sigma^{i}_{N}, ~g^{i}(x)\boldsymbol{\lambda}^{i}_{T} = -\boldsymbol{\sigma}^{i}_{T},
$$
where $\boldsymbol{u}$ is the solution of problem \nameref{4:prb:var.ieq}.
\end{theorem}
\begin{proof}
Since the \nameref{4:prb:mix} is equivalent to the \nameref{4:prb:Lagrangian}, based on the first equation of (\ref{4:eq:mix}) and \textcolor{black}{Green's formula} (\ref{4:eq:Green}), the following relation can be obtained:
$$
\mathcal{G}(\boldsymbol{\lambda},\boldsymbol{v}) = - \sum_{i=1}^{n-1}\left( \left\langle\sigma_N^i, \left[v_N^i\right]\right\rangle _c + \left\langle\boldsymbol{\sigma}_T^i, \left[\boldsymbol{v}_T^i\right] \right\rangle _c \right).
$$
The above equation holds $\forall \boldsymbol{v} \in V$, so $\lambda_N^i=-\sigma_N^i(\boldsymbol{\omega}) \text { and } g^i(x)\boldsymbol{\lambda}_T^i=-\boldsymbol{\sigma}_T^i(\boldsymbol{\omega})$. 

Then, the second inequality of (\ref{4:eq:mix}) should be analyzed. By introducing $\mu^{i}_{N} = 0$, $\boldsymbol{\mu}^{i}_{T} = \boldsymbol{\lambda}^{i}_{T}$, and $\mu^{i}_{N} = 2\lambda^{i}_{N}$, $\boldsymbol{\mu}^{i}_{T} = \boldsymbol{\lambda}^{i}_{T}$, we have
\begin{align*}
& \sum_{i=1}^{n-1}\left\langle\lambda_N^i,\left[w_N^i\right]\right\rangle _c \leqslant 0 \text{ and } -\sum_{i=1}^{n-1}\left\langle\lambda_N^i,\left[w_N^i\right]\right\rangle _c \leqslant 0 \\
\Rightarrow & \sum_{i=1}^{n-1}\left\langle\lambda_N^i,\left[w_N^i\right]\right\rangle _c = 0 \Rightarrow \sum_{i=1}^{n-1}\left\langle\mu_N^i,\left[w_N^i\right]\right\rangle _c \leqslant 0,~~ \forall {\mu}^{i}_{N} \in \Lambda^{i}_{N}. 
\end{align*}
Therefore, it can be deduced that $\left[w_N^i\right] \leqslant 0$ a.e. on $\Gamma^{i}_{c}$, which implies $\boldsymbol{w}\in\mathcal{K}$. Similarly, by introducing $\mu^{i}_{N} = \lambda^{i}_{N}$ and $\boldsymbol{\mu}_T^i = \boldsymbol{\mu}_T^i$ or $\boldsymbol{\lambda}_T^i$ into the inequality of (\ref{4:eq:mix}), we have
\begin{align*}
& \left\langle g^i(x) \boldsymbol{\mu}_T^i,\left[\boldsymbol{w}_T^i\right] \right\rangle _c \leqslant \left\langle g^i(x) \boldsymbol{\lambda}_T^i,\left[\boldsymbol{w}_T^i\right] \right\rangle _c, ~\forall \boldsymbol{\mu}_T^i \in \Lambda^{i}_{T}, \\
\Rightarrow & \int_{\Gamma^{i}_{c}} g^{i}(x) \left| \left[\boldsymbol{w}_T^i\right] \right| ds \leqslant \left\langle g^i(x) \boldsymbol{\lambda}_T^i,\left[\boldsymbol{w}_T^i\right] \right\rangle _c.
\end{align*}

Then, by introducing $\boldsymbol{v} = \boldsymbol{v} - \boldsymbol{w}$ into the equation of (\ref{4:eq:mix}), it can be obtained that the following inequality holds for all $\boldsymbol{v}\in V$:
\begin{align*}
& a(\boldsymbol{w}, \boldsymbol{v} - \boldsymbol{w}) + \mathcal{G}(\boldsymbol{\lambda}, \boldsymbol{v}- \boldsymbol{w}) =L(\boldsymbol{v} - \boldsymbol{w})\\
\Leftrightarrow & a(\boldsymbol{w}, \boldsymbol{v} - \boldsymbol{w}) + \sum_{i=1}^{n-1}\left(\left\langle\lambda_N^i,\left[v_N^i\right] - \left[w_N^i\right]\right\rangle_c+\left\langle g^i(x) \boldsymbol{\lambda}_T^i, \left[\boldsymbol{v}_T^i\right] - \left[\boldsymbol{w}_T^i\right]\right\rangle_c\right) = L(\boldsymbol{v} - \boldsymbol{w}) \\
\Rightarrow & a(\boldsymbol{w}, \boldsymbol{v} - \boldsymbol{w}) + \sum_{i=1}^{n-1}\left(\left\langle\lambda_N^i,\left[v_N^i\right]\right\rangle_c + \left\langle g^i(x) \boldsymbol{\lambda}_T^i, \left[\boldsymbol{v}_T^i\right]\right\rangle_c - \int_{\Gamma_c^i} g^i(x)\left|\left[\boldsymbol{w}_T^i\right]\right| d s \right) \geqslant L(\boldsymbol{v} - \boldsymbol{w}).
\end{align*}
Then, let $\boldsymbol{v} \in \mathcal{K}$ and since $\lambda^{i}_{N}\in\Lambda^{i}_{N}$ and $\boldsymbol{\lambda}^{i}_{T}\in\Lambda^{i}_{T}$, we have
$$
a(\boldsymbol{w}, \boldsymbol{v} - \boldsymbol{w}) + \sum_{i=1}^{n-1}\left( \int_{\Gamma_c^i} g^i(x)\left|\left[\boldsymbol{v}_T^i\right]\right| d s - \int_{\Gamma_c^i} g^i(x)\left|\left[\boldsymbol{w}_T^i\right]\right| d s \right) \geqslant L(\boldsymbol{v} - \boldsymbol{w}), \forall \boldsymbol{v} \in \mathcal{K},
$$
where $\boldsymbol{w} \in \mathcal{K}$. Therefore, $\boldsymbol{w}$ is the solution of \nameref{4:prb:var.ieq}, that is $\boldsymbol{w} = \boldsymbol{u}$.

Conversely, let $\boldsymbol{u}\in \mathcal{K}$ be the solution of problem \nameref{4:prb:var.ieq}. Based on \textcolor{black}{Lemma \ref{4:lem:sigma}}, it can be verified that $\left\{\boldsymbol{u}^{i}, -\sigma^{i}_{N}, - \boldsymbol{\sigma}^{i}_{T}/g^{i}(x)\right\} \in V^{i} \times \Lambda^{i}_{N} \times \Lambda^{i}_{T}$ is the solution of problem \nameref{4:prb:Lagrangian}, where if $g^{i}(x) = 0$, let $- \boldsymbol{\sigma}^{i}_{T}/g^{i}(x) = \boldsymbol{0}$. So, the equivalence between problem \nameref{4:prb:mix} and problem \nameref{4:prb:var.ieq} is proved. Finally, according to the previous work \cite{zhang2022variational}, since the solution of variational inequality problem \nameref{4:prb:var.ieq} exists and is unique, the solution $\{\boldsymbol{w}, \boldsymbol{\lambda}\}$ to problem \nameref{4:prb:mix} also exists and is unique.
\end{proof}

\section{Finite element approximation and convergence analysis}

\subsection{Finite element spaces and discrete problem}

Before using the finite element method to solve the approximate solution of the problem \nameref{4:prb:mix}, the corresponding finite element spaces need to be constructed. 
Let $E^{i}_{h} = \left\{T^{i}_{1},\ldots,T^{i}_{N^{i}_{h}}\right\}$ be a non-degenerate quasi-uniform subdivision (triangular or tetrahedral elements) of the $i$-th layer elastic body $\bar{\Omega}^{i}$, where $T_j^i\left(i=1,2, \ldots, n; j=1,2, \ldots, N_h^i\right)$ is a regular finite element, $h$ is the maximum diameter of all elements in the system, and $N^{i}_{h}$ is the number of elements in $\bar{\Omega}^{i}$. 
Then $E_{h} = \cup^{n}_{i=1} E^{i}_{h}$ represents the finite element division of the multi-layer elastic system. It should be noted that on the potential contact boundaries $\Gamma^{i}_{3}$ and $\Gamma^{i+1}_{2}$, the finite element divisions $E^{i}_{h}$ and $E^{i+1}_{h}$ may not be compatible.
Therefore, let $F^{i}_{H} = \left\{S^{i}_{1},\ldots,S^{i}_{N^{i}_{H}}\right\}$ denote a regular subdivision (linear or triangular elements) of the contact boundary $\bar{\Gamma}^{i}_{c}$, independent of $E_{h}$, where $S_j^i\left(i=1, \ldots, n-1; j=1, \ldots, N_H^i\right)$ represents a regular contact element, $H$ denotes the maximum diameter of elements in all contact zones, and $N^{i}_{H}$ is the total number of elements in $\bar{\Gamma}^{i}_{c}$.  
Let $F_{H} = \bigcup_{i=1}^{n-1} F^{i}_{H}$. It is worth noting that, under the assumption of a regular finite element division, \textcolor{black}{there exists} a positive constant $c_{H}$ such that $H^{i}_{j}/H \geqslant c_{H}$ for all $j \in \{1, \ldots, N^{i}_{H}\}$ and $i \in \{1, \ldots, n-1\}$, where $H^{i}_{j}$ denotes the diameter of the contact element $S^{i}_{j}$.
Finally, it is agreed that when $H = h$, the finite element subdivisions $E^{i}_{h}$ and $E^{i+1}_{h}$ are compatible on the contact zone $\Gamma^{i}_{c}$, and consequently, the finite element subdivision $F^{i}_{h}$ of $\Gamma^{i}_{c}$ is derived from $E_{h}$. 
Furthermore, we assume that all aforementioned finite element subdivisions are compatible with the boundaries of the support sets $\Gamma^{i}_{g}$, for $i=1,\ldots,n-1$.

Upon completion of the finite element division of the multi-layer system, the corresponding finite element spaces can now be defined as follows:
\begin{align}
& V_h^i=\left\{\boldsymbol{v}_h^i \in C\left(\bar{\Omega}^i\right)^d ~\Big|~ \boldsymbol{v}_h^i|_{T_j^i} \in P_{k_{v}}\left(T_j^i\right)^d, \forall T^{i}_{j}\in E^{i}_{h}, \boldsymbol{v}_h^i=0 \text { on } \Gamma_1^i\right\}, \label{4:space:Vih}\\
& {W}^{i\prime}_{cH}=\left\{\boldsymbol{\mu}^{i}_{H} \in C\left(\bar{\Gamma}^{i}_{c}\right)^{d+1} ~\Big|~ \boldsymbol{\mu}^{i}_{H}|_{S^{i}_{j}} \in P_{k_{b}}\left(S^{i}_{j}\right)^{d+1}, \forall S^{i}_{j}\in F^{i}_{H} \right\}, \label{4:space:WicH}\\
& \Lambda^{i}_{HN}=\left\{\mu^{i}_{hN} \in C\left(\bar{\Gamma}^{i}_{c}\right) ~\Big|~ \mu^{i}_{hN}|_{S^{i}_{j}} \in P_{k_{b}}\left(S^{i}_{j}\right), \mu^{i}_{hN} \geqslant 0 \text { on } \Gamma^{i}_{c}, \forall S^{i}_{j} \in F^{i}_{H}\right\}, \label{4:space:LiNH}\\
& \Lambda^{i}_{HT}=\left\{\boldsymbol{\mu}^{i}_{hT} \in C\left(\bar{\Gamma}^{i}_{c}\right)^d ~\Big|~ \boldsymbol{\mu}^{i}_{hT}|_{S^{i}_{j}} \in P_{k_{b}} \left(S^{i}_{j}\right)^d, \left|\boldsymbol{\mu}^{i}_{hT}\right| \leqslant 1 \text { on } \Gamma^{i}_{c}, \forall S^{i}_{j} \in F^{i}_{H}\right\}, \label{4:space:LiTH}\\
& V_{h} = V_h^1 \times \cdots \times V_h^n,~ {W}^{\prime}_{cH} = {W}^{1\prime}_{cH} \times \cdots \times {W}^{n-1\prime}_{cH},  \label{4:space:VhWH}\\
& \Lambda^{i}_{H} = \Lambda^{i}_{N H} \times \Lambda^{i}_{T H}, ~ \Lambda_{H} = \Lambda^{1}_{H} \times \cdots \times \Lambda^{n-1}_{H}, \label{4:space:LH}
\end{align}
where $P_{k_{v}}\left(T_j^i\right)$ and $P_{k_{b}}\left(S^{i}_{j}\right)$ are the space of $k_{v}$-order polynomial functions on $T_j^i$ and space of $k_{b}$-order polynomial functions on $S^{i}_{j}$, respectively. In particular, $P_0$ represents a space of constant functions. To ensure homogeneity between the contact element space of the boundary and the gradient of the displacement, it can be defined that $k_{b} = k_{v}-1$.
It is easy to verify that $V_{h}\subset V$, ${W}^{\prime}_{cH}\subset {W}^{\prime}_{c}$ and $\Lambda_{H} \subset \Lambda$, which will be used as finite element approximation spaces of problem \nameref{4:prb:mix}. Furthermore, the finite element approximation space of $\mathcal{K}$ is defined by
\begin{equation}\label{4:eq:KhH}
\mathcal{K}_{hH} = \left\{ \boldsymbol{v}_{h}\in V_{h} ~\Big|~ \left\langle \mu^{i}_{HN}, \left[v^{i}_{hN}\right] \right\rangle _{c} \leqslant 0, ~\forall \mu^{i}_{HN} \in \Lambda^{i}_{HN}, ~i=1,\ldots,n-1 \right\}.
\end{equation}

Then, the discrete form of problem \nameref{4:prb:mix} can be expressed as:

\begin{problem}[$\mathcal{P}_{hH}$]\label{4:prb:mixhH}
Find a solution $\{\boldsymbol{u}_{h}, \boldsymbol{\lambda}_{H}\} \in V_{h} \times \textcolor{black}{\Lambda_{H}}$ which satisfies:
\begin{equation}\label{4:eq:mixhH}
\left\{\begin{aligned}
& a(\boldsymbol{u}_{h}, \boldsymbol{v}_{h})+\mathcal{G}(\boldsymbol{\lambda}_{H}, \boldsymbol{v}_{h})=L(\boldsymbol{v}_{h}), \quad \forall \boldsymbol{v}_{h} \in V_{h} \\
& \mathcal{G}(\boldsymbol{\mu}_{H}-\boldsymbol{\lambda}_{H}, \boldsymbol{u}_{h}) \leqslant 0, \quad \forall \boldsymbol{\mu}_{H} \in \Lambda_{H}
\end{aligned}\right. .
\end{equation}
\end{problem}
Based on the equivalence between \nameref{4:prb:var.ieq} and \nameref{4:prb:mix}, it can be verified that the first component of the numerical solution $\{\boldsymbol{u}_{h}, \boldsymbol{\lambda}_{H}\}$ to problem \nameref{4:prb:mixhH} exists and is unique. 

\begin{remark}[]\label{4:rem:condition.unique}
In order to ensure the uniqueness of the second component $\boldsymbol{\lambda}_{H}$, the following condition is assumed:
\begin{equation}\label{4:eq:lambdaH.unique}
\sup _{V_h} \frac{ \mathcal{G} \left(\boldsymbol{\mu}_H, \boldsymbol{v}_h\right)}{\left\|\boldsymbol{v}_h\right\|_{V}} \geqslant c_{hH} \left\|\boldsymbol{\mu}_{H}\right\|_{W'_{c}},~~ \forall \boldsymbol{\mu}_H \in \Lambda_H .
\end{equation}
The above condition ensures that 
$$
\mathcal{G} \left(\boldsymbol{\mu}_H, \boldsymbol{v}_h\right) = 0 ~ \forall\boldsymbol{v}_h\in V_{h}   \Rightarrow \boldsymbol{\mu}_H =0, 
$$
then $\boldsymbol{\lambda}_{H}$ is uniquely determined. 

When the finite element space is a piecewise linear polynomial space, that is $k_{b}=1$, and $H=h$, it can be verified that the condition (\ref{4:eq:lambdaH.unique}) is true. However, if the finite element mesh nodes on the contact boundary do not match, the assessment of the condition (\ref{4:eq:lambdaH.unique}) will become more complex and require a case-by-case discussion. Nevertheless, based on the relationship between $H$ and $h$, it can be confirmed that when $h/H\leqslant c_{H}$, the condition (\ref{4:eq:lambdaH.unique}) must hold true. Lastly, if the finite element space of the contact boundary $\Gamma^{i}_{c}$ is a piecewise constant function space, the condition (\ref{4:eq:lambdaH.unique}) above is satisfied if $h/H$ is sufficiently small. Detailed analysis and examples can be found in Ref.\cite{haslinger1981mixed,haslinger1982approximation}.
\end{remark}

\subsection{Convergence analysis}

Without making any assumptions about the regularity and smoothness of $\boldsymbol{u}$, the convergence of $\left\{\boldsymbol{u}_h, \boldsymbol{\lambda}_{H}\right\}$ to $\left\{\boldsymbol{u}, \boldsymbol{\lambda}\right\}$ will be discussed. In the finite element framework, it is assumed that $h\to 0$ if and only if $H\to 0$. Subsequently, the following convergence theorem is presented:
\begin{theorem}[]\label{4:thm:convergence}
If the finite element space $\mathcal{K}_{hH}$ satisfies that $\forall \boldsymbol{v} \in \mathcal{K}$, $\exists$ a sequence $\left\{ \boldsymbol{v}_{h} \right\}$, where $\boldsymbol{v}_{h} \in \mathcal{K}_{hH}$, such that $\boldsymbol{v}_{h}\to\boldsymbol{v}$ in $V$ for \textcolor{black}{$h\to 0+$}, then:
\begin{align*}
&\lim_{h\to 0+}\boldsymbol{u}_{h} \to \boldsymbol{u} \text{ in } V,\\
&\lim_{H\to 0+}\boldsymbol{\lambda}^{i}_{HT} \to \boldsymbol{\lambda}^{i}_{T} \text{ in } W^{i\prime}_{cT},~i=1,\ldots,n-1
\end{align*}
where $\left\{\boldsymbol{u},\boldsymbol{\lambda}\right\}$ and $\left\{\boldsymbol{u}_{h},\boldsymbol{\lambda}_{H}\right\}$ are the solutions of problem \nameref{4:prb:mix} and problem \nameref{4:prb:mixhH}, respectively.
\end{theorem}

\begin{proof}
Let $\left\{\boldsymbol{u}_{h},\boldsymbol{\lambda}_{H}\right\}$ be a solution of the discrete problem \nameref{4:prb:mixhH}. Based on the characteristics of problems \nameref{4:prb:mixhH} and \nameref{4:prb:mix}, it can be confirmed that $\boldsymbol{u}_{h}\in \mathcal{K}_{hH}$ and is a solution of the following variational inequality:
\begin{equation}\label{4:ieq:var.hH}
a(\boldsymbol{u}_{h},\boldsymbol{v}_{h}-\boldsymbol{u}_{h}) + \sum_{i=1}^{n-1} \left\langle g^{i}(x)\boldsymbol{\lambda}^{i}_{HT}, \left[ \boldsymbol{v}^{i}_{hT} \right] - \left[ \boldsymbol{u}^{i}_{hT} \right] \right\rangle _{c} \geqslant L\left( \boldsymbol{v}_{h} - \boldsymbol{u}_{h} \right), ~\forall \boldsymbol{v}_{h}\in \mathcal{K}_{hH}.
\end{equation}
According to the definition of $\Lambda^{i}_{HT}$, sequence $\left\{\boldsymbol{\lambda}^{i}_{HT}\right\}$ is bounded in $L^{2}(\Gamma^{i}_{c})^{d}$, so the sequence $\left\{\boldsymbol{u}^{i}_{h}\right\}$ is bounded in $V^{i}$. Therefore, there exists a subsequence $\left\{\boldsymbol{\lambda}^{i}_{\bar{H}T}\right\} \subset \left\{\boldsymbol{\lambda}^{i}_{HT}\right\}$ and $\left\{\boldsymbol{u}^{i}_{\bar{h}}\right\} \subset \left\{\boldsymbol{u}^{i}_{h}\right\}$ such that
\begin{align}
& \boldsymbol{\lambda}^{i}_{\bar{H}T} \rightharpoonup \bar{\boldsymbol{\lambda}}^{i}_{T},~\text{ in } L^{2}(\Gamma^{i}_{c})^{d},~i=1,\ldots,n-1, \label{4:eq:proof:thm:convergence.1}\\
& \boldsymbol{u}^{i}_{\bar{h}}	\rightharpoonup \bar{\boldsymbol{u}}^{i},~\text{ in } V^{i}, ~i=1,\ldots,n, \label{4:eq:proof:thm:convergence.2}
\end{align}
where $\left\{\bar{\boldsymbol{u}}^{i}, \bar{\boldsymbol{\lambda}}^{i}_{T}\right\} \in V^{i} \times L^{2}(\Gamma^{i}_{c})^{d}$.

Based on the definition of $\Lambda_{HT}$, it can be verified that $\bar{\boldsymbol{\lambda}}^{i}_{T}\in \Lambda^{i}_{T}$. Then, it should be proved that $\bar{\boldsymbol{u}}\in \left(\bar{\boldsymbol{u}}^{1},\ldots,\bar{\boldsymbol{u}}^{n}\right)\in \mathcal{K}$. To this end, it is sufficient to verify that $\left[\bar{u}^{i}_{N}\right]\leqslant 0$ a.e. on $\Gamma^{i}_{c}$, which can be equivalently expressed as
\begin{equation}\label{4:eq:proof:thm:convergence.3}
\left\langle \mu^{i}_{N}, \left[\bar{u}^{i}_{N}\right] \right\rangle _c \leqslant 0, ~~\forall \mu^{i}_{N}\in L^{2}_{+}\left(\Gamma^{i}_{c}\right),~~i=1,\ldots,n-1,
\end{equation}
where $L^{2}_{+}\left(\Gamma^{i}_{c}\right) = \left\{ \mu^{i}_{N} \in L^{2}\left(\Gamma^{i}_{c}\right) ~\Big|~ \mu^{i}_{N} \geqslant 0 \text{ a.e. on } \Gamma^{i}_{c} \right\}$. Therefore, let $\mu^{i}_{N}$ be an arbitrary, fixed element in $L^{2}_{+}\left(\Gamma^{i}_{c}\right)$. Then, it can be found that there exists a sequence $\left\{ \mu^{i}_{HN} \right\}$, $\mu^{i}_{HN} \in \Lambda^{i}_{HN}$ such that 
\begin{equation}\label{4:eq:proof:thm:convergence.4}
\mu^{i}_{HN} \to \mu^{i}_{N},~H\to 0+ \text{ in } L^{2}\left(\Gamma^{i}_{c}\right).
\end{equation}
Since $\boldsymbol{u}_{\bar{h}} \in \mathcal{K}_{hH}$, it can be deduced that 
$$
\left\langle \mu^{i}_{HN}, \left[u^{i}_{\bar{h}N}\right] \right\rangle _c \leqslant 0. 
$$
Based on (\ref{4:eq:proof:thm:convergence.2}) and (\ref{4:eq:proof:thm:convergence.4}), we have 
$$
\lim_{\bar{h},H\to 0+} \left\langle \mu^{i}_{HN}, \left[u^{i}_{\bar{h}N}\right] \right\rangle _c = \left\langle \mu^{i}_{N}, \left[\bar{u}^{i}_{N}\right] \right\rangle _c \leqslant 0,
$$
so (\ref{4:eq:proof:thm:convergence.3}) is proved.

Then, the conditions verifying that $\bar{\boldsymbol{u}}$ is the solution of problem \nameref{4:prb:mix} are satisfied. According to the assumption, $\forall \boldsymbol{v}\in \mathcal{K}$, there exists a sequence $\{\boldsymbol{v}_{h}\}$, where $\boldsymbol{v}_{h}\in \mathcal{K}_{hH}$, such that
\begin{equation}\label{4:eq:proof:thm:convergence.5}
\boldsymbol{v}_{h} \to \boldsymbol{v} \text{ in } V.
\end{equation}
According to \textcolor{black}{Definition \ref{4:def:1}}, it can be verified that the mapping $\gamma^{i}_{r}: V^{i} \to L^{2}\left(\Gamma^{2}_{r}\right)^{d}$ ($r=2,3$) is completely continuous. By utilizing formulas (\ref{4:ieq:var.hH}), (\ref{4:eq:proof:thm:convergence.1}), (\ref{4:eq:proof:thm:convergence.2}), and (\ref{4:eq:proof:thm:convergence.5}), and considering the limit $\bar{h},\bar{H}\to 0+$, the following inequality can be deduced:
\begin{equation}\label{4:eq:proof:thm:convergence.6}
a(\bar{\boldsymbol{u}},\boldsymbol{v}-\bar{\boldsymbol{u}}) + \sum_{i=1}^{n-1} \left\langle g^{i}(x)\bar{\boldsymbol{\lambda}}^{i}_{T}, \left[ \boldsymbol{v}^{i}_{T} \right] - \left[ \bar{\boldsymbol{u}}^{i}_{T} \right] \right\rangle _{c} \geqslant L\left( \boldsymbol{v} - \bar{\boldsymbol{u}} \right).
\end{equation}
By setting $\mu^{i}_{HN} = \lambda^{i}_{HN}$ or $\boldsymbol{\mu}^{i}_{HT} = \boldsymbol{\lambda}^{i}_{HT}$ in problem \nameref{4:prb:mixhH}, we can easily derive the following relationship:
\begin{equation}\label{4:eq:proof:thm:convergence.7}
\left\langle g^{i}(x)\boldsymbol{\mu}^{i}_{HT}, \left[ \boldsymbol{u}^{i}_{HT} \right] \right\rangle _{c} \leqslant \left\langle g^{i}(x)\boldsymbol{\lambda}^{i}_{HT}, \left[ \boldsymbol{u}^{i}_{HT} \right] \right\rangle _{c}, ~~\forall \boldsymbol{\mu}^{i}_{HT} \in \Lambda^{i}_{HT}.
\end{equation}
According to the definitions of $\Lambda^{i}_{HT}$ and $\Lambda^{i}_{T}$, it can be verified that $\Lambda^{i}_{HT}$ is dense in $\Lambda^{i}_{T}$ with $H\to0$ and norm $\|\cdot\|_{0}$. So, a sequence $\{\boldsymbol{\mu}^{i}_{HT}\}$, $\boldsymbol{\mu}^{i}_{HT}\in\Lambda^{i}_{HT}$ can be constructed, which satisfies 
$$
\boldsymbol{\mu}^{i}_{HT} \to \frac{\left[ \bar{\boldsymbol{u}}^{i}_{T} \right]}{\left| \left[ \bar{\boldsymbol{u}}^{i}_{T} \right] \right|}.
$$
Therefore, based on the above formula and (\ref{4:eq:proof:thm:convergence.7}), it can be deduced that:
\begin{equation}\label{4:eq:proof:thm:convergence.8}
\int_{\Gamma^{i}_{c}} g^{i}(x) \left| \left[ \bar{\boldsymbol{u}}^{i}_{T} \right] \right| ds \leqslant \left\langle g^{i}(x)\bar{\boldsymbol{\lambda}}^{i}_{T}, \left[ \bar{\boldsymbol{u}}^{i}_{T} \right] \right\rangle _{c}. 
\end{equation}
Furthermore, according to the definition of $\Lambda^{i}_{T}$, it is easy to find the following inequality:
$$
\left\langle g^{i}(x)\bar{\boldsymbol{\lambda}}^{i}_{T}, \left[ \boldsymbol{v}^{i}_{T} \right] \right\rangle _{c} \leqslant \int_{\Gamma^{i}_{c}} g^{i}(x) \left| \left[ \boldsymbol{v}^{i}_{T} \right] \right| ds.
$$
By combining the above inequality and (\ref{4:eq:proof:thm:convergence.8}) with (\ref{4:eq:proof:thm:convergence.6}), the following variational inequality can be deduced:
$$
a\left(\bar{\boldsymbol{u}}, \boldsymbol{v} - \bar{\boldsymbol{u}}\right) + \sum_{i=1}^{n-1} \int_{\Gamma^{i}_{c}} g^{i}(x) \left( \left| \left[ \boldsymbol{v}^{i}_{T} \right] \right| - \left| \left[ \bar{\boldsymbol{u}}^{i}_{T} \right] \right| \right) ds \geqslant L(\boldsymbol{v} - \bar{\boldsymbol{u}}), ~\forall \boldsymbol{v} \in \mathcal{K}.
$$
Then, combined with $\bar{\boldsymbol{u}}\in \mathcal{K}$, it has been proved that $\bar{\boldsymbol{u}} = \boldsymbol{u}$. Referring to the inequality (\ref{4:eq:proof:thm:convergence.6}), it can also be verified that $\boldsymbol{\sigma}^{i}_{T} = -g^{i}(x) \bar{\boldsymbol{\lambda}}^{i}_{T}$ on $\Gamma^{i}_{c}$. 
Therefore, it can be concluded that the sequences $\left\{\boldsymbol{u}_{h}\right\}$ and $\left\{\boldsymbol{\lambda}^{i}_{HT}\right\}$ weakly converge to $\boldsymbol{u}$ and $\boldsymbol{\lambda}^{i}_{T}$, respectively, as the solution $\boldsymbol{u}$ of the problem \nameref{4:prb:var.ieq} is uniquely determined. 

It remains \textcolor{black}{to prove that} $\left\{\boldsymbol{u}_{h},\boldsymbol{\lambda}^{i}_{HT}\right\}$ \textcolor{black}{tends} strongly to $\left\{\boldsymbol{u},\boldsymbol{\lambda}^{i}_{T}\right\}$. Because the imbeding of $L^{2}\left(\Gamma^{i}_{c}\right)$ into $W^{i\prime}_{cT}$ is compact, it can be proved that $\boldsymbol{\lambda}^{i}_{TH}$ \textcolor{black}{tends} strongly to $\boldsymbol{\lambda}^{i}_{T}$ in $W^{i\prime}_{cT}$. In order to prove the strong convergence of $\left\{\boldsymbol{u}_{h}\right\}$, the following quadratic functional is defined:
\begin{equation}\label{4:eq:proof:thm:convergence.9}
F_{H} \left(\boldsymbol{v}\right) = \frac{1}{2} a(\boldsymbol{v},\boldsymbol{v}) - L(\boldsymbol{v}) + \sum_{i=1}^{n-1}\left\langle g^i(x) \boldsymbol{\lambda}_{HT}^{i},\left[\boldsymbol{v}_T^i\right]\right\rangle _c. 
\end{equation}
Just as the problems \nameref{4:prb:var.ieq} and \nameref{4:prb:optimization} are equivalent, it can be verified that the solution $\boldsymbol{u}_{h}$ to the variational inequality (\ref{4:ieq:var.hH}) is also the solution to the following optimization problem:
$$
\textcolor{black}{\boldsymbol{u}_{h} = \arg\min_{ \boldsymbol{v}_{h}\in \mathcal{K}_{hH}} F_{H}(\boldsymbol{v}_{h}).}
$$
Then, according to the Taylor expansion of $F_{H}(\boldsymbol{u}_{h})$ at $\boldsymbol{u}$, the following inequality can be deduced:
\begin{align*}
& F_{H}(\boldsymbol{v}_{h}) \geqslant F_{H}(\boldsymbol{u}_{h})\\
=& \frac{1}{2}a(\boldsymbol{u},\boldsymbol{u}) - L(\boldsymbol{u}) + \sum_{i=1}^{n-1}\left\langle g^i(x) \boldsymbol{\lambda}_{HT}^{i},\left[\boldsymbol{u}_T^i\right]\right\rangle _c + a(\boldsymbol{u},\boldsymbol{u}_{h} - \boldsymbol{u}) - L(\boldsymbol{u}_{h}-\boldsymbol{u}) \\ 
&+ \sum_{i=1}^{n-1}\left\langle g^i(x) \boldsymbol{\lambda}_{HT}^{i},\left[\boldsymbol{u}_{hT}^i\right] - \left[\boldsymbol{u}_{T}^i\right]\right\rangle _c + \frac{1}{2} a(\boldsymbol{u}_{h} - \boldsymbol{u},\boldsymbol{u}_{h} - \boldsymbol{u})\\
\geqslant& \frac{1}{2}a(\boldsymbol{u},\boldsymbol{u}) - L(\boldsymbol{u}) + a(\boldsymbol{u},\boldsymbol{u}_{h} - \boldsymbol{u}) - L(\boldsymbol{u}_{h}-\boldsymbol{u}) + \sum_{i=1}^{n-1}\left\langle g^i(x) \boldsymbol{\lambda}_{HT}^{i},\left[\boldsymbol{u}_{hT}^i\right]\right\rangle _c + \frac{a_{m}}{2} \left\| \boldsymbol{u}_{h}-\boldsymbol{u}\right\|^{2}_{V}\\
=& F_{H}(\boldsymbol{u}) + \sum_{i=1}^{n-1}\left\langle g^i(x) \boldsymbol{\lambda}_{HT}^{i},\left[\boldsymbol{u}_{hT}^i\right] - \left[\boldsymbol{u}_{T}^i\right]\right\rangle _c + a(\boldsymbol{u},\boldsymbol{u}_{h} - \boldsymbol{u}) - L(\boldsymbol{u}_{h}-\boldsymbol{u}) + \frac{a_{m}}{2} \left\| \boldsymbol{u}_{h}-\boldsymbol{u}\right\|^{2}_{V}
\end{align*}
where $\forall \boldsymbol{v}_{h}\in \mathcal{K}_{hH}$ and $a_{m}>0$ decided by the elastic operators $A^{i}$. By taking $\boldsymbol{v}_{h}\to \boldsymbol{u}$, based on the continuity of $F_{H}$, $\boldsymbol{u}_{h}\rightharpoonup \boldsymbol{u}$ and $\boldsymbol{\lambda}^{i}_{HT}\rightharpoonup \boldsymbol{\lambda}^{i}_{T}$, it can be found that:
\begin{align*}
\lim_{h,H\to 0+}\left\| \boldsymbol{u}_{h}-\boldsymbol{u}\right\|^{2}_{V} \leqslant & \lim_{h,H\to 0+}\frac{2}{a_{m}} \Big( F_{H}(\boldsymbol{v}_{h}) - F_{H}(\boldsymbol{u}) - \sum_{i=1}^{n-1}\left\langle g^i(x) \boldsymbol{\lambda}_{HT}^{i},\left[\boldsymbol{u}_{hT}^i\right] - \left[\boldsymbol{u}_{T}^i\right]\right\rangle _c\\ 
&- a(\boldsymbol{u},\boldsymbol{u}_{h} - \boldsymbol{u}) + L(\boldsymbol{u}_{h}-\boldsymbol{u}) \Big) \to 0.
\end{align*}
\end{proof}

\subsection{Error estimates}

In this subsection, the distance between the finite element numerical solution $\boldsymbol{u}_{h}$ and the exact solution $\boldsymbol{u}$ will be estimated. In the problems \nameref{4:prb:mix} and \nameref{4:prb:mixhH}, the stress on the contact zone $\Gamma^{i}_{c}$ will also be calculated. Therefore, the distance between the exact stress field $\boldsymbol{\lambda}$ and the approximate stress field $\boldsymbol{\lambda}_{H}$ can be estimated as well.

In order to deduce the error relation, an equivalent variational problem of \nameref{4:prb:mix} should be provided. Firstly, the Hilbert space $M=V\times W^{\prime}_{c}$ and the corresponding norm is defined by
$$
\|\mathfrak{v}\|_{M} = \left( \|\boldsymbol{v}\|_{V}^2 + \|\boldsymbol{\mu}\|_{W^{\prime}_{c}}^{2} \right)^{1/2},~~\mathfrak{v} = \left( \boldsymbol{v}, \boldsymbol{\mu}\right) \in M
$$
Then, the bilinear operator $\mathfrak{H}$ and the linear functional $\mathfrak{L}$ are defined as follows:
\begin{align}
& \mathfrak{H}(\mathfrak{u}, \mathfrak{v}) = a(\boldsymbol{u},\boldsymbol{v}) + \mathcal{G}(\boldsymbol{\lambda},\boldsymbol{v}) - \mathcal{G}(\boldsymbol{\mu},\boldsymbol{u}),~\mathfrak{u}=\left(\boldsymbol{u},\boldsymbol{\lambda}\right),~\mathfrak{v}=\left(\boldsymbol{v},\boldsymbol{\mu}\right) \in M, \label{4:eq:fH}\\
& \mathfrak{L}(\mathfrak{v}) = L(\boldsymbol{v}),~~\mathfrak{v}=\left(\boldsymbol{v},\boldsymbol{\mu}\right) \in M. \label{4:eq:fL}
\end{align}
Based on the definition of $\mathfrak{H}$ in formula (\ref{4:eq:fH}), it can be found that:
\begin{equation}\label{4:eq:afHv}
a_{m}\|\boldsymbol{v}\|_{V}^{2} \leqslant a(\boldsymbol{v},\boldsymbol{v}) = \mathfrak{H}(\mathfrak{v},\mathfrak{v}),~\forall \mathfrak{v} \in M,
\end{equation}
where the constant $a_{m}$ is determined by the elastic operator $A^{i}$, that is $a_{m}=\min_{1\leqslant i \leqslant n} \{a^{i}_{m}\}$. Moreover, based on the properties of the elastic operator $A^{i}$ in formula (\ref{4:character:A}), the inequality (\ref{4:eq:trace}), the definitions of the operators $g$ in formula (\ref{4:eq:def:g}) and $\mathfrak{H}$ in formula (\ref{4:eq:fH}), there exists a positive constant $c_{M}$ such that
\begin{equation}\label{4:ieq:fH.cM}
\left|\mathfrak{H}(\mathfrak{u}, \mathfrak{v}) \right| \leqslant c_{M} \|\mathfrak{u}\|_{M}\cdot \|\mathfrak{v}\|_{M},~ \forall \mathfrak{u}, \mathfrak{v} \in M.
\end{equation}
Then, it can be verified that $\textcolor{black}{\mathfrak{J}=V \times \Lambda}$ is the closed and convex subspace of $M$. Therefore, the following optimization problem can be introduced:

\begin{problem}[$\mathcal{P}_{mv}$]\label{4:prb:mv}
Find $\mathfrak{u} = \left( \boldsymbol{u}, \boldsymbol{\lambda}\right) \in \textcolor{black}{\mathfrak{J}}$ such that
\begin{equation}\label{4:eq:mv}
\mathfrak{H}(\mathfrak{u}, \mathfrak{v} - \mathfrak{u}) \geqslant \mathfrak{L}(\mathfrak{v}-\mathfrak{u}),~~\forall \mathfrak{v}\in \textcolor{black}{\mathfrak{J}}.
\end{equation}
\end{problem}

The problem \nameref{4:prb:mix} and problem \nameref{4:prb:mv} are equivalent. Indeed, the variational inequality (\ref{4:eq:mv}) can be expressed as:
$$
a(\boldsymbol{u},\boldsymbol{v}-\boldsymbol{u}) + \mathcal{G}(\boldsymbol{\lambda},\boldsymbol{v}-\boldsymbol{u}) - \mathcal{G}(\boldsymbol{\mu}-\boldsymbol{\lambda},\boldsymbol{u}) \geqslant L(\boldsymbol{v}-\boldsymbol{u}),
$$
which is equivalent with the mix problem (\ref{4:eq:mix}). Naturally, the finite element approximation problem corresponding to the \textcolor{black}{problem} \nameref{4:prb:mv} can be defined as follows:
\begin{problem}[$\mathcal{P}_{vhH}$]\label{4:prb:vhH}
Find $\mathfrak{u}_{hH} = \left( \boldsymbol{u}_{h}, \boldsymbol{\lambda}_{H}\right) \in \textcolor{black}{\mathfrak{J}}_{hH}$ such that
\begin{equation}\label{4:eq:vhH}
\mathfrak{H}(\mathfrak{u}_{hH}, \mathfrak{v}_{hH} - \mathfrak{u}_{hH}) \geqslant \mathfrak{L}(\mathfrak{v}_{hH} - \mathfrak{u}_{hH}),~~\forall \mathfrak{v}_{hH}\in \textcolor{black}{\mathfrak{J}}_{hH},
\end{equation}
where $\textcolor{black}{\mathfrak{J}}_{hH} = V_{h}\times \Lambda_{H}$ is the subset of $\textcolor{black}{\mathfrak{J}}$.
\end{problem}

Similarly, it can be proved that discrete problems \nameref{4:prb:mixhH} and \nameref{4:prb:vhH} are equivalent. Then, the solutions of problems \nameref{4:prb:mv} and \nameref{4:prb:vhH} satisfy the following lemma:

\begin{lemma}[]\label{4:lem:mv.vhH}
The solutions of problems \nameref{4:prb:mv} and \nameref{4:prb:vhH} are denoted as $\mathfrak{u} = (\boldsymbol{u}, \boldsymbol{\lambda})$ and $\mathfrak{u}_{hH}= (\boldsymbol{u}_{h}, \boldsymbol{\lambda}_{H})$, respectively. Then, the following inequality holds:
\begin{equation}\label{4:ieq:u.uh}
\|\boldsymbol{u} - \boldsymbol{u}_{h}\|^{2}_{V} \leqslant \frac{1}{a_{m}} \left( \mathfrak{H}(\mathfrak{u}-\mathfrak{u}_{hH}, \mathfrak{u}-\mathfrak{v}_{hH}) + \mathfrak{H}(\mathfrak{u},\mathfrak{v}_{hH} - \mathfrak{u}) + \mathfrak{L}(\mathfrak{u}-\mathfrak{v}_{hH})\right), ~\forall \mathfrak{v}_{hH}\in \textcolor{black}{\mathfrak{J}}_{hH}.
\end{equation} 
\end{lemma}
\begin{proof}
Based on formula (\ref{4:eq:afHv}) and the definitions of \nameref{4:prb:mv} and \nameref{4:prb:vhH}, it can be deduced that
\begin{align*}
&a_{m} \left\| \boldsymbol{u}-\boldsymbol{u}_{h} \right\|_{V}^{2} \leqslant \mathfrak{H}\left(\mathfrak{u}-\mathfrak{u}_{hH}, \mathfrak{u}-\mathfrak{u}_{hH} \right) \\
=& \mathfrak{H}\left(\mathfrak{u}-\mathfrak{u}_{hH}, \mathfrak{u}-\mathfrak{v}_{hH} \right) + \mathfrak{H}\left(\mathfrak{u}, \mathfrak{v}_{hH}-\mathfrak{u} \right) - \mathfrak{H}\left(\mathfrak{u}, \mathfrak{u}_{hH}- \mathfrak{u}\right) - \mathfrak{H}\left(\mathfrak{u}_{hH}, \mathfrak{v}_{hH}-\mathfrak{u}_{hH} \right)\\
\leqslant& \mathfrak{H}\left(\mathfrak{u}-\mathfrak{u}_{hH}, \mathfrak{u}-\mathfrak{v}_{hH} \right) + \mathfrak{H}\left(\mathfrak{u}, \mathfrak{v}_{hH}-\mathfrak{u} \right) - \mathfrak{L}(\mathfrak{u}_{hH}- \mathfrak{u}) - \mathfrak{L}(\mathfrak{v}_{hH}-\mathfrak{u}_{hH})\\
=& \mathfrak{H}\left(\mathfrak{u}-\mathfrak{u}_{hH}, \mathfrak{u}-\mathfrak{v}_{hH} \right) + \mathfrak{H}\left(\mathfrak{u}, \mathfrak{v}_{hH}-\mathfrak{u} \right) + \mathfrak{L}(\mathfrak{u}-\mathfrak{v}_{hH}),~~\forall \mathfrak{v}_{hH}\in \textcolor{black}{\mathfrak{J}}_{hH}.
\end{align*}
\end{proof}

The premise for obtaining the error analysis of the distance between the exact solution $\boldsymbol{\lambda}$ and the finite element solution $\boldsymbol{\lambda}_{H}$ is that the solution $\boldsymbol{\lambda}_{H}$ in problem \nameref{4:prb:vhH} is unique.
\textcolor{black}{Therefore, under the assumption that condition (\ref{4:eq:lambdaH.unique}) holds, the following lemma can be derived:}

\begin{lemma}[]\label{4:lem:lambda.lambdaH}
If the condition (\ref{4:eq:lambdaH.unique}) holds, then the following inequality can be derived:
\begin{equation}\label{4:ieq:lambdaH.Wc}
\left\|\boldsymbol{\lambda}-\boldsymbol{\lambda}_{H} \right\|_{W^{\prime}_{c}} \leqslant c_{p} \left( \left\|\boldsymbol{u}-\boldsymbol{u}_{h} \right\|_{V} + \inf_{\boldsymbol{\mu}_{H}\in \Lambda_{H}} \left\|\boldsymbol{\lambda}-\boldsymbol{\mu}_{H} \right\|_{W^{\prime}_{c}} \right).
\end{equation}
\end{lemma}
\begin{proof}
For any $\boldsymbol{\mu}_{H} = \left(\boldsymbol{\mu}^{i}_{H}, \ldots, \boldsymbol{\mu}^{n-1}_{H}\right)$, based on the problems \nameref{4:prb:mix} and \nameref{4:prb:mixhH}, it can be deduced that:
\begin{align*}
&\mathcal{G}(\boldsymbol{\mu}_{H}-\boldsymbol{\lambda}_{H},\boldsymbol{v}_{h}) = a(\boldsymbol{u}_{h},\boldsymbol{v}_{h}) - L(\boldsymbol{v}_{h}) + \mathcal{G}(\boldsymbol{\mu}_{H},\boldsymbol{v}_{h})\\
=& \mathcal{G}(\boldsymbol{\mu}_{H},\boldsymbol{v}_{h}) + a(\boldsymbol{u}_{h}-\boldsymbol{u},\boldsymbol{v}_{h}) + a(\boldsymbol{u},\boldsymbol{v}_{h}) - L(\boldsymbol{v}_{h})\\
=& a(\boldsymbol{u}_{h}-\boldsymbol{u},\boldsymbol{v}_{h}) + \mathcal{G}(\boldsymbol{\mu}_{H},\boldsymbol{v}_{h}) - \mathcal{G}(\boldsymbol{\lambda},\boldsymbol{v}_{h})\\
=& a(\boldsymbol{u}_{h}-\boldsymbol{u},\boldsymbol{v}_{h}) + \sum_{i=1}^{n-1} \left(\left\langle \mu^{i}_{HN} - \lambda^{i}_{N}, [v^{i}_{hN}] \right\rangle _{c} + \left\langle g^{i}({x})\left(\boldsymbol{\mu}^{i}_{HT} - \boldsymbol{\lambda}^{i}_{T}\right), [\boldsymbol{v}^{i}_{hT}] \right\rangle _{c} \right)\\
\leqslant& a_{M} \|\boldsymbol{u}_{h}-\boldsymbol{u}\|_{V}\cdot\|\boldsymbol{v}_{h}\|_{V} + \sqrt{2}c_{g} \|\boldsymbol{\lambda} - \boldsymbol{\mu}_{H}\|_{W^{\prime}_{c}}\cdot\sum_{i=1}^{n-1}\left( \|[v^{i}_{hN}]\|_{1/2,\Gamma^{i}_{c}} + \|[\boldsymbol{v}^{i}_{hT}]\|_{1/2,\Gamma^{i}_{c}} \right)\\
\leqslant& a_{M} \|\boldsymbol{u}_{h}-\boldsymbol{u}\|_{V}\cdot\|\boldsymbol{v}_{h}\|_{V} + 4nc_{t}c_{g}\|\boldsymbol{\lambda} - \boldsymbol{\mu}_{H}\|_{W^{\prime}_{c}}\cdot\left\|\boldsymbol{v}_{h} \right\|_{V}\\
\leqslant& c_{l} \left( \|\boldsymbol{u}_{h}-\boldsymbol{u}\|_{V} + \|\boldsymbol{\lambda} - \boldsymbol{\mu}_{H}\|_{W^{\prime}_{c}} \right) \cdot\left\|\boldsymbol{v}_{h} \right\|_{V},
\end{align*}
where constants $c_{l} = \max\{a_{M}, 4nc_{t}c_{g}\}$, $a_{M} = \max_{1\leqslant i \leqslant n}\{a^{i}_{M}\}$, $c_{t} = \max_{1\leqslant i \leqslant n}\{c^{i}_{t}\}$ and $c_{g} = \max_{1\leqslant i \leqslant n}\{g^{i}(x),1\}$. Then, according to the condition (\ref{4:eq:lambdaH.unique}), it can be obtained that:
$$
c_{hH} \left\|\boldsymbol{\lambda}_{H} - \boldsymbol{\mu}_{H}\right\|_{W'_{c}} \leqslant\sup _{V_h} \frac{\mathcal{G}\left(\boldsymbol{\lambda}_{H} - \boldsymbol{\mu}_{H}, \boldsymbol{v}_h\right)}{\left\|\boldsymbol{v}_h\right\|_{V}} \leqslant c_{l} \left( \|\boldsymbol{u}_{h}-\boldsymbol{u}\|_{V} + \|\boldsymbol{\lambda} - \boldsymbol{\mu}_{H}\|_{W^{\prime}_{c}} \right).
$$
Finally, based on the triangle inequality, it can be verified that
\begin{align*}
\|\boldsymbol{\lambda} - \boldsymbol{\lambda}_{H}\|_{W^{\prime}_{c}} \leqslant& \|\boldsymbol{\lambda} - \boldsymbol{\mu}_{H}\|_{W^{\prime}_{c}} + \|\boldsymbol{\mu}_{H} - \boldsymbol{\lambda}_{H}\|_{W^{\prime}_{c}}\\
\leqslant & \frac{c_{l}}{c_{hH}} \left( \|\boldsymbol{u}_{h}-\boldsymbol{u}\|_{V} + \|\boldsymbol{\lambda} - \boldsymbol{\mu}_{H}\|_{W^{\prime}_{c}} \right) + \|\boldsymbol{\lambda} - \boldsymbol{\mu}_{H}\|_{W^{\prime}_{c}}\\
\leqslant & c_{p} \left( \left\|\boldsymbol{u}-\boldsymbol{u}_{h} \right\|_{V} + \inf_{\boldsymbol{\mu}_{H}\in \Lambda_{H}} \left\|\boldsymbol{\lambda}-\boldsymbol{\mu}_{H} \right\|_{W^{\prime}_{c}} \right).
\end{align*}
\end{proof}

Based on the above lemmas, the convergence rate of $\mathfrak{u}_{hH}=(\boldsymbol{u}_{h},\boldsymbol{\lambda}_{H})$ can be discussed. By utilizing the inequalities $2ab\leqslant \epsilon a^{2} +1/\epsilon b^{2}$ and formula (\ref{4:ieq:fH.cM}), formula (\ref{4:ieq:u.uh}) can be expressed as
\begin{equation}\label{4:ieq:u.uh.1}
\begin{aligned}
\|\boldsymbol{u} - \boldsymbol{u}_{h}\|^{2}_{V} \leqslant& \frac{c_{M}}{a_{m}} \Big( \epsilon\|\boldsymbol{u} -\boldsymbol{u}_{h} \|^{2}_{V} + \epsilon\| \boldsymbol{\lambda} - \boldsymbol{\lambda}_{H} \|^{2}_{W^{\prime}_{c}} + 1/\epsilon\|\boldsymbol{u} -\boldsymbol{v}_{h} \|^{2}_{V} + 1/\epsilon\| \boldsymbol{\lambda} - \boldsymbol{\mu}_{H} \|^{2}_{W^{\prime}_{c}} + \\
&a(\boldsymbol{u},\boldsymbol{v}_{h} - \boldsymbol{u}) + \mathcal{G}(\boldsymbol{\lambda}, \boldsymbol{v}_{h} - \boldsymbol{u}) - \mathcal{G}(\boldsymbol{\mu}_{H}- \boldsymbol{\lambda}, \boldsymbol{u}) + L(\boldsymbol{u}-\boldsymbol{v}_{h}) 
\Big).
\end{aligned}
\end{equation}
Then, based on \textcolor{black}{Green's formula} (\ref{4:eq:Green}), it can be deduced that:
\begin{align*}
a(\boldsymbol{u},\boldsymbol{v}_{h} - \boldsymbol{u}) =& \sum_{i=1}^{n}\left( \left\langle \boldsymbol{\sigma}^{i}\cdot\boldsymbol{\alpha}^{i}, \gamma^{i}_{2}(\boldsymbol{v}^{i}_{h} - \boldsymbol{u}^{i}) \right\rangle _2 + \left\langle \boldsymbol{\sigma}^{i}\cdot\boldsymbol{\beta}^{i}, \gamma^{i}_{3}(\boldsymbol{v}^{i}_{h} - \boldsymbol{u}^{i}) \right\rangle _3\right) - \sum_{i=1}^{n} \left( \nabla\cdot\boldsymbol{\sigma}^{i}, \boldsymbol{v}^{i}_{h} - \boldsymbol{u}^{i} \right)_{0}\\
=& \sum_{i=1}^{n} \left( \nabla\cdot\boldsymbol{\sigma}^{i}, \boldsymbol{u}^{i} - \boldsymbol{v}^{i}_{h}\right)_{0} - \left\langle \boldsymbol{\sigma}^{1}\cdot\boldsymbol{\alpha}^{1}, \gamma^{1}_{2}(\boldsymbol{u}^{1} -\boldsymbol{v}^{1}_{h}) \right\rangle _2\\ 
&+ \sum_{i=1}^{n-1} \left( \left\langle {\sigma}^{i}_{N}, [{v}^{i}_{hN}] - [{u}^{i}_{N}] \right\rangle _c + \left\langle \boldsymbol{\sigma}^{i}_{T} , [\boldsymbol{v}^{i}_{hT}] - [\boldsymbol{u}^{i}_{T}] \right\rangle _c \right).
\end{align*}
According to the stress boundary condition and the balanced function of body force, it can be inferred that
$$
\nabla\cdot\boldsymbol{\sigma}^{i} = -\boldsymbol{f}^{i}_{0} \text{ a.e. on } \Omega^{i} \text{ and } \boldsymbol{\sigma}^{1}\cdot\boldsymbol{\alpha}^{1} = \boldsymbol{f}_{1} \text{ a.e. on } \Gamma^{1}_{2}.
$$
Therefore, based on the above relation, the definition of the linear operator $L$ in problem \nameref{4:prb:var.ieq} and \textcolor{black}{Therorem} \ref{4:thm:saddle.point}, the inequality (\ref{4:ieq:u.uh.1}) can be expressed as
\begin{equation}\label{4:ieq:u.uh.2}
\begin{aligned}
\|\boldsymbol{u} - \boldsymbol{u}_{h}\|^{2}_{V} \leqslant& \frac{c_{M}}{a_{m}} \Big( \epsilon\|\boldsymbol{u} -\boldsymbol{u}_{h} \|^{2}_{V} + \epsilon\| \boldsymbol{\lambda} - \boldsymbol{\lambda}_{H} \|^{2}_{W^{\prime}_{c}} + 1/\epsilon\|\boldsymbol{u} -\boldsymbol{v}_{h} \|^{2}_{V} + 1/\epsilon\| \boldsymbol{\lambda} - \boldsymbol{\mu}_{H} \|^{2}_{W^{\prime}_{c}} + \\
& - \mathcal{G}(\boldsymbol{\mu}_{H}- \boldsymbol{\lambda}, \boldsymbol{u}) 
\Big).
\end{aligned}
\end{equation}
Then, if the constant $\epsilon > 0$ is sufficiently small, by combining inequalities (\ref{4:ieq:u.uh.2}) and (\ref{4:ieq:lambdaH.Wc}), the following inequalities can be deduced:
\begin{equation}\label{4:ieq:u.uh.3}
\|\boldsymbol{u} - \boldsymbol{u}_{h}\|^{2}_{V} \leqslant c_{\epsilon} \inf_{V_{h}\times \Lambda_{H}} \left( \|\boldsymbol{u} -\boldsymbol{v}_{h} \|^{2}_{V} + \| \boldsymbol{\lambda} - \boldsymbol{\mu}_{H} \|^{2}_{W^{\prime}_{c}} - \mathcal{G}(\boldsymbol{\mu}_{H}- \boldsymbol{\lambda}, \boldsymbol{u}) \right),
\end{equation}
where the constant $c_{\epsilon}>0$ is dependent on $c_{M}$, $a_{m}$, $c_{p}$ and $\epsilon$.

Finally, according to the inequality (\ref{4:ieq:u.uh.3}), the convergence rate of the finite element approximation solution $\left(\boldsymbol{u}_{h}, \boldsymbol{\lambda}_{H}\right)$ to the exact solution $\left(\boldsymbol{u}, \boldsymbol{\lambda}\right)$ can be estimated, that is the following theorem holds.

\begin{theorem}[]\label{4:thm:error.estimate}
Assume that the inequality (\ref{4:eq:lambdaH.unique}) holds and $\boldsymbol{u}^{i}\in \mathcal{H}^{2}(\Omega^{i})^{d}$, $\sigma^{i}_{N}\in \mathcal{H}^{3/2}(\Gamma^{i}_{c})$, $\boldsymbol{\sigma}^{i}_{T}\in \mathcal{H}^{3/2}(\Gamma^{i}_{c})^d$, $g^{i}(x) \in C^{\infty}(\Gamma^{i}_{c})$. Then, the following error estimation formulas are established:
\begin{equation}\label{4:eq:u.uh.error}
\| \boldsymbol{u} - \boldsymbol{u}_{h}\|_{V} = O(H^{\hat{k}}),~ H\to 0+
\end{equation}
and
\begin{equation}\label{4:eq:lambda.lambdaH.error}
\| \boldsymbol{\lambda} - \boldsymbol{\lambda}_{H}\|_{W^{\prime}_{c}} = O(H^{\hat{k}}),~ H\to 0+,
\end{equation}
where $\hat{k} = 3/4$.
\end{theorem}

\begin{proof}
Firstly, the piecewise Lagrange interpolation operator is denoted as $\Xi_{h}$. According to the properties of the operator $\Xi_{h}$ \cite{ciarlet2002finite,chouly2023finite}, it can be known that $\Xi_{h} \boldsymbol{u}\in V_{h}$ and
\begin{equation*}
\|\boldsymbol{u}-\Xi_{h} \boldsymbol{u} \|_{V} \leqslant ch \|\boldsymbol{u}\|_{2}.
\end{equation*}
According to the assumptions that $\lambda^{i}_{N} = -\sigma^{i}_{N} \in \mathcal{H}^{3/2}(\Gamma^{i}_{c})$ and  $g^{i}(x)\boldsymbol{\lambda}^{i}_{T} = -\boldsymbol{\sigma}^{i}_{T} \in \mathcal{H}^{3/2}(\Gamma^{i}_{c})^{d}$, it can be deduced that 
$$
\left(\varrho^{i*}_{c}\right)^{-1} \boldsymbol{\lambda}^{i} = \boldsymbol{\phi}^{i} \in \mathcal{H}^{3/2}(\Gamma^{i}_{c})^d.
$$
Then, the approximation solution $\boldsymbol{\phi}^{i}_{H} \in \mathcal{H}^{-1/2}(\Gamma^{i}_{c})^d$ is defined by $\boldsymbol{\phi}^{i}_{H}= \Xi_{H} \boldsymbol{\phi}^{i}$. Since $\lambda^{i}_{N}\geqslant 0$ and $|\boldsymbol{\lambda}^{i}_{T}| \leqslant 1$, these conditions are satisfied for $\varrho^{i*}_{c} \boldsymbol{\phi}^{i}_{H}$. Therefore, it can be deduced that $\varrho^{i*}_{c} \boldsymbol{\phi}^{i}_{H} \in \Lambda^{i}_{H}$ and 
$$
\|\boldsymbol{\lambda}^{i} - \varrho^{i*}_{c} \boldsymbol{\phi}^{i}_{H}\|_{W^{i\prime}_{c}} = \| \left(\varrho^{i*}_{c}\right)^{-1}\boldsymbol{\lambda}^{i} -  \boldsymbol{\phi}^{i}_{H}\|_{-1/2,\Gamma^{i}_{c}} \leqslant ch^{3/2} \|\boldsymbol{\lambda}^{i}\|_{3/2,\Gamma^{i}_{c}}.
$$
Based on the above condition, when $\boldsymbol{\mu}_{H} = \left(\varrho^{1*}_{c} \boldsymbol{\phi}^{1}_{H},\ldots,\varrho^{n-1*}_{c} \boldsymbol{\phi}^{n-1}_{H}\right)$, the following inequality can be verified:
\begin{align*}
\mathcal{G}(\boldsymbol{\mu}_{H}- \boldsymbol{\lambda}, \boldsymbol{u}) &= \sum_{i=1}^{n-1} \left( \left\langle \phi^{i}_{HN} - \lambda^{i}_{N}, [u^{i}_{N}] \right\rangle _{c} + \left\langle g^{i}(x)\left(\boldsymbol{\phi}^{i}_{HT} - \boldsymbol{\lambda}^{i}_{T}\right), [\boldsymbol{u}^{i}_{T}] \right\rangle _{c} \right)\\
&\leqslant \sum_{i=1}^{n-1} \left( \left\|\phi^{i}_{HN} - \lambda^{i}_{N} \right\|_{W^{i\prime}_{N}} \cdot \left\| [u^{i}_{N}] \right\|_{W^{i}_{N}} + c_{g}\left\|\boldsymbol{\phi}^{i}_{HT} - \boldsymbol{\lambda}^{i}_{T} \right\|_{W^{i\prime}_{T}} \cdot \left\| [\boldsymbol{u}^{i}_{T}] \right\|_{W^{i}_{T}}\right) \\
&\leqslant c\sum_{i=1}^{n-1} \left\| \varrho^{i*}_{c}\boldsymbol{\phi}^{i}_{H} - \boldsymbol{\lambda}^{i} \right\|_{W^{i\prime}_{c}} \cdot \left\| \boldsymbol{u} \right\|_{V} \leqslant ch^{3/2} \left\| \boldsymbol{u} \right\|_{V}\cdot \sum_{i=1}^{n-1} \|\boldsymbol{\lambda}^{i}\|_{3/2,\Gamma^{i}_{c}}.
\end{align*}
Hence, based on the inequality (\ref{4:ieq:u.uh.3}), the error estimation formula (\ref{4:eq:u.uh.error}) of $\boldsymbol{u}$ is derived. Then, according to (\ref{4:ieq:lambdaH.Wc}), the error estimation formula (\ref{4:eq:lambda.lambdaH.error}) of $\boldsymbol{\lambda}$ can be deduced.
\end{proof}

\begin{remark}[]\label{4:rem:error.regularity}
According to the assumptions and conclusions of Theorem \ref{4:thm:error.estimate}, it can be found that when the regularity of the exact solution is strong, the convergence rate of the numerical solution is faster. When the assumptions of the displacement field and stress field on the contact zone $\Gamma^{i}_{c}$ in Theorem \ref{4:thm:error.estimate} are changed to $\boldsymbol{u}^{i}\in \mathcal{H}^{1+k}(\Omega^{i})^{d}$ and $\sigma^{i}_{N}\in L^{2}(\Gamma^{i}_{c})$, $g^{i}(x) \in L^{\infty}(\Gamma^{i}_{c})$ respectively, it can be obtained that $\left(\varrho^{i*}_{c}\right)^{-1} \boldsymbol{\lambda}^{i} = \boldsymbol{\phi}^{i} \in L^{2}(\Gamma^{i}_{c})^d$. Then, by defining $\Lambda^{i}_{H}$ as a piecewise constant space and $\boldsymbol{\phi}^{i}_{H}$ as an approximate solution satisfying that $\boldsymbol{\phi}^{i}_{H}|_{S^{i}_{l}}$ is the $L^{2}(S^{i}_{l})$-projection of $\boldsymbol{\phi}^{i}|_{S^{i}_{l}}$ onto $P_{0}(S^{i}_{l})^{d}$, $S^{i}_{l}\in F^{i}_{H}$, the convergence rate estimate can be obtained in the $L^{2}(\Gamma^{i}_{c})$ space, and $\hat{k}=\min\{k,1/4\}$.
\end{remark}

\section{Implementation of Algorithm}

In this section, the algebraic form of the mixed finite element method for a multi-layer elastic system will be introduced to ensure that \textcolor{black}{the algorithm} is executable on the computer.

First, to simplify the construction of algebraic forms, the finite element spaces corresponding to the function spaces $V^i$ and $\Lambda^{i}$ are defined as follows:
\begin{align}
& V_h^i=\left\{\boldsymbol{v}_h^i \in C\left(\bar{\Omega}^i\right)^d ~\Big|~ \boldsymbol{v}_h^i|_{T_j^i} \in P_{1}\left(T_j^i\right)^d, \forall T^{i}_{j}\in E^{i}_{h}, \boldsymbol{v}_h^i=0 \text { on } \Gamma_1^i\right\}, \label{4:space:Vih.P1}\\
& \Lambda^{i*}_{HN}=\left\{\mu^{i}_{hN} \in L^2\left(\Gamma^{i}_{c}\right) ~\Big|~ \mu^{i}_{hN}|_{S^{i}_{j}} \in P_{0}\left(S^{i}_{j}\right), \mu^{i}_{hN} \geqslant 0 \text { on } \Gamma^{i}_{c}, \forall S^{i}_{j} \in F^{i}_{H}\right\}, \label{4:space:LiNH.P0}\\
& \Lambda^{i*}_{HT}=\left\{\boldsymbol{\mu}^{i}_{hT} \in L^2\left(\Gamma^{i}_{c}\right)^d ~\Big|~ \boldsymbol{\mu}^{i}_{hT}|_{S^{i}_{j}} \in P_{0} \left(S^{i}_{j}\right)^d, \left|\boldsymbol{\mu}^{i}_{hT}\right| \leqslant g^{i}(x) \text { on } \Gamma^{i}_{c}, \forall S^{i}_{j} \in F^{i}_{H}\right\}, \label{4:space:LiTH.P0}\\
& \Lambda^{i**}_{HN}=\left\{\mu^{i}_{hN} \in C\left(\bar{\Gamma}^{i}_{c}\right) ~\Big|~ \mu^{i}_{hN}|_{S^{i}_{j}} \in P_{1}\left(S^{i}_{j}\right), \mu^{i}_{hN} \geqslant 0 \text { on } \Gamma^{i}_{c}, \forall S^{i}_{j} \in F^{i}_{H} \right\}, \label{4:space:LiNH.P1}\\
& \Lambda^{i**}_{HT}=\left\{\boldsymbol{\mu}^{i}_{hT} \in C\left(\bar{\Gamma}^{i}_{c}\right)^d ~\Big|~ \boldsymbol{\mu}^{i}_{hT}|_{S^{i}_{j}} \in P_{1} \left(S^{i}_{j}\right)^d, \left|\boldsymbol{\mu}^{i}_{hT}\right| \leqslant g^{i}(x) \text { on } \Gamma^{i}_{c}, \forall S^{i}_{j} \in F^{i}_{H}\right\}, \label{4:space:LiTH.P1}
\end{align}
where $h=H$ is specified in the subsequent algorithm. Then $V_{h}$, $\Lambda^{i*}_{H}$, $\Lambda^{*}_{H}$, $\Lambda^{i**}_{H}$, $\Lambda^{**}_{H}$ are defined similarly to (\ref{4:space:VhWH}) and (\ref{4:space:LH}). Numerical experiments will be performed in both finite element spaces $\Lambda^{*}_{H}$ and $\Lambda^{**}_{H}$. To simplify the notation, both are denoted as $\Lambda_{H}$ and will be explained only when distinction is needed. Based on the definition of the above finite element spaces, the discrete form of the mixed finite element problem \nameref{4:prb:mixhH} can be rewritten as

\begin{problem}[$\mathcal{P}^{*}_{hH}$]\label{4:prb:mixhH*}
Find a solution $\{\boldsymbol{u}_{h}, \boldsymbol{\lambda}_{H}\} \in V_{h} \times \Lambda_{h}$ which satisfies:
\begin{equation}\label{4:eq:mixhH*}
\left\{\begin{aligned}
& a(\boldsymbol{u}_{h}, \boldsymbol{v}_{h})+\mathcal{G}^{*}(\boldsymbol{\lambda}_{H}, \boldsymbol{v}_{h})=L(\boldsymbol{v}_{h}), \quad \forall \boldsymbol{v}_{h} \in V_{h} \\
& \mathcal{G}^{*}(\boldsymbol{\mu}_{H}-\boldsymbol{\lambda}_{H}, \boldsymbol{u}_{h}) \leqslant 0, \quad \forall \boldsymbol{\mu}_{H} \in \Lambda_{H}
\end{aligned}\right. ,
\end{equation}
where $\mathcal{G}^{*}(\boldsymbol{\mu}_{H}, \boldsymbol{v}_{h})$ is defined as
\begin{equation}\label{4:eq:def:g*}
\mathcal{G}^{*}(\boldsymbol{\mu}_{H}, \boldsymbol{v}_{H}) = \sum_{i=1}^{n-1} \left(\left\langle \mu^{i}_{NH}, [v^{i}_{NH}] \right\rangle_{c} + \left\langle \boldsymbol{\mu}^{i}_{TH}, [\boldsymbol{v}^{i}_{TH}] \right\rangle_{c} \right).
\end{equation}
\end{problem}

In the finite element meshing, the number of finite elements and the number of finite element nodes on the elastic layer $\bar{\Omega}^{i}$ are recorded as $N^{i}_{h}$ and $N^{i}_v$ respectively; the number of finite elements and the number of finite element nodes on the contact zone $\Gamma^{i}_{c}$ are recorded as $N^{i}_{H}$ and $N^{i}_c$ respectively. Then, the basis function sets of the finite element space $V^{i}_{h}$, $\Lambda^{i}_{HN}$ and $\Lambda^{i}_{HT}$ are denoted as $\{\phi^{i}_{k}\}_{k=1}^{dN^{i}_{v}}$, $\{\psi^{i}_{Nk}\}_{k=1}^{N^{i}_{l}}$, $\{\psi^{i}_{Tk}\}_{k=1}^{dN^{i}_{l}}$, where $N^{i}_{l} = N^{i}_{H}$ if $\Lambda_{H} = \Lambda^{*}_{H}$ and $N^{i}_{l} = N^{i}_{c}$ if $\Lambda_{H} = \Lambda^{**}_{H}$. The diagonal basis function matrices composed of basis functions are recorded as $\Phi^{i}$, $\Psi^{i}_N$ and $\Psi^{i}_T$ respectively. Therefore, according to the finite element theory, it can be deduced that:
\begin{align*}
& a^{i}(\boldsymbol{u}^{i}_{h},\boldsymbol{v}^{i}_{h}) = (\bar{\boldsymbol{u}}^{i})^{\top} \cdot K^{i} \cdot \bar{\boldsymbol{v}}^{i}, && \bar{\boldsymbol{u}}^{i}, \bar{\boldsymbol{v}}^{i} \in \mathbb{R}^{dN^{i}_{v}},\\
& a(\boldsymbol{u}_{h},\boldsymbol{v}_{h}) = (\bar{\boldsymbol{u}})^{\top} \cdot K \cdot \bar{\boldsymbol{v}}, && \bar{\boldsymbol{u}}, \bar{\boldsymbol{v}} \in \mathbb{R}^{dN_{v}},\\
& \int_{\Omega^{i}} \boldsymbol{f}^{i}_{0} \boldsymbol{v}^{i}_{h} dx = (\bar{\boldsymbol{f}}^{i}_{0})^{\top} \cdot \bar{\boldsymbol{v}}^{i},  && \bar{\boldsymbol{f}}^{i}_{0} \in \mathbb{R}^{dN^{i}_{v}}, \\
& \int_{\Gamma^{1}_{2}} \boldsymbol{f}_{1} \boldsymbol{v}^{1}_{h} ds = (\bar{\boldsymbol{f}}_{1})^{\top} \cdot \bar{\boldsymbol{v}}^{1}, && \bar{\boldsymbol{f}}_{1} \in \mathbb{R}^{dN^{1}_{v}},\\
& L(\boldsymbol{v}_{h}) = (\bar{\boldsymbol{f}})^{\top} \cdot \bar{\boldsymbol{v}}, && \bar{\boldsymbol{f}} \in \mathbb{R}^{dN_{v}},
\end{align*}
where $K^{i}\in \mathbb{R}^{dN^{i}_{v}\times dN^{i}_{v}}$ is a positive definite symmetric stiffness matrix, $K = diag\{ K^{1}, \cdots, K^{n}\} \in \mathbb{R}^{dN_{v}\times dN_{v}}$, $\bar{\boldsymbol{f}} = (\bar{\boldsymbol{f}}^{1}_{0} + \bar{\boldsymbol{f}}_{1}; \bar{\boldsymbol{f}}^{2}_{0}; \cdots; \bar{\boldsymbol{f}}^{n}_{0}) \in \mathbb{R}^{dN^{i}_{v}}$ and $N_{v} = \sum_{i=1}^{n} N^{i}_{v}$. Furthermore, it can be verified that:
\begin{align*}
&\left\langle \mu^{i}_{HN}, [v^{i}_{HN}] \right\rangle_{c} = \left\langle \mu^{i}_{HN}, \boldsymbol{v}^{i}_{H} \cdot \beta^{i} + \boldsymbol{v}^{i+1}_{H} \cdot \alpha^{i} \right\rangle_{c} = \int_{\Gamma^{i}_{c}} \mu^{i}_{HN} \cdot \boldsymbol{v}^{i}_{H} \cdot \beta^{i} ds + \int_{\Gamma^{i}_{c}} \mu^{i}_{HN} \cdot \boldsymbol{v}^{i+1}_{H} \cdot \alpha^{i+1} ds \\
&= \int_{\Gamma^{i}_{3}} \mu^{i}_{HN} \cdot \boldsymbol{v}^{i}_{H} \cdot \beta^{i} ds - \int_{\Gamma^{i+1}_{2}} \mu^{i}_{HN} \cdot \boldsymbol{v}^{i+1}_{H} \cdot \alpha^{i+1} ds = (\bar{\boldsymbol{\mu}}^{i}_{N})^{\top} \cdot G^{i}_{N3} \cdot \bar{\boldsymbol{v}}^{i} - (\bar{\boldsymbol{\mu}}^{i}_{N})^{\top} \cdot G^{i+1}_{N2} \cdot \bar{\boldsymbol{v}}^{i+1} \\
&= (\bar{\boldsymbol{\mu}}^{i}_{N})^{\top} \cdot \Big( \mathbf{0}, ~ G^{i}_{N3},~ -G^{i+1}_{N2},~ \mathbf{0} \Big) \cdot \bar{\boldsymbol{v}} = (\bar{\boldsymbol{\mu}}^{i}_{N})^{\top} \cdot G^{i}_{N} \cdot \bar{\boldsymbol{v}},\\
& \left\langle \boldsymbol{\mu}^{i}_{HT}, [\boldsymbol{v}^{i}_{HT}] \right\rangle_{c} = \left\langle \boldsymbol{\mu}^{i}_{HT}, \boldsymbol{v}^{i}_{H\eta} - \boldsymbol{v}^{i+1}_{H\tau} \right\rangle_{c} = \left\langle \boldsymbol{\mu}^{i}_{HT}, \boldsymbol{v}^{i}_{H} - \boldsymbol{v}^{i}_{H}\cdot\beta^{i} \cdot\beta^{i} \right\rangle_{c} - \left\langle \boldsymbol{\mu}^{i}_{HT}, \boldsymbol{v}^{i+1}_{H} - \boldsymbol{v}^{i+1}_{H} \cdot\alpha^{i+1} \cdot\alpha^{i+1} \right\rangle_{c} \\
&=\int_{\Gamma^{i}_{c}} \boldsymbol{\mu}^{i}_{HT}\cdot \boldsymbol{v}^{i}_{H\eta} ds - \int_{\Gamma^{i}_{c}} \boldsymbol{\mu}^{i}_{HT}\cdot \boldsymbol{v}^{i+1}_{H\tau} ds =\int_{\Gamma^{i}_{3}} \boldsymbol{\mu}^{i}_{HT}\cdot \boldsymbol{v}^{i}_{H\eta} ds + \int_{\Gamma^{i+1}_{2}} \boldsymbol{\mu}^{i}_{HT}\cdot \boldsymbol{v}^{i+1}_{H\tau} ds\\ 
&= (\bar{\boldsymbol{\mu}}^{i}_{T})^{\top} \cdot G^{i}_{T3} \cdot \bar{\boldsymbol{v}}^{i} + (\bar{\boldsymbol{\mu}}^{i}_{T})^{\top} \cdot G^{i+1}_{T2} \cdot \bar{\boldsymbol{v}}^{i+1} = (\bar{\boldsymbol{\mu}}^{i}_{T})^{\top} \cdot G^{i}_{T} \cdot \bar{\boldsymbol{v}},
\end{align*}
where $\bar{\boldsymbol{\mu}}^{i}_{N}\in \mathbb{R}^{N^{i}_{l}}$, 
$\bar{\boldsymbol{\mu}}^{i}_{T}\in \mathbb{R}^{dN^{i}_{l}}$, $G^{i}_{N} \in \mathbb{R}^{N^{i}_{l}\times dN_{v}}$ and $G^{i}_{T} \in \mathbb{R}^{dN^{i}_{l}\times dN_{v}}$.

\begin{remark}[]\label{4:rem:MatrixForm}
The matrices $G^{i}_{N}$ and $G^{i}_{T}$ can be constructed directly as follows:
\begin{align*}
& \left\langle \mu^{i}_{NH}, [v^{i}_{NH}] \right\rangle_{c} = \int_{\Gamma^{i}_{c}} \mu^{i}_{NH}\cdot [v^{i}_{NH}] ds = (\bar{\boldsymbol{\mu}}^{i}_{N})^{\top} \cdot \int_{\Gamma^{i}_{c}} \Psi^{i}_{N} E^{i}_{N} \bar{M}^{i}_{N} \bar{M}^{i}_{\alpha} \Phi ds \cdot \bar{\boldsymbol{v}} = (\bar{\boldsymbol{\mu}}^{i}_{N})^{\top} \cdot G^{i}_{N} \cdot \bar{\boldsymbol{v}},\\
& \left\langle \boldsymbol{\mu}^{i}_{TH}, [\boldsymbol{v}^{i}_{TH}] \right\rangle_{c} = \int_{\Gamma^{i}_{c}} \boldsymbol{\mu}^{i}_{TH}\cdot [\boldsymbol{v}^{i}_{TH}] ds = (\bar{\boldsymbol{\mu}}^{i}_{T})^{\top} \cdot \int_{\Gamma^{i}_{c}} \Psi^{i}_{T} E^{i}_{T} \bar{M}^{i}_{T} \bar{M}^{i}_{\beta} \Phi ds \cdot \bar{\boldsymbol{v}} = (\bar{\boldsymbol{\mu}}^{i}_{T})^{\top} \cdot G^{i}_{T} \cdot \bar{\boldsymbol{v}},
\end{align*}
where $\Phi = diag\{\Phi^{1}, \ldots, \Phi^{n}\} \in \mathbb{R}^{N_{v}\times N_{v}}$; 
$E^{i}_{N}\in\mathbb{R}^{N^{i}_{l}\times N^{i}_{c}}$ and $E^{i}_{T}\in\mathbb{R}^{dN^{i}_{l}\times dN^{i}_{c}}$ are matrices with all elements being $1$; 
$\bar{M}^{i}_{\alpha}\in \mathbb{R}^{N_{v}\times dN_{v}}$ is a matrix composed of the unit external normal vectors of the finite element nodes on the boundaries $\Gamma^{i}_{3}$, $\Gamma^{i+1}_{2}$ ($i=1,\ldots,n-1$) and the row vector is $0$ if the corresponding node is not on this boundaries, which can be used to calculate normal displacement; the matrix $\bar{M}^{i}_{\beta} = I_{dN_{v}} - (\bar{M}^{i}_{\alpha})^{\top} \bar{M}^{i}_{\alpha} \in \mathbb{R}^{dN_{v}\times dN_{v}}$ can be used to calculate the displacement on the tangent plane of the finite element node located at the boundaries $\Gamma^{i}_{3}$, $\Gamma^{i+1}_{2}$; 
the function of the matrix $\bar{M}^{i}_{N}\in \mathbb{R}^{N^{i}_{c}\times N_{v}}$ is to add the normal displacements of the finite element nodes on both sides of the contact zone $\Gamma^{i}_{c}$, so each row has only two $1$ elements and corresponds to the finite element node on $\Gamma^{i}_{c}$; 
the function of matrix $\bar{M}^{i}_{T}\in \mathbb{R}^{dN^{i}_{c}\times dN_{v}}$ is to subtract the tangential displacements of the finite element nodes on both sides of the contact zone $\Gamma^{i}_{c}$, so the elements in each row corresponding to the finite element nodes at this position on $\Gamma^{i}_{3}$ and $\Gamma^{i+1}_{2}$ are $1$ and $-1$, respectively. 
\end{remark}

Based on the above relation, the matrix form of $\mathcal{G}^{*}(\boldsymbol{\mu}_{H}, \boldsymbol{v}_{H})$ can be expressed as:
\begin{align*}
\mathcal{G}^{*}(\boldsymbol{\mu}_{H}, \boldsymbol{v}_{H}) = \sum_{i=1}^{n-1} \left( (\bar{\boldsymbol{\mu}}^{i}_{N})^{\top} \cdot G^{i}_{N} \cdot \bar{\boldsymbol{v}} + (\bar{\boldsymbol{\mu}}^{i}_{T})^{\top} \cdot G^{i}_{T} \cdot \bar{\boldsymbol{v}} \right) = \bar{\boldsymbol{\mu}}^{\top} \cdot G \cdot \bar{\boldsymbol{v}},
\end{align*}
where $\bar{\boldsymbol{\mu}}_{N} = \{ \bar{\boldsymbol{\mu}}^{1}_{N};\ldots; \bar{\boldsymbol{\mu}}^{n-1}_{N} \} \in \mathbb{R}^{N_{l}}$, $\bar{\boldsymbol{\mu}}_{T} = \{ \bar{\boldsymbol{\mu}}^{1}_{T};\ldots; \bar{\boldsymbol{\mu}}^{n-1}_{T} \} \in \mathbb{R}^{dN_{l}}$, $G_{N} = \{G^{1}_{N};\ldots; G^{n-1}_{N}\}\in \mathbb{R}^{N_{l}\times dN_{v}}$, $G_{T} = \{G^{1}_{T};\ldots; G^{n-1}_{T}\}\in \mathbb{R}^{dN_{l}\times dN_{v}}$, $\bar{\boldsymbol{\mu}} = \{ \bar{\boldsymbol{\mu}}_{N}; \bar{\boldsymbol{\mu}}_{T} \} \in \mathbb{R}^{(d+1)N_{l}}$, $G = \{G_{N}; G_{T}\}\in \mathbb{R}^{(d+1)N_{l}\times dN_{v}}$, $N_{l} = \sum_{i=1}^{n-1} N^{i}_{l}$. Therefore, the algebraic form of the discrete mixed finite element problem \nameref{4:prb:mixhH*} can be expressed as

\begin{problem}[$\mathcal{P}_{M}$]\label{4:prb:M}
Find a solution $\{\bar{\boldsymbol{u}}, \bar{\boldsymbol{\lambda}}\} \in \mathbb{R}^{dN_{v}} \times \Lambda_{M}$ which satisfies:
\begin{equation}\label{4:eq:var.matrix}
\left\{\begin{aligned}
& (\bar{\boldsymbol{u}})^{\top} \cdot K + \bar{\boldsymbol{\lambda}}^{\top} \cdot G = (\bar{\boldsymbol{f}})^{\top}\\
& (\bar{\boldsymbol{\mu}} - \bar{\boldsymbol{\lambda}})^{\top} \cdot G \cdot \bar{\boldsymbol{u}} \leqslant 0, ~\forall \bar{\boldsymbol{\mu}} \in \Lambda_{M} = \Lambda_{MN} \times \Lambda_{MT}
\end{aligned}\right. ,
\end{equation}
where
\begin{align*}
& \Lambda_{MN} = \left\{ \bar{\boldsymbol{\mu}}_{N} \in \mathbb{R}^{N_{l}} \Big| \Psi_{N}\cdot \bar{\boldsymbol{\mu}}_{N} \in \Lambda_{HN} \right\} ,~~ \Lambda_{MT} = \left\{ \bar{\boldsymbol{\mu}}_{T} \in \mathbb{R}^{dN_{l}} \Big| \Psi_{T}\cdot \bar{\boldsymbol{\mu}}_{T} \in \Lambda_{HT} \right\} \\
& \Psi_{N} = diag\{\Psi^{1}_{N},\ldots,\Psi^{n-1}_{N}\},~\Psi_{T} = diag\{\Psi^{1}_{T},\ldots,\Psi^{n-1}_{T}\},\\
& \Lambda_{HN}= \Lambda_{HN}^{1} \times \cdots \times \Lambda_{HN}^{n-1},~ \Lambda_{HT}= \Lambda_{HT}^{1} \times \cdots \times \Lambda_{HT}^{n-1}.
\end{align*}
\end{problem}

By evaluating $\bar{\boldsymbol{u}}$ from the equation in formula (\ref{4:eq:var.matrix}) and substituting it into the inequality in formula (\ref{4:eq:var.matrix}), the algebraic dual formulation of the mixed finite element method can be deduced as follows \cite{dostal2002implementation}:

\begin{problem}[$\mathcal{P}_{D}$]\label{4:prb:D}
Find a solution $\bar{\boldsymbol{\lambda}} \in \Lambda_{M}$ which satisfies:
\begin{equation}\label{4:eq:var.matrix.dual}
\left\{\begin{aligned}
& \bar{\boldsymbol{\lambda}} = \arg\min \frac{1}{2} \bar{\boldsymbol{\mu}}^{\top} GK^{-1}G^{\top} \bar{\boldsymbol{\mu}} - \bar{\boldsymbol{\mu}}^{\top} GK^{-1} \bar{\boldsymbol{f}}, \\
& \text{s.t.}~~ \bar{\boldsymbol{\mu}} \in \Lambda_{M}.
\end{aligned}\right. 
\end{equation}
\end{problem}
After the solution $\bar{\boldsymbol{\lambda}} \in \Lambda_{M}$ is calculated, the solution $\bar{\boldsymbol{u}}$ can be deduced using the following formula:
\begin{equation*}
\bar{\boldsymbol{u}} = K^{-1}(\bar{\boldsymbol{f}} - G^{\top} \bar{\boldsymbol{\lambda}}).
\end{equation*}

\section{Numerical examples and comparisons}

In this section, the three-dimensional three-layer elastic contact system will be used as the experimental model, and the displacement field of this model will be calculated under the action of external forces. Moreover, the layer decomposition method (LDM) and the mixed finite element method (MFEM) proposed in this paper will be used for numerical calculations, respectively, and the numerical results will be compared. These algorithms are developed based on the GetFEM finite element calculation module on Python \cite{renard2020getfem}, and the computing platform is the ISDM-MESO High-Performance Computing Center.

\subsection{Experimental model}

\begin{figure}[!t]
\centering
\includegraphics[width=4in]{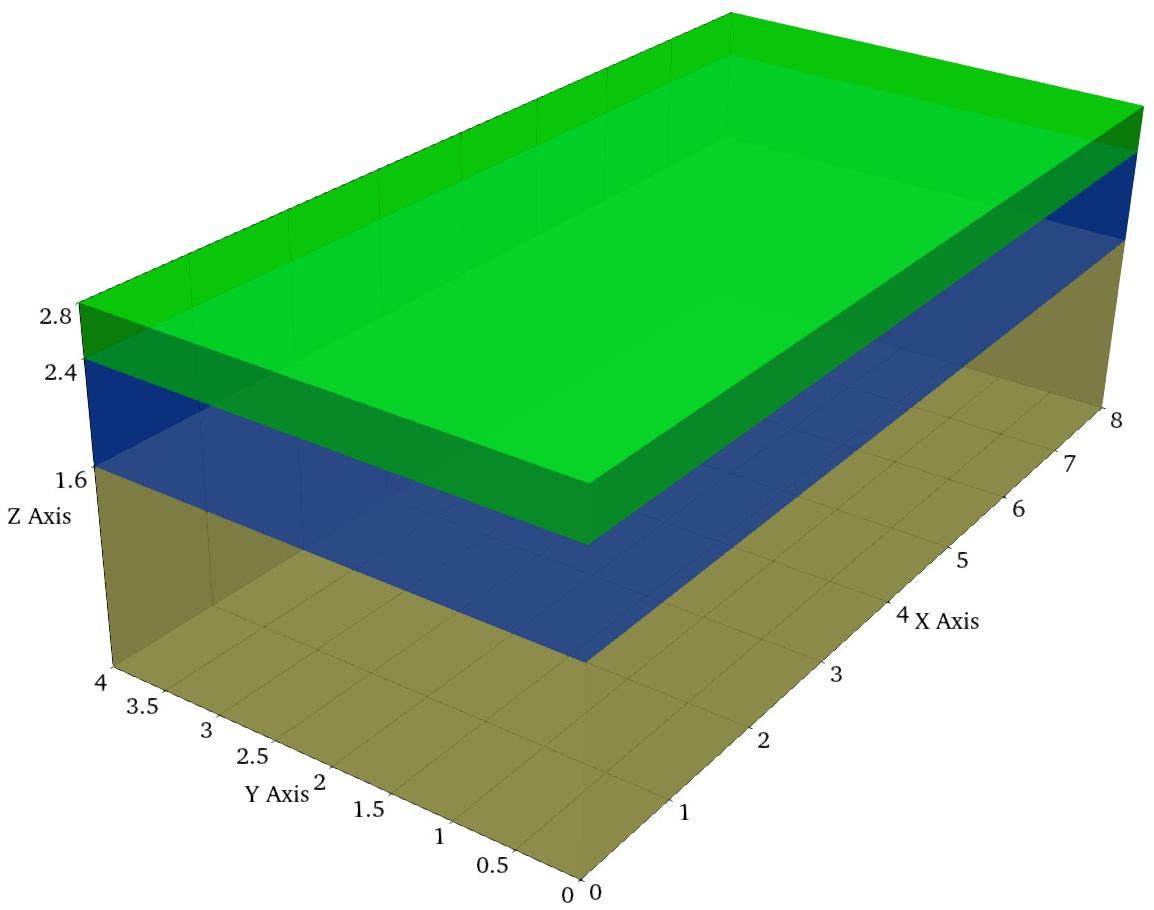}
\caption{The physical model of three-layer elastic contact system.}
\label{4:fig2}
\end{figure}

In order to simulate the mechanical response of the pavement under vehicle load, a physical model of the three-layer elastic contact system is designed, as illustrated in \textcolor{black}{Fig}.\ref{4:fig2}. The system comprises three elastic bodies, denoted as $\Omega^1$, $\Omega^2$, and $\Omega^3$, with their spatial regions defined as follows:
$$
\begin{aligned}
& \Omega^1=(0,8) \times(0,4) \times(2.4,2.8), \\
& \Omega^2=(0,8) \times(0,4) \times(1.6,2.4), \\
& \Omega^3=(0,8) \times(0,4) \times(0,1.6).
\end{aligned}
$$
In this model, all elastic bodies satisfy homogeneity and isotropy. To make the model more representative of real pavement \cite{ma2021analytical,kim2011numerical}, the elastic modulus $E^{i}$ and Poisson's ratio $P^{i}$ of the elastic body $\Omega^{i}$ are defined as follows:
$$
\begin{array}{ll}
E^1=5 \cdot 10^3, & P^1=0.25 \\
E^2=2 \cdot 10^3, & P^2=0.25 \\
E^3=2 \cdot 10^2, & P^3=0.4
\end{array}
$$

The boundaries of the model are then defined as follows:
\begin{align*}
& \Gamma_2^1=(0,8) \times(0,4) \times\{2.8\};~ \Gamma_3^1=(0,8) \times(0,4) \times\{2.4\};~ \Gamma_1^1=\partial \Omega^1 / \overline{\Gamma_2^1 \cup \Gamma_3^1}; \\
& \Gamma_2^2=(0,8) \times(0,4) \times\{2.4\};~ \Gamma_3^2=(0,8) \times(0,4) \times\{1.6\};~ \Gamma_1^2=\partial \Omega^2 / \overline{\Gamma_2^2 \cup \Gamma_3^2}; \\
& \Gamma_2^3=(0,8) \times(0,4) \times\{1.6\};~ \Gamma_1^3=\partial \Omega^3 / \overline{\Gamma_2^3} ;
\end{align*}
And the friction functions on the contact zones $\Gamma_c^1$ and $\Gamma_c^2$ are respectively defined as:
$$
g^1(x)=0.2 \text { on } \Gamma_c^1,~ g^2(x)=0.05 \text { on } \Gamma_c^2 .
$$
The three-layer system will then be subjected to a gravitational force $\boldsymbol{f}_0 = [0,0,-0.05]^{\top}$ and a surface force $\boldsymbol{f}_1 = [0,-4.5,-22.5]^{\top}$, applied specifically on the boundary $\Gamma^{1}_{2}$. The surface force is exerted over the rectangular area defined by the coordinates $(3.8,4.4)\times(1.8,2.2)$.

\subsection{Mesh and Finite Element Space}

\begin{table}[t]
\renewcommand{\arraystretch}{2}
\centering
\begin{tabular}{p{3cm}|p{3cm}|p{3cm}|p{3cm}}
$H$  & $0.2$  & $0.16$  & $0.1$ \\ \hline
$\|\bar{\boldsymbol{u}}^{1}_{0}-\bar{\boldsymbol{u}}^{1}_{1}\|$ & $4.30*10^{-5}$ & $4.13*10^{-5}$ & $3.25*10^{-5}$\\
$\|\bar{\boldsymbol{u}}^{2}_{0}-\bar{\boldsymbol{u}}^{2}_{1}\|$ & $7.08*10^{-5}$ & $6.64*10^{-5}$ & $5.26*10^{-5}$\\
$\|\bar{\boldsymbol{u}}^{3}_{0}-\bar{\boldsymbol{u}}^{3}_{1}\|$ & $8.45*10^{-5}$ & $6.07*10^{-5}$ & $4.26*10^{-5}$\\
$\|\bar{\boldsymbol{u}}^{1}_{0}-\bar{\boldsymbol{u}}^{1}_{1}\|/ \|\bar{\boldsymbol{u}}^{1}_{0}\|$ & $6.89*10^{-4}$ & $5.22*10^{-4}$ & $2.47*10^{-4}$\\
$\|\bar{\boldsymbol{u}}^{2}_{0}-\bar{\boldsymbol{u}}^{2}_{1}\|/ \|\bar{\boldsymbol{u}}^{2}_{0}\|$ & $9.11*10^{-4}$ & $6.74*10^{-4}$ & $3.22*10^{-4}$\\
$\|\bar{\boldsymbol{u}}^{3}_{0}-\bar{\boldsymbol{u}}^{3}_{1}\|/ \|\bar{\boldsymbol{u}}^{3}_{0}\|$ & $1.383*10^{-3}$ & $7.84*10^{-4}$ & $3.32*10^{-4}$\\
\end{tabular}
\caption{Error between displacement vector $\bar{\boldsymbol{u}}^{i}_{0}$ and $\bar{\boldsymbol{u}}^{i}_{1}$}
\label{4:tab:1}
\end{table}

In this numerical experiment, \textcolor{black}{the elastic bodies} $\Omega^1$, $\Omega^2$ and $\Omega^3$ are all meshed using standard tetrahedral elements. Specifically, in this section, the displacement field of this physical model is calculated respectively when $H=0.2$, $H=0.16$ and $H=0.1$. Then, as defined in formula (\ref{4:space:Vih.P1}), the first-order piecewise linear polynomial space serves as a finite element space. At this time, the contact element space can be defined as $\Lambda^{*}_{H}$ or $\Lambda^{**}_{H}$, and the specific definition can be found in formulas (\ref{4:space:LiNH.P0})-(\ref{4:space:LiTH.P1}). 
It is worth explaining that the $\Lambda^{i*}_{H}$ space and the gradient of the displacement field are homogeneous. However, the degrees of freedom in the $\Lambda^{**}_{H}$ space are smaller, which is advantageous for calculations. To compare the errors in numerical solution results caused by these two types of boundary contact element spaces, the absolute errors $\|\bar{\boldsymbol{u}}^{i}_{0}-\bar{\boldsymbol{u}}^{i}_{1}\|$ and relative errors $\|\bar{\boldsymbol{u}}^{i}_{0}-\bar{\boldsymbol{u}}^{i}_{1}\|/ \|\bar{\boldsymbol{u}}^{i}_{0}\|$ of the displacement vectors of the finite element nodes under different mesh divisions were calculated separately. The results are shown in Table \ref{4:tab:1}, where $\bar{\boldsymbol{u}}^{i}_{0}$ and $\bar{\boldsymbol{u}}^{i}_{1}$ are the displacement vectors obtained using the $\Lambda^{*}_{H}$ and $\Lambda^{**}_{H}$ contact element spaces, respectively. Upon comparison, it is evident that the error decreases with the subdivision of the finite element mesh. Therefore, the numerical results under $H=0.1$ and the $\Lambda^{**}_{H}$ space will be utilized in subsequent comparative experiments.

\begin{figure}[ht!]
\centering
\begin{minipage}{0.32\linewidth}
    \centering
    \includegraphics[width=0.9\linewidth]{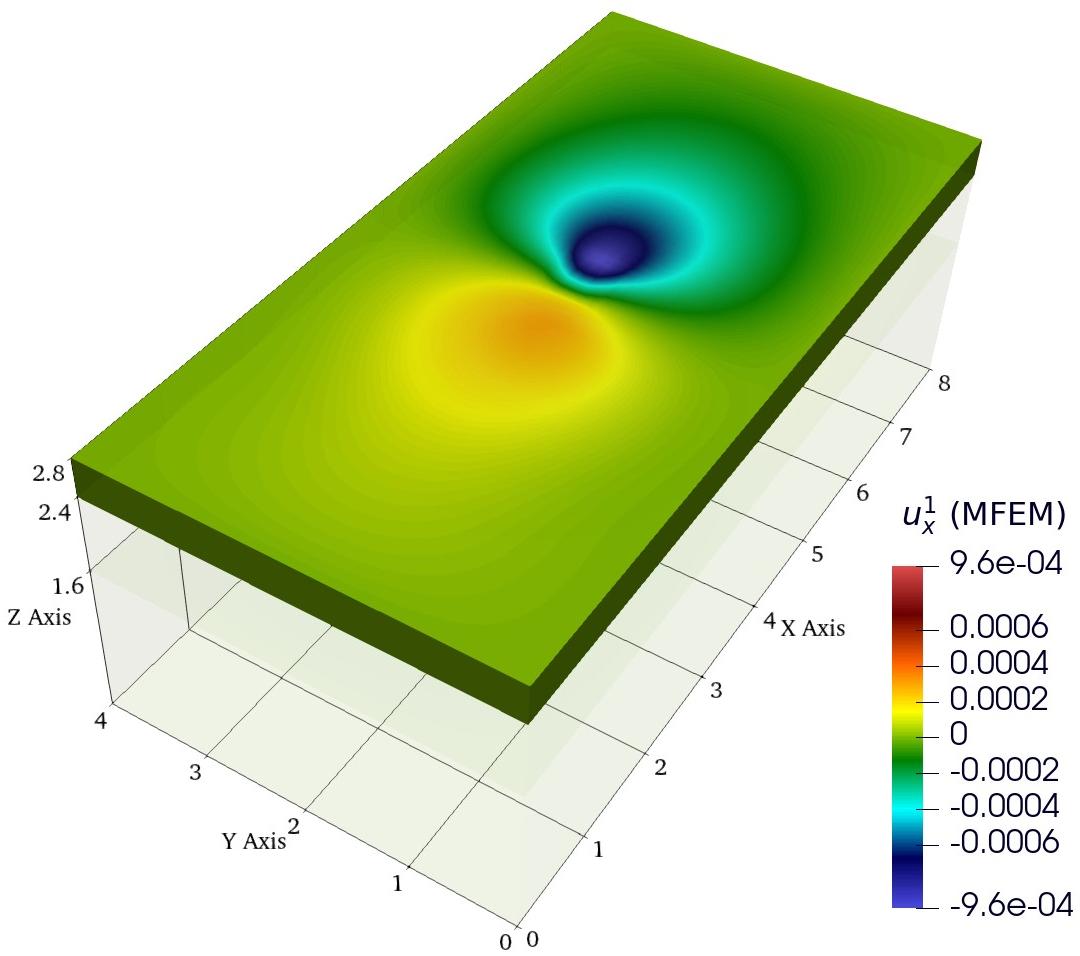}
    \caption{Displacement nephogram along the X-axis of the $\Omega^1$ obtained by MFEM.}
    \label{4:MFEM1x}
\end{minipage}
\begin{minipage}{0.32\linewidth}
    \centering
    \includegraphics[width=0.9\linewidth]{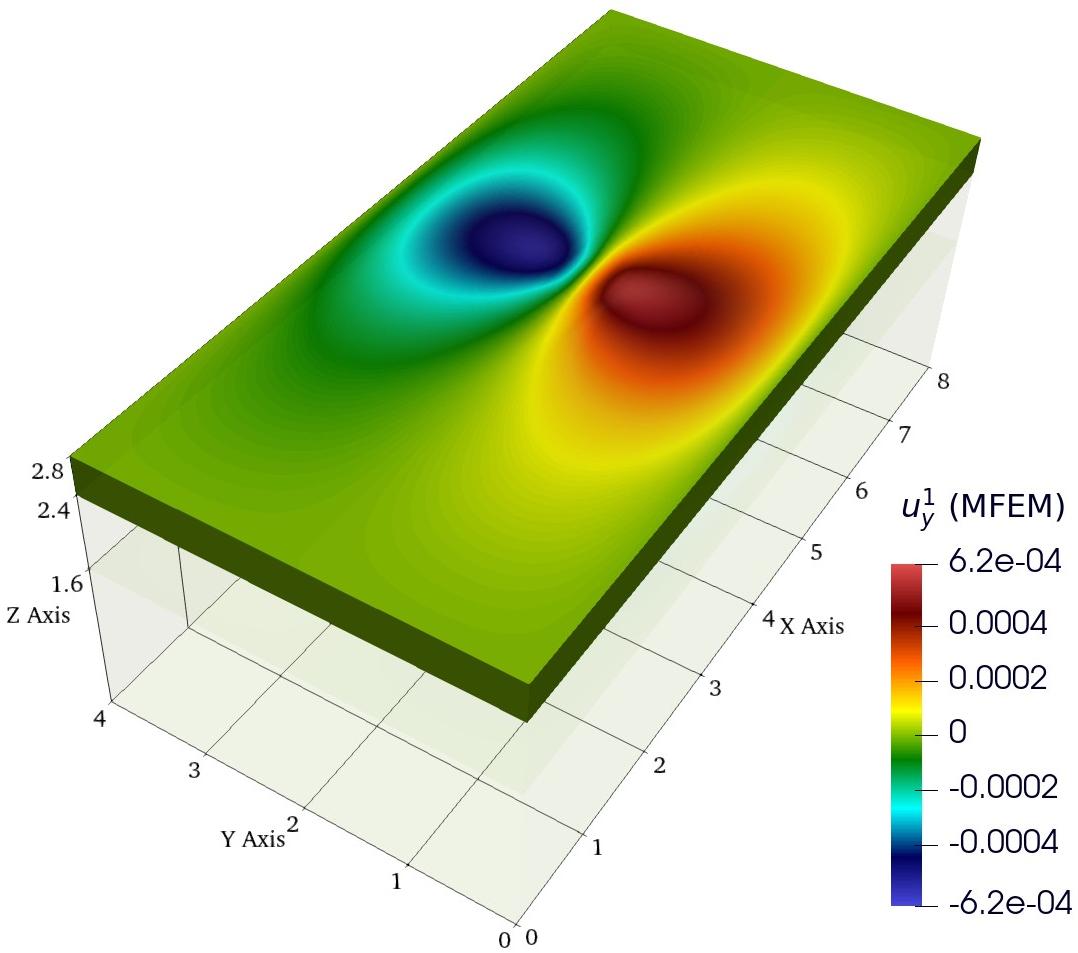}
    \caption{Displacement nephogram along the Y-axis of the $\Omega^1$ obtained by MFEM.}
    \label{4:MFEM1y}
\end{minipage}
\begin{minipage}{0.32\linewidth}
    \centering
    \includegraphics[width=0.9\linewidth]{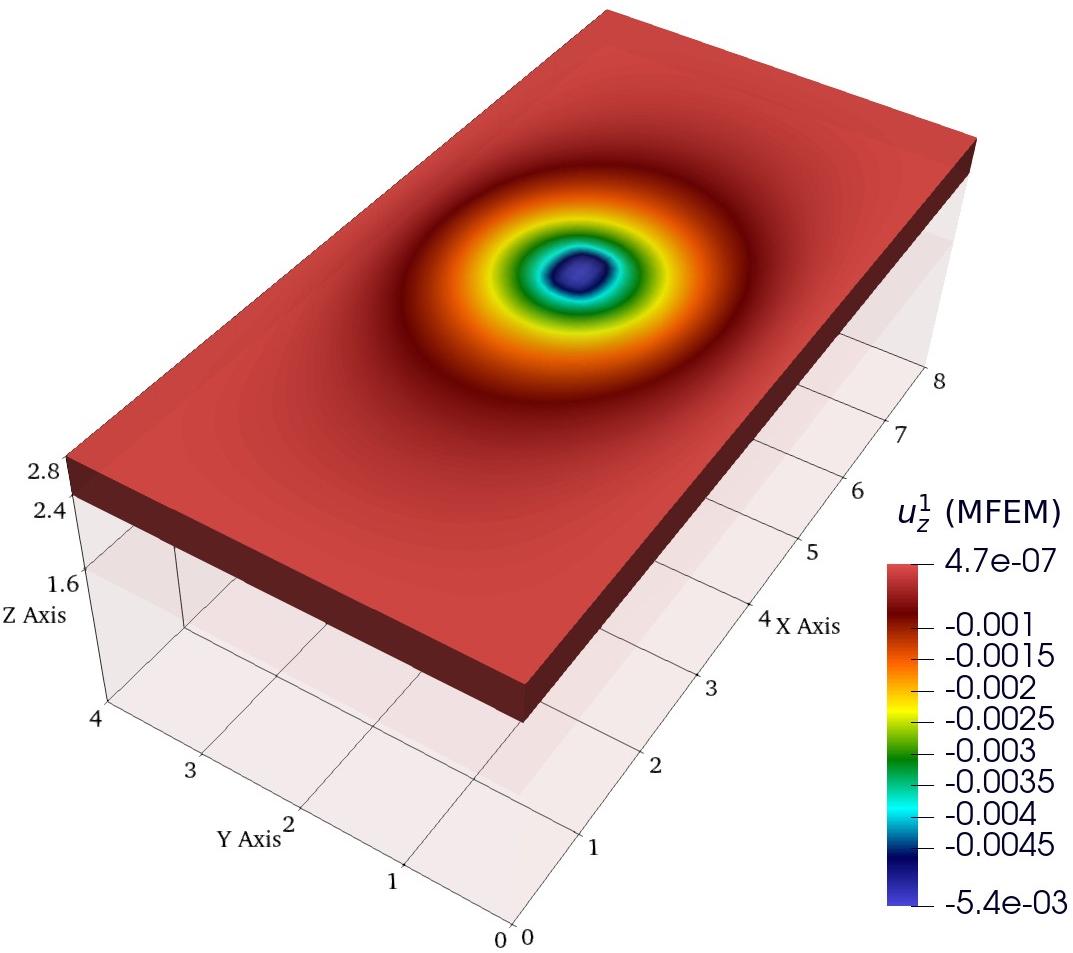}
    \caption{Displacement nephogram along the Z-axis of the $\Omega^1$ obtained by MFEM.}
    \label{4:MFEM1z}
\end{minipage}

\begin{minipage}{0.32\linewidth}
    \centering
    \includegraphics[width=0.9\linewidth]{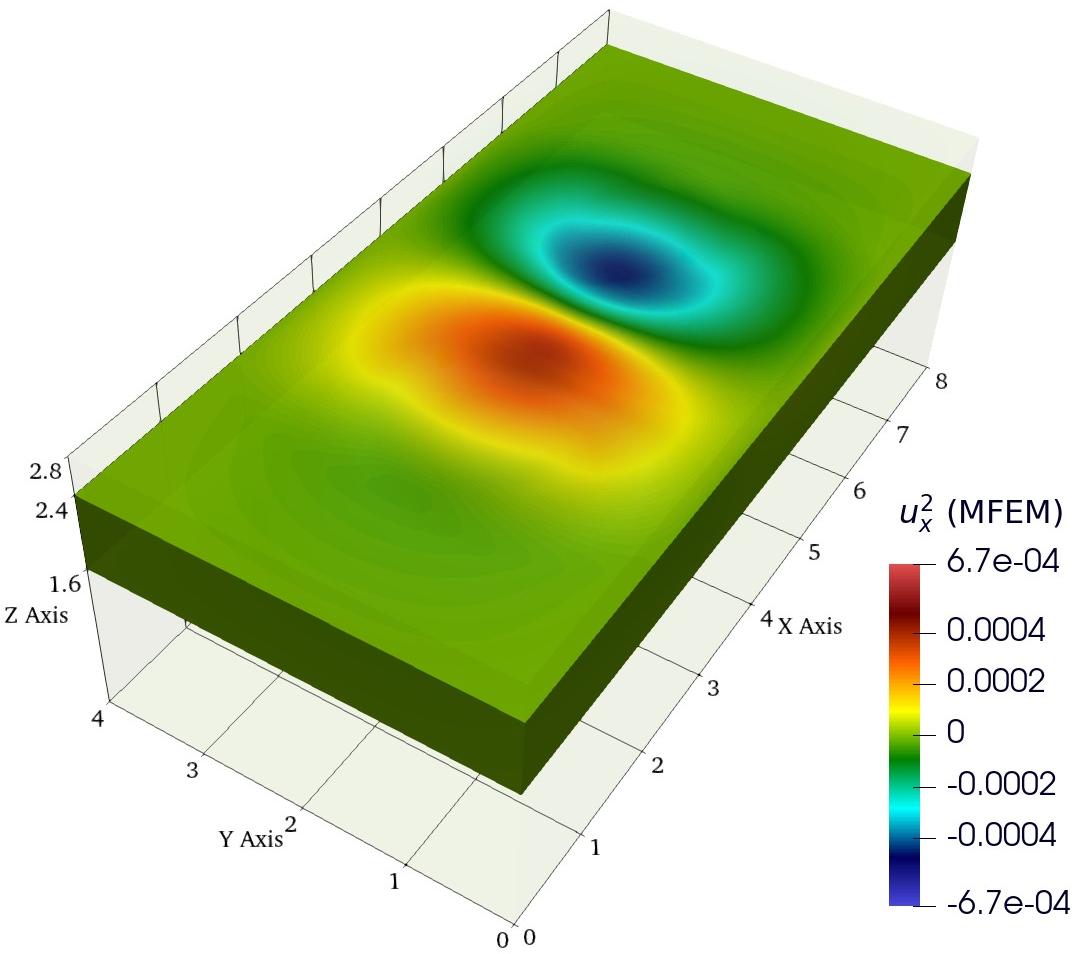}
    \caption{Displacement nephogram along the X-axis of the $\Omega^2$ obtained by MFEM.}
    \label{4:MFEM2x}
\end{minipage}
\begin{minipage}{0.32\linewidth}
    \centering
    \includegraphics[width=0.9\linewidth]{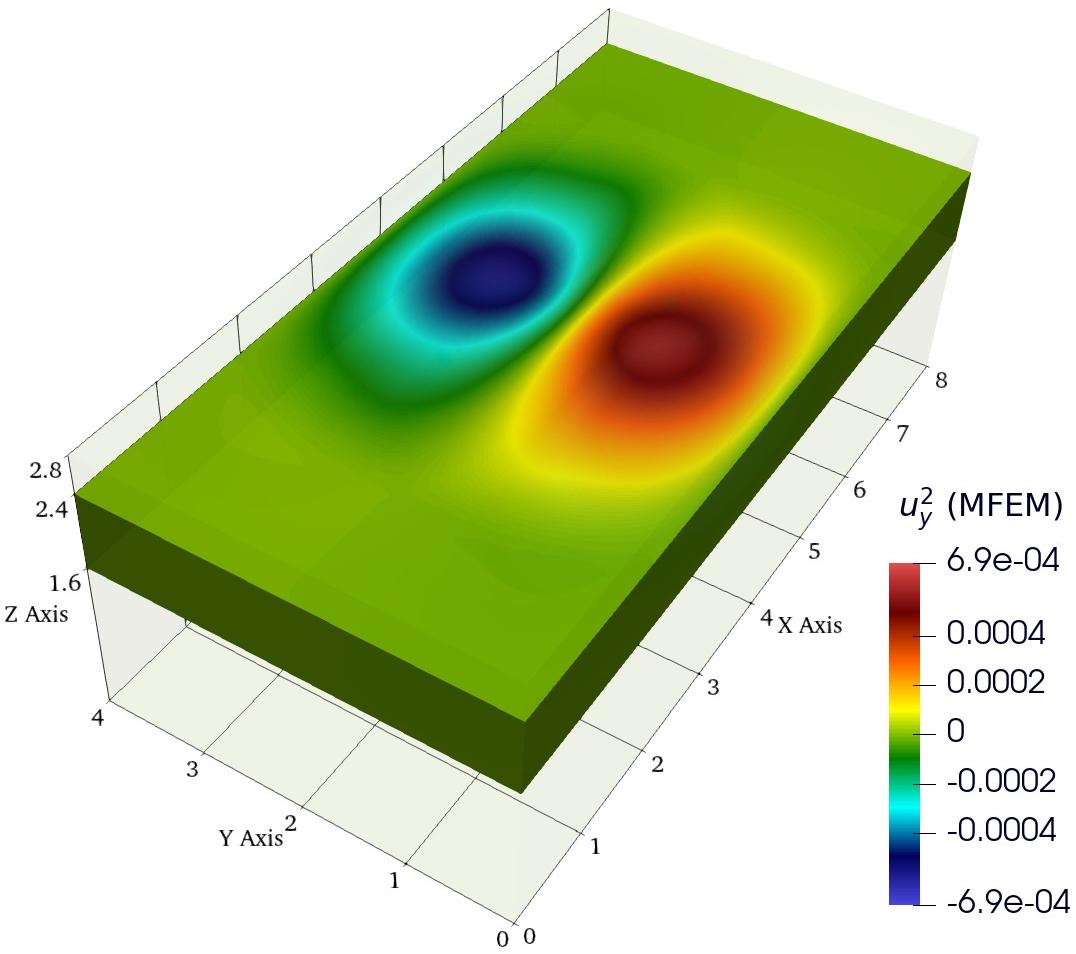}
    \caption{Displacement nephogram along the Y-axis of the $\Omega^2$ obtained by MFEM.}
    \label{4:MFEM2y}
\end{minipage}
\begin{minipage}{0.32\linewidth}
    \centering
    \includegraphics[width=0.9\linewidth]{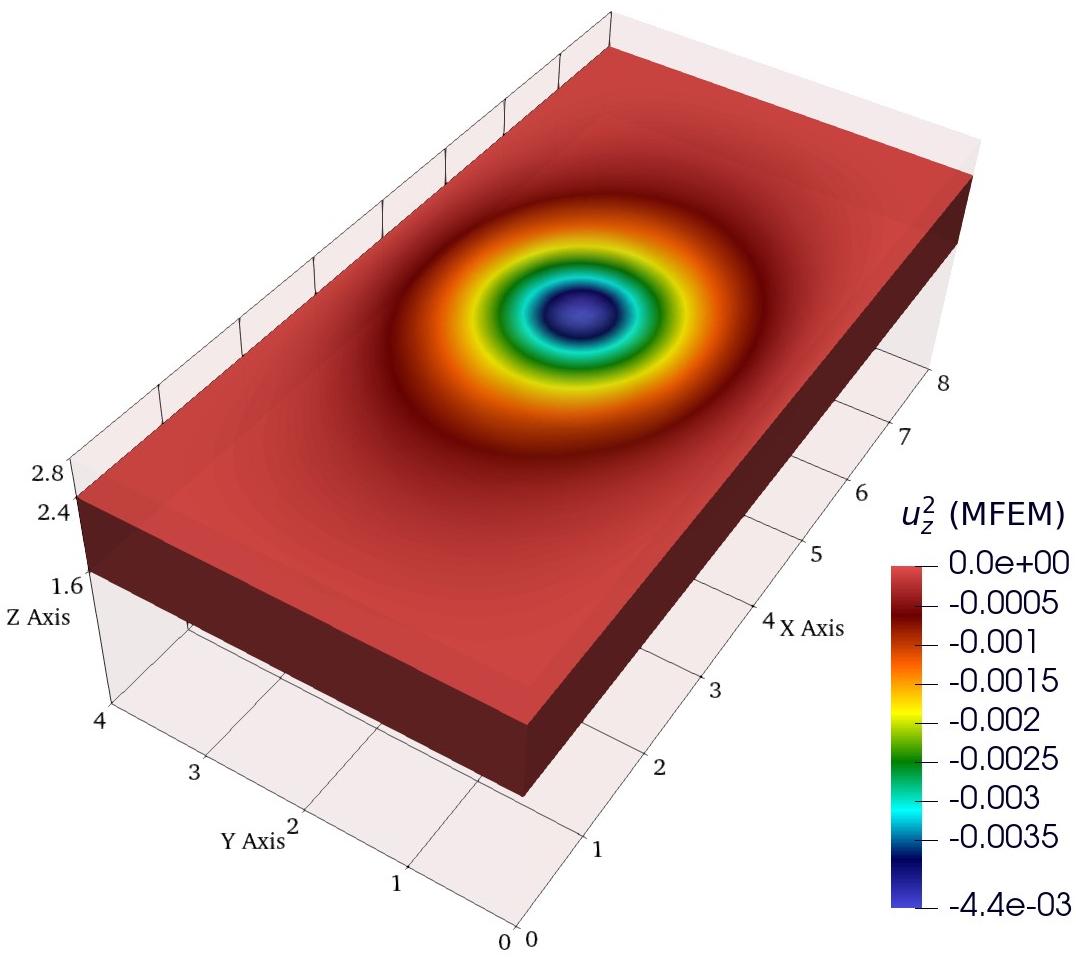}
    \caption{Displacement nephogram along the Z-axis of the $\Omega^2$ obtained by MFEM.}
    \label{4:MFEM2z}
\end{minipage}

\begin{minipage}{0.32\linewidth}
    \centering
    \includegraphics[width=0.9\linewidth]{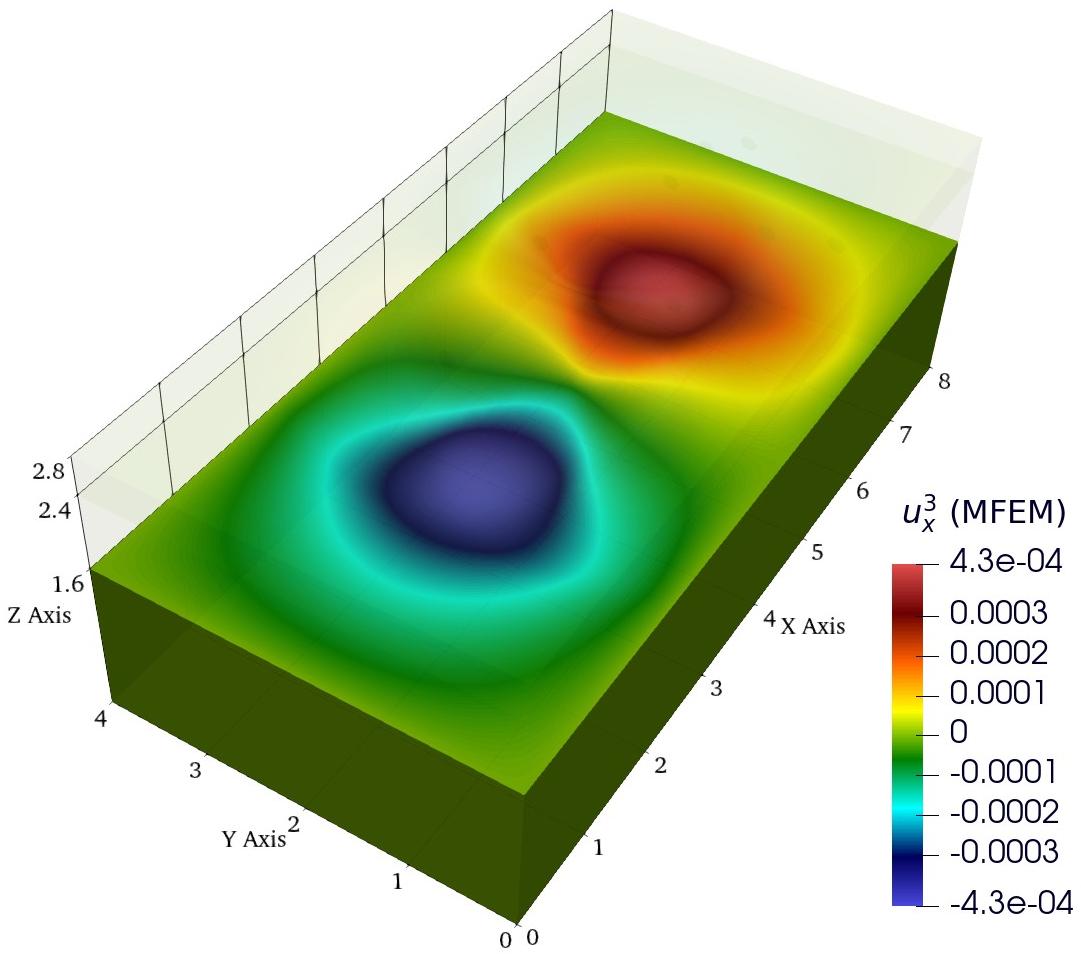}
    \caption{Displacement nephogram along the X-axis of the $\Omega^3$ obtained by MFEM.}
    \label{4:MFEM3x}
\end{minipage}
\begin{minipage}{0.32\linewidth}
    \centering
    \includegraphics[width=0.9\linewidth]{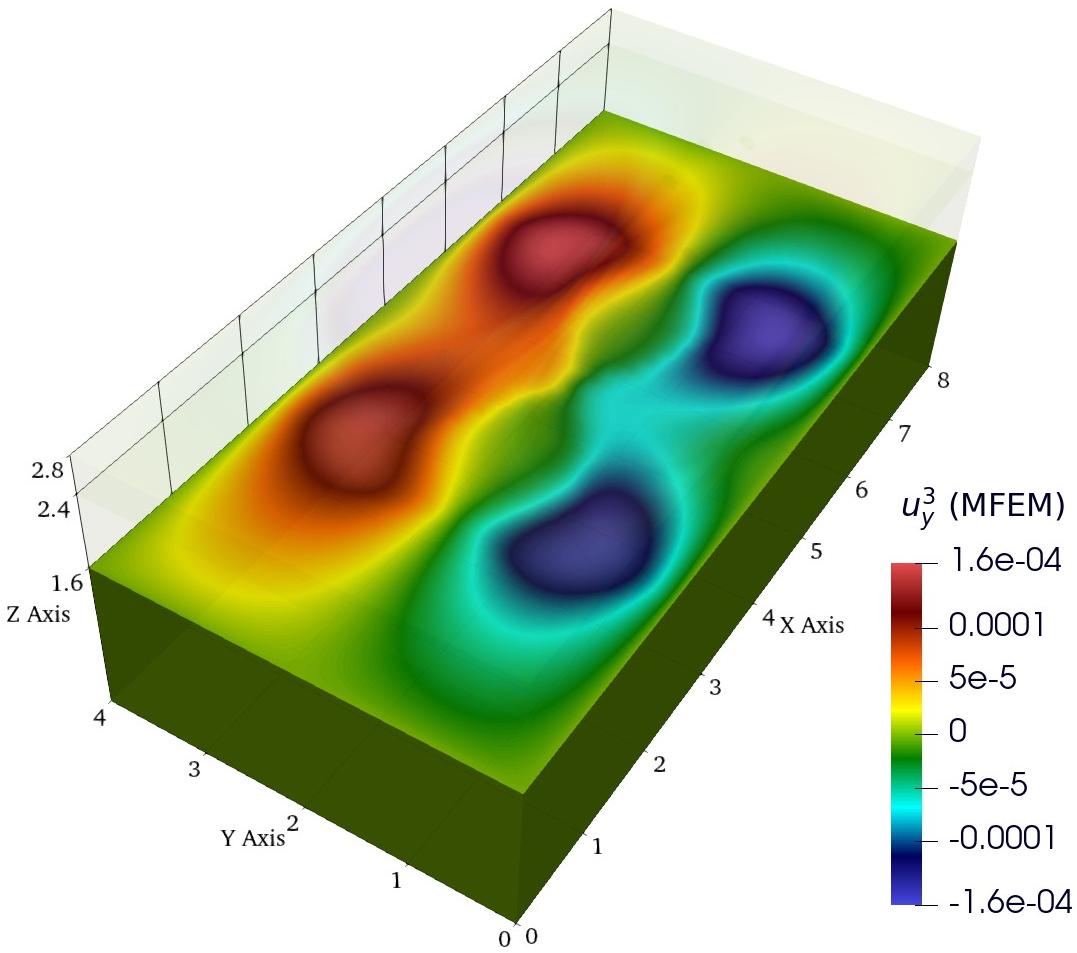}
    \caption{Displacement nephogram along the Y-axis of the $\Omega^3$ obtained by MFEM.}
    \label{4:MFEM3y}
\end{minipage}
\begin{minipage}{0.32\linewidth}
    \centering
    \includegraphics[width=0.9\linewidth]{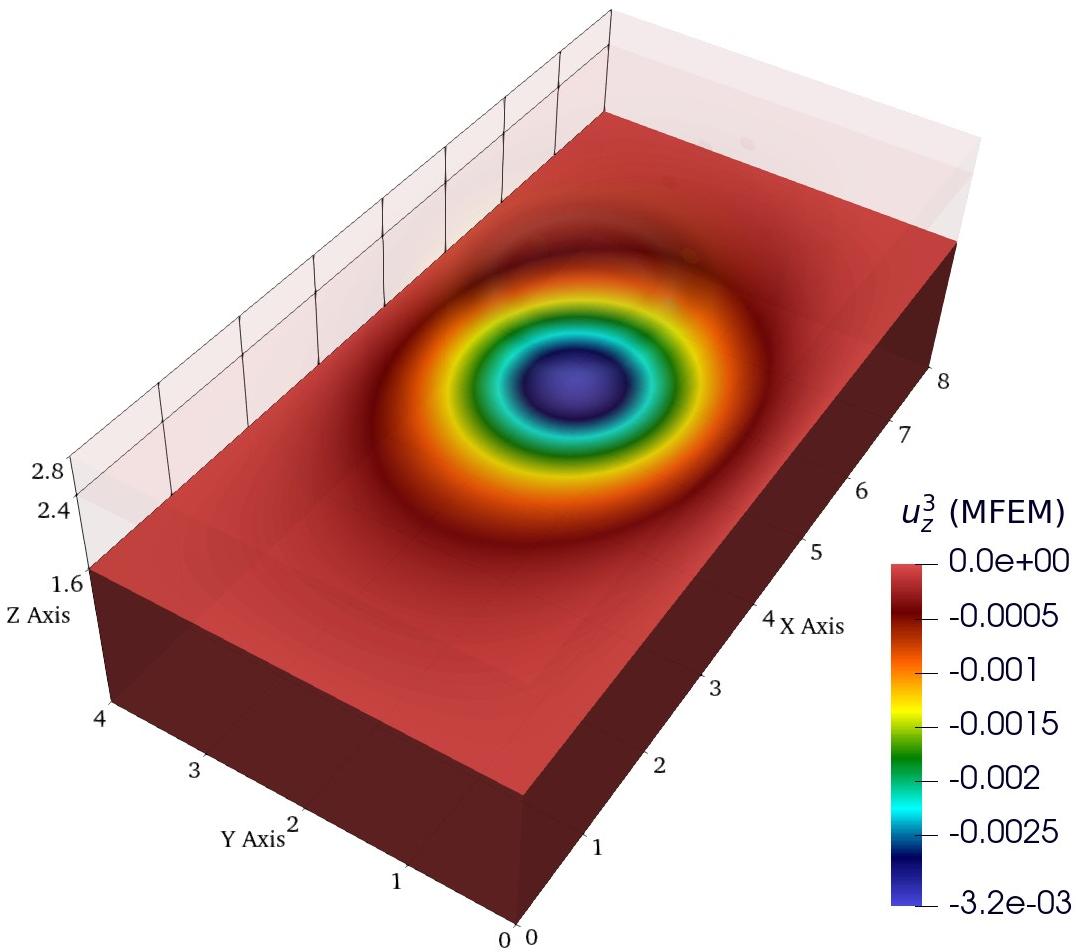}
    \caption{Displacement nephogram along the Z-axis of the $\Omega^3$ obtained by MFEM.}
    \label{4:MFEM3z}
\end{minipage}
\end{figure}

\begin{figure}[!ht]
\centering
\begin{minipage}{0.45\linewidth}
    \centering
    \includegraphics[width=0.9\linewidth]{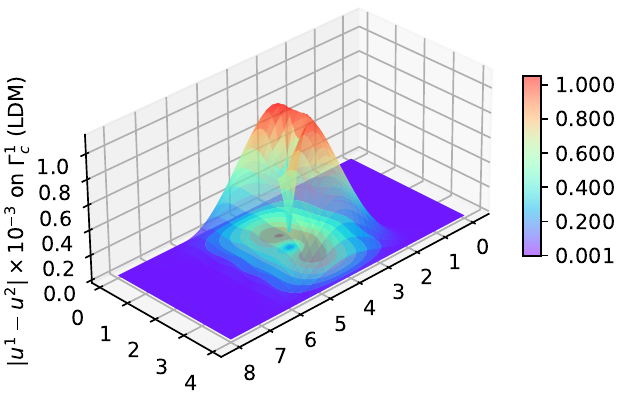}
    \caption{The difference between the displacement of $\Omega^1$ and $\Omega^2$ in the contact zone $\Gamma^1_{c}$, which calculated by LDM.}
    \label{4:LDM_1322}
\end{minipage}
\centering
\begin{minipage}{0.45\linewidth}
    \centering
    \includegraphics[width=0.9\linewidth]{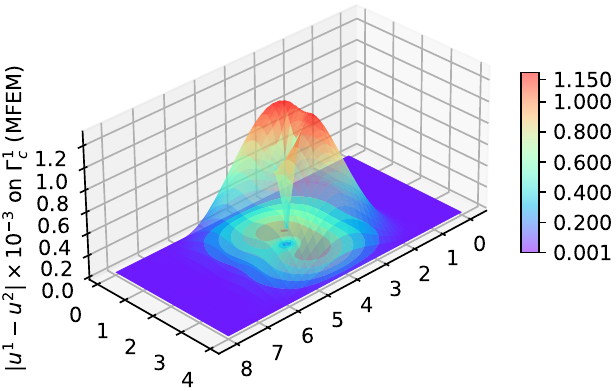}
    \caption{The difference between the displacement of $\Omega^1$ and $\Omega^2$ in the contact zone $\Gamma^1_{c}$, which calculated by MFEM.}
    \label{4:MFEM_1322}
\end{minipage}

\begin{minipage}{0.45\linewidth}
    \centering
    \includegraphics[width=0.9\linewidth]{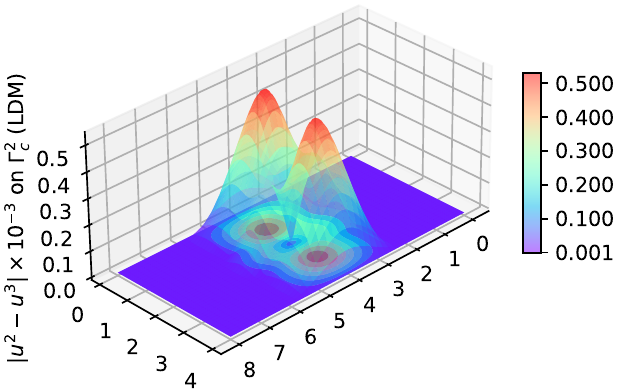}
    \caption{The difference between the displacement of $\Omega^2$ and $\Omega^3$ in the contact zone $\Gamma^2_{c}$, which calculated by LDM.}
    \label{4:LDM_2332}
\end{minipage}
\begin{minipage}{0.45\linewidth}
    \centering
    \includegraphics[width=0.9\linewidth]{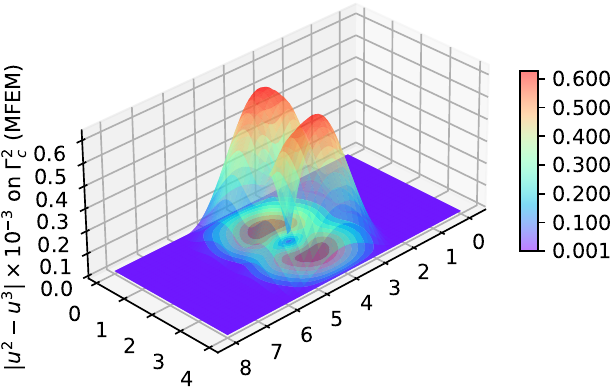}
    \caption{The difference between the displacement of $\Omega^2$ and $\Omega^3$ in the contact zone $\Gamma^2_{c}$, which calculated by MFEM.}
    \label{4:MFEM_2332}
\end{minipage}
\end{figure}

\subsection{Comparison between LDM and MFEM}

In Fig.\ref{4:MFEM1x}-Fig.\ref{4:MFEM3z}, the displacement fields $\boldsymbol{u}^{i}=(u^{i}_x,u^{i}_y,u^{i}_z)$ $(i=1,2,3)$ of the three-layer system obtained by MFEM are shown. It can be observed that the deformation characteristics of the elastic bodies in each layer are significantly different, which is attributed to the nonlinear contact conditions between the layers.
Furthermore, to compare the slip states calculated by the MFEM and LDM on the contact zones $\Gamma_{c}^{1}$ and $\Gamma_{c}^{2}$ for consistency, the displacement differences $|\boldsymbol{u}^{1}-\boldsymbol{u}^{2}|$ and $|\boldsymbol{u}^{2}-\boldsymbol{u}^{3}|$ are illustrated in Fig.\ref{4:LDM_1322}-Fig.\ref{4:MFEM_2332}. After preliminary comparison, it is confirmed that in the contact zone $\Gamma_{c}^{1}$, the slip area and bonding area obtained by the two algorithms are essentially consistent. However, there is a slight difference in the results obtained by the two algorithms in the contact zone $\Gamma_{c}^{2}$. Due to the limitations of mesh division and finite element space, both algorithms can only approximate the exact solution. Therefore, it is expected that there will be minor differences in the results obtained.

Finally, from the perspective of algorithm implementation, the two algorithms are suitable for different computing scenarios. First of all, MFEM can also be called a global method, which only needs to deal with one optimization problem during calculation. However, the matrix dimension in this optimization problem is determined by the finite element degrees of freedom of the entire system, so the calculation efficiency is high and the computational memory requirement is greater. On the other hand, the LDM algorithm decomposes the variational inequality problem of the entire system into multiple sub-optimization problems according to the number of layers. The matrix dimensions in each sub-problem are determined by the finite element degrees of freedom of a single elastic body, so these subproblems have smaller computational memory requirements. Moreover, the iterative approximation solution strategy requires repeated operations on these sub-problems, which allows the algorithm to handle larger-scale contact problems while being less computationally efficient than MFEM.

\subsection{Numerical validation of the error estimation}

\begin{figure}[!t]
\centering
\includegraphics[width=0.7\linewidth]{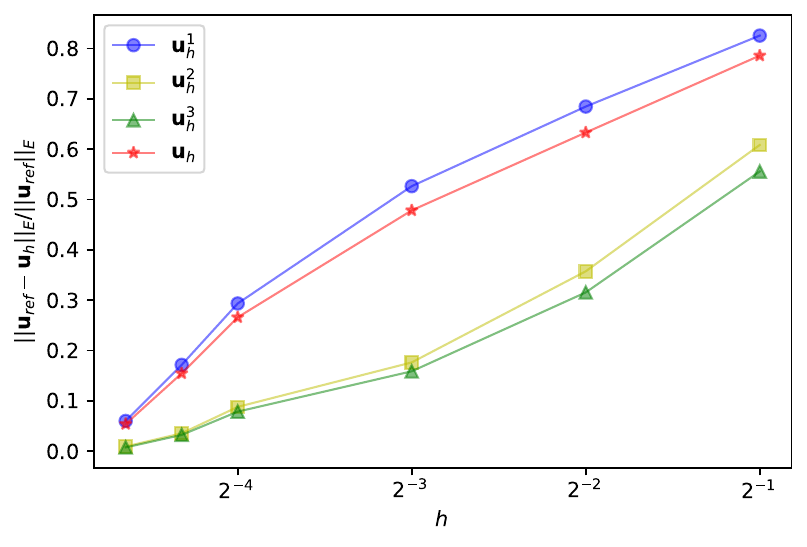}
\caption{Relative Numerical errors: $\|\boldsymbol{u}_{h}^{i}-\boldsymbol{u}_{r}^{i}\|_{E}/\|\boldsymbol{u}_{r}^{i}\|_{E}$ and $\|\boldsymbol{u}_{h}-\boldsymbol{u}_{r}\|_{E}/\|\boldsymbol{u}_{r}\|_{E}$}.
\label{4:Num_err}
\end{figure}

\begin{figure}[t]
\centering
\begin{minipage}{0.45\linewidth}
    \centering
    \includegraphics[width=0.9\linewidth]{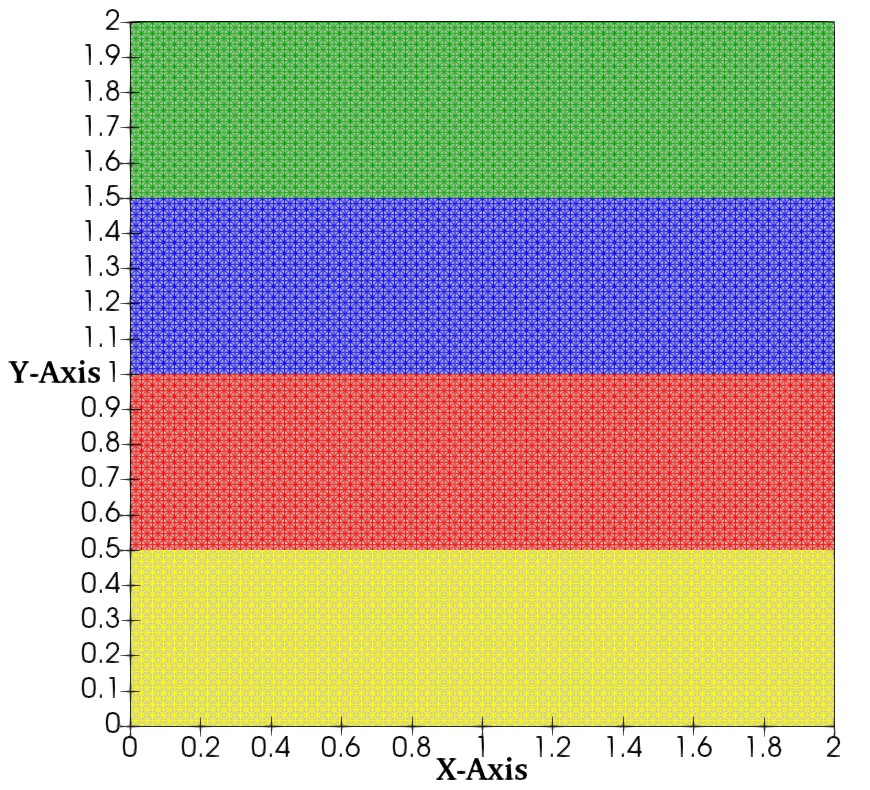}
    \caption{The physical model of two-dimensional four-layer contact system.}
    \label{4:fig:Model_2D_4L}
\end{minipage}
\centering
\begin{minipage}{0.47\linewidth}
    \centering
    \includegraphics[width=0.9\linewidth]{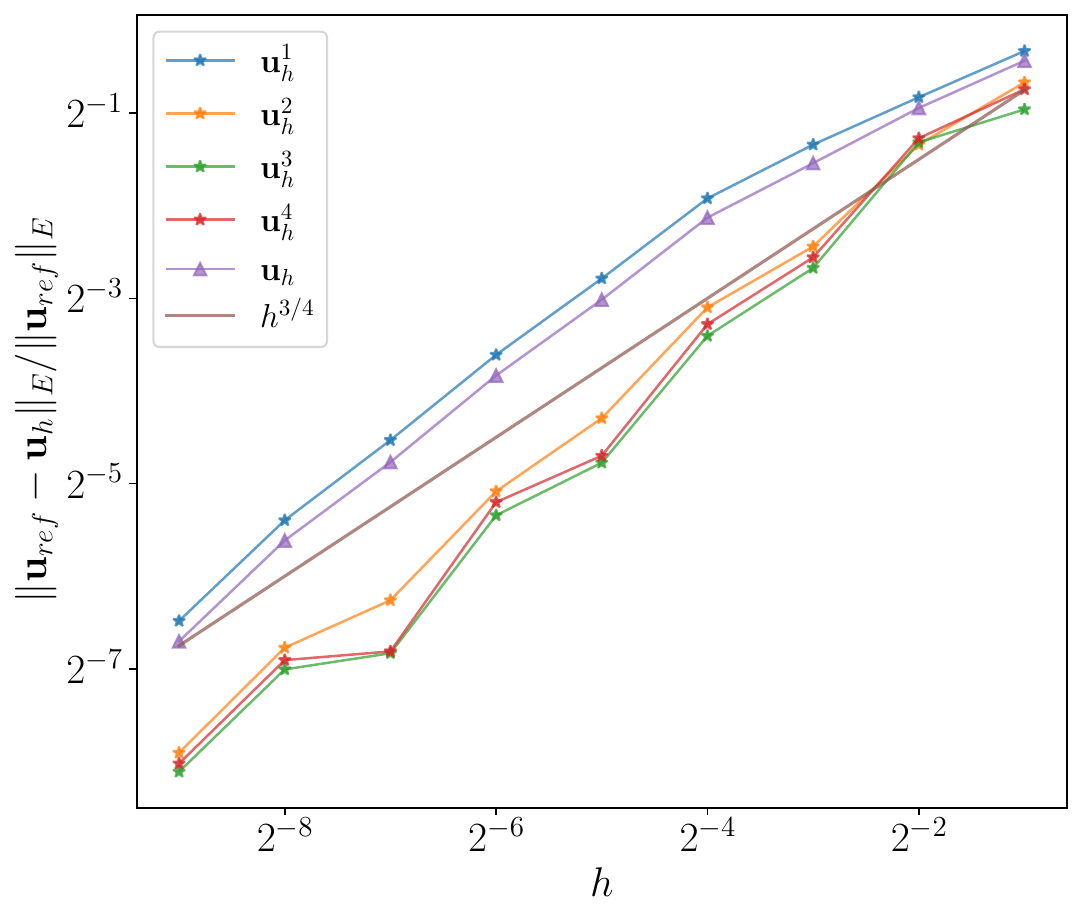}
    \caption{Relative numerical error in the two-dimensional experiment}
    \label{4:figure:Error_2D}
\end{minipage}
\end{figure}

In order to verify the convergence of the mixed finite element method, a sequence of numerical solutions of a new model with respect to the spatial discretization parameter $H=h$ is calculated. For this purpose, the energy norm $\|\cdot\|_{E}$ is defined as:
$$
\left\|\boldsymbol{v}_{h}^{i}\right\|_E=\frac{1}{\sqrt{2}}\left(a^{i}\left(\boldsymbol{v}_{h}^{i},\boldsymbol{v}_{h}^{i}\right)\right)^{1/2},~~
\left\|\boldsymbol{v}_{h}\right\|_E=\frac{1}{\sqrt{2}}\left(a\left(\boldsymbol{v}_{h},\boldsymbol{v}_{h}\right)\right)^{1/2}.
$$
Therefore, the absolute error of the displacement field of the $i$-th layer of elastic body and the error of the total displacement field in this model are recorded as $\|\boldsymbol{u}_{h}^{i}-\boldsymbol{u}^{i}\|_{E}$ ($i=1,2,3$) and $\|\boldsymbol{u}_{h}-\boldsymbol{u}\|_{E}$ respectively. 
However, since the exact solution $\boldsymbol{u}$ of this model cannot be calculated analytically, a reference numerical $\boldsymbol{u}_{ref}$ solution needs to be calculated and approximately replaces the exact solution. We start with $H=1/2$ which are successively divided by $2$.
Here, the numerical solution $\boldsymbol{u}_{h}$ corresponding to $H=1/32$ is used as the reference solution $\boldsymbol{u}_{ref}$, which corresponds to a problem with $1282905$ degrees of freedom and $2359296$ finite elements. The simulation of reference solution runs in $30409.87$ (expressed in seconds) CPU time and the relationship between the relative error $\|\boldsymbol{u}_{h}-\boldsymbol{u}_{ref}\|_{E}/\|\boldsymbol{u}_{ref}\|_{E}$ and the spatial discretization parameter $H$ is shown in Fig.\ref{4:Num_err}. 
\textcolor{black}{According to the error result, it can be inferred that the finite element numerical solution $\boldsymbol{u}_{h}$ will converge as the parameter $H\to 0$.
To more accurately verify whether the convergence order of the numerical solution aligns with the conclusion of Theorem \ref{4:thm:error.estimate}, a more precise reference solution $\boldsymbol{u}_{ref}$ needs to be computed. To this end, a numerical experiment was conducted on the two-dimensional four-layer contact system illustrated in Fig. \ref{4:fig:Model_2D_4L}, with the parameters set as follows:}
\begin{align*}
& E^1= 1 \cdot 10^4, ~ E^2= 1.5 \cdot 10^4, ~ E^3= 8 \cdot 10^3, ~ E^4= 6.5 \cdot 10^3,\\
& P^1= 0.3, ~ P^2=0.3, ~ P^3=0.3, ~ P^4=0.35, \\
& g^{1}(x) = 20 \text{ on } \Gamma^{1}_c,~ g^{2}(x) = 15 \text{ on } \Gamma^{2}_c,~ g^{3}(x) = 5 \text{ on } \Gamma^{3}_c,\\
&\boldsymbol{f}_{0} = [0,0]^{\top},~ \boldsymbol{f}_{1} = [0,-450]^{\top} \text{ on } (0.9,1.1) \times \{2\}.
\end{align*}
In this experiment, the numerical solution with $H=1/2^{10}$ was used as the reference solution $\boldsymbol{u}_{ref}$, which had a CPU computation time of $8672.63$ seconds. 
The final error results are shown in Fig.\ref{4:figure:Error_2D}, indicating that the convergence rate is approximately $h^{3/4}$, which is consistent with the conclusion of Theorem \ref{4:thm:error.estimate}.

\section{Conclusion}

In summary, the above study explores and discusses the mixed finite element method for multi-layer elastic contact systems.
First, in Section 2, the physical model of this multi-layer elastic contact system is introduced. Subsequently, the corresponding variational inequalities, equivalent optimization problems, and saddle point problems are discussed and proven, respectively.
In Section 3, the convergence of numerical solutions to discrete mixed problems in finite element space is analyzed, and the corresponding error estimates are explored.
In order to ensure the deployment of this mixed method on computers, the algebraic form of the algorithm is derived in Section 4 for specific finite element space and boundary contact element space.
Finally, based on the algebraic form of the algorithm, numerical experiments on the displacement field of the three-layer elastic contact system are carried out in Section 5.
In particular, the mixed finite element method is compared with our previously proposed layer decomposition method, and the advantages and disadvantages of the two methods are also discussed.





\section*{Acknowledgements}
This work was supported by the National Key Research and Development Program of China under Grant No.2020YFA0714300, the National Natural Science Foundations of China under Grant No.12326311 and No.61673111.

\bibliographystyle{elsarticle-num}
\bibliography{reference.bib}


\end{document}